\documentclass{article}

\usepackage[utf8]{inputenc}
\usepackage{amsmath}
\usepackage{amsfonts}
\usepackage{enumerate}
\usepackage{subfigure}
\usepackage{textcomp}
\usepackage{tabularx}
\usepackage{numprint}
\usepackage{cancel}
\usepackage{booktabs}
\usepackage{multicol}
\usepackage{tikz}
\usepackage{hyperref}
\usepackage{verbatim}
\usepackage{rotating}
\usetikzlibrary{arrows,arrows.meta,decorations.markings}
\usetikzlibrary{calc,bending}
\usepackage{pdflscape}
 \usepackage[pagewise]{lineno}
 \usepackage[normalem]{ulem}
 \usepackage{multirow}

 \usepackage[a4paper, total={6in, 8in}]{geometry}

 \usepackage{xcolor,cancel}

\newcommand{\Frac}[2]{\displaystyle \frac{#1}{#2}}
\newcommand{\deriv}[3]{\Frac{\partial^{#1} #2}{\partial {#3}^{#1}}}
\newcommand{\Deriv}[3]{\Frac{d^{#1} #2}{d {#3}^{#1}}}

\renewcommand{\vector}[1]{\underline{\boldsymbol{#1}}}
\renewcommand{\matrix}[1]{\underline{\underline{\boldsymbol{#1}}}}

\newcommand{\tr}[1]{\mbox{tr}\left( #1\right)}

\newcommand{\rev}[1]{{\textcolor{black}{#1}}}

\numberwithin{equation}{section}

\makeatletter
\renewcommand{\@biblabel}[1]{#1\hfill \hspace{-0.2cm}}
\makeatother

\usepackage{cite}
\usepackage{enumerate}

\begin{document}

\title{Quantifying Ocular Surface Changes with Contact Lens Wear}
\author{ {\sc Lucia Carichino$^1$,  Kara L.~Maki$^1$, David S.~Ross$^{1,2}$,} \\
{\sc Riley K.~Supple$^1$ and Evan Rysdam$^1$}\\[6pt]
$^1$ School of Mathematics and Statistics, \\
Rochester Institute of Technology, \\
85 Lomb Memorial Drive, Rochester, NY 14623, USA.\\
$^2$ Clinical Pharmacology Science Department, Eisai Inc., \\
200 Metro Blvd, Nutley, NJ 07110, USA..\\[6pt]
}
\pagestyle{headings}
\markboth{L. CARICHINO}{\rm Ocular Surface Changes with Contact Lens Wear}
\maketitle

\begin{abstract} 
Over 140 million people worldwide and over 45 million people in the United \rev{S}tates wear contact lenses; it is estimated \rev{that} $12\%$--$27.4\%$ contact lens users stop wearing them due to discomfort. Contact lens mechanical interactions with the ocular surface have been found to affect the ocular surface itself. The\rev{se} mechanical interactions are difficult to measure and calculate in a clinical setting, and the research in this field is limited. This paper presents the first mathematical model that \rev{captures} the interaction\rev{s} between the contact lens and the open eye, where the contact lens configuration, the contact lens suction pressure, and the deformed ocular shape are all emergent properties of the model. The non-linear coupling between the contact lens and the eye is achieved by assuming that the suction pressure under the lens is applied directly to the ocular surface \rev{through} the post-lens tear film layer. The contact lens \rev{mechanics} \rev{are} modeled using a previous published model. We consider homogeneous and heterogeneous linear elastic eye model\rev{s}, different ocular shapes, different lens shapes and thickness profiles, and extract lens deformation\rev{s}, suction pressure profiles, and ocular deformations and stresses for all the considered scenarios. The model predicts higher ocular deformations and stresses at the center of the eye and in the limbal/scleral regions. Accounting for heterogeneous material eye parameters increases \rev{the magnitude of} such deformations and stresses. The ocular displacements and stresses non-linearly increase as we increase the stiffness of the contact lens. Inserting a steeper contact lens on the eye results in a reduction of the ocular displacement at the center of the eye and a larger displacement at the edge of the contact lens. The model predictions are compared \rev{with} experimental data and previously developed mathematical models. 
\end{abstract}

\maketitle

\section{Introduction}
There are over 140 million people worldwide and over 45 million people in the United States who benefit from wearing contact lenses~\cite{website}. These benefits include vision improvements and one's overall self-perception~\cite{Pucker2020}.  Unfortunately, many people stop wearing contact lenses. In 2020, Pucker and Tichenor~\cite{Pucker2020} conducted a literature review and found that the frequency of contact lens dropout ranged between 12\% and 27.4\%, with the top cited reason being discomfort.  Successful contact lens wear is growing in public health importance.  Two contact lens-based treatments have recently showed promise in slowing down myopic progression \rev{in children}, therefore lessening the chance of developing vision threatening conditions later in life~\cite{Pucker2020}. 

When \rev{placed} on \rev{an} eye, a contact lens interacts with the ocular surface. In particular\rev{,} the contact lens interacts with the corneal tissue, at the center of the eye, and with the limbal and scleral tissue\rev{s} close to the edge of the lens. The cornea is the transparent domed tissue that lets light into the eye, the sclera is the white part of the eye, and the limbus is the transition region between the cornea and the sclera~\cite{snell2013clinical}. There is a thin layer of tears between the contact lens and the ocular surface called the post-lens tear film. Contact lens interactions with the ocular surface have been found to have a significant effect on the anterior ocular surface~\cite{Liesegang2002,alonso2012using,efron2013tfos}. Some of these ocular changes can result in adverse events.  For example, contact lens\rev{es} can limit the supply of oxygen to the ocular surface, and before the invention of silicon hydrogel contact lenses (softer lenses), the \rev{duration} of contact lens wear was restricted to prevent hypoxia. Liesegang~\cite{Liesegang2002} classified these adverse events into four categories: (i)~hypoxia-mediated events, (ii)~immune events, (iii)~mechanical events, and (iv)~osmotic events.  In this work, we focus on quantifying contact lens-induced mechanical events.       

Contact lens design efforts focus on minimizing \rev{their} interference with ocular surface functions and metabolism, while maximizing the patient's comfort and quality of vision~\cite{stapleton2017impact}.  Hall et al.~\cite{Halletal2011} showed that the shape of the ocular surface, measured by optical coherence tomography (OCT), is significantly correlated with the contact lens fit and comfort. \rev{The shapes of the ocular surface and the contact lens affect contact lens-induced mechanical events; the effects of these mechanical interactions have been measured experimentally \cite{efron2013tfos,alonso2012using}. For example,} experimental works have shown that all kinds of contact lenses (soft and stiff) are capable of inducing minimal changes in the cornea topography, and these changes have been found to be statistically significant~\cite{efron2013tfos}. \rev{More specifically}, experimentalists found a small degree of cornea\rev{l} flattening over 9 to 12 months \rev{of continuous wear}~\cite{efron2013tfos}, and significant changes in the surface of the corneal and limbal/scleral regions after 6 hours of soft contact lens wear~\cite{alonso2012using}. There is evidence \rev{that tight-fitting soft contact lenses cause} ``limbal indentations" (i.e., indent\rev{ations in the limbal region of} the ocular surface)~\cite{alonso2012using}. In this work, we build a mathematical model to predict how contact lens wear mechanically affects the ocular surface. 

Various mathematical models have been developed to describe some of the different dynamics that arise from the interactions between contact lenses and the eye~\cite{funkenbusch1996conformity,maki2014new,maki2014exchange,ross2016existence,wu2021biomechanical,zhao2023biomechanical,wu2024fea,ramasubramanian2024influence}. \rev{In Table~\ref{tab:comp_models}, we compare the assumptions of the current work (last row) with the
previous modeling efforts which studied contact lens mechanics when placed on the eye}. In 1996, Funkenbusch and Benson~\cite{funkenbusch1996conformity} published a model to study how the contact lens conforms to the surface of the eye \rev{when the eye is closed}. They model\rev{ed} the contact lens as a thin, axially symmetric elastic shell, and the eye \rev{as} a rigid structure composed of a convex cornea and sclera and a concave limb\rev{al} region. The contact lens \rev{was} applied to the eye under a constant external pressure to mimic the force applied to the lens during a blink or when the eye is closed. \rev{The model predicted the conformity of the lens to the eye, the ocular reaction pressure when the lens and ocular surface are in contact, and the contact lens internal and interfacial loads.} Additionally, in 2024, Ramasubramanian et al.~\cite{ramasubramanian2024influence} studied how a soft contact lens conforms to the ocular surface in the closed\rev{-}eye configuration. \rev{Unlike Funkenbusch and Benson~\cite{funkenbusch1996conformity}, they considered the elastic properties of the eye tissue and fluid pressure inside the eye. Specifically, they} modeled the cornea as a Mooney-Rivlin material that is deformed by \rev{both} the lens and the intraocular pressure (IOP). Similar to Funkenbusch and Benson~\cite{funkenbusch1996conformity}, the lens \rev{was} placed on the eye using \rev{an applied} constant external pressure to mimic the closed \rev{eyelid. Additionally, the post-tear film was modeled by applying a constant pressure on both the contact lens posterior surface and the ocular surface. The model predicted both contact lens and ocular displacements and stresses.}
\begin{table}[h]
    \centering
    \caption{\rev{Comparison between the assumptions of the current work (last row of the table) and previous modeling efforts studying contact lens mechanics when inserted on the eye.}}
    \begin{tabular}{l l l l l l}
    \hline
     Reference & Eye & Cornea model & IOP & Contact lens model\\
\hline    
\multirow{2}{*}{\cite{funkenbusch1996conformity}} & closed, eyelid & \multirow{2}{*}{rigid} & \multirow{2}{*}{neglected} & \multirow{2}{*}{elastic thin shell} \\
 & pressure $1.7$~mmHg & & & \\

\hline 
\multirow{2}{*}{\cite{maki2014new,maki2014exchange,ross2016existence}} & \multirow{2}{*}{open} & \multirow{2}{*}{rigid} & \multirow{2}{*}{neglected} &  linear elastic thin shell \\
& & & & conforming to the eye\\
\hline
\multirow{2}{*}{\cite{wu2021biomechanical}} & closed, eyelid  & \multirow{2}{*}{linear elastic} & \multirow{2}{*}{$13$~mmHg} & not specified\\
 & pressure $9$~mmHg & & & Young's modulus reported\\
 \hline 
\multirow{2}{*}{\cite{zhao2023biomechanical}} & closed, eyelid & \multirow{2}{*}{linear elastic} & \multirow{2}{*}{$13$~mmHg} & not specified\\
& pressure $9$~mmHg & & & Young's modulus reported\\
 \hline 
\multirow{2}{*}{\cite{wu2024fea}} & closed, eyelid & \multirow{2}{*}{linear elastic} & \multirow{2}{*}{$8$--$24$~mmHg} & \multirow{2}{*}{linear elastic}\\
 & pressure $8$--$16$~mmHg & & & \\
 \hline 
\multirow{2}{*}{\cite{ramasubramanian2024influence}} & closed, eyelid & \multirow{2}{*}{Mooney–Rivlin} & \multirow{2}{*}{$15$~mmHg} & \multirow{2}{*}{neo-Hookean}\\
 & pressure $8$~mmHg & & & \\
\hline
 \multirow{2}{*}{This work} & \multirow{2}{*}{open} & \multirow{2}{*}{linear elastic} & \multirow{2}{*}{neglected} & linear elastic thin shell\\
 &  & & &  conforming to the eye\\
 \hline
    \end{tabular}
    
    \label{tab:comp_models}
\end{table}

Multiple modeling efforts have been published to study the interactions between the eye and a type of stiff lens, called an orthokeratology (ortho-k) lens\rev{, which is} worn at night. In 2021, Wu et al.~\cite{wu2021biomechanical} developed a model to study the biomechanical response of the cornea to \rev{an} ortho-k contact lens. \rev{Similar to Ramasubramanian et al.~\cite{ramasubramanian2024influence}, the} lens and the eye \rev{were} modeled as elastic materials and coupled by applying a constant eyelid pressure to the outside of the lens and a constant IOP to the inside of the cornea. \rev{Wu et al.~\cite{wu2021biomechanical} focused on the corneal response and considered corneas of different thicknesses and curvatures. Unlike prior works, Wu et al.~\cite{wu2021biomechanical} implemented a three-layer corneal model by considering a linear elastic model with nonhomogeneous material properties. Details on the lens and cornea interactions were not provided. The biomechanical response of the cornea was characterized for different shaped ortho-k lenses.}
Additionally, Zhao et al.~\cite{zhao2023biomechanical} studied the interaction between \rev{an} ortho-k \rev{lens} and \rev{an eye} using a model similar to~\cite{wu2021biomechanical}. \rev{They extended the work of Wu et al.~\cite{wu2021biomechanical} by considering aspheric corneas and different ortho-k lens designs. The  interactions between the lens and ocular surface were such that the lens could slide and deform the ocular surface. They reported corneal displacements and stresses.}

\rev{In 2024,} Wu et al.~\cite{wu2024fea} \rev{used} a three-part model that include\rev{d} the lens, eye, and eyelid. The eye  \rev{was} modeled as a linear elastic material that deform\rev{ed} due to IOP and the lens, and the lens\rev{, modeled as a linear elastic material,} conform\rev{ed} to the eye due to the eyelid pressure. \rev{Similar to Zhao et al.~\cite{zhao2023biomechanical}, the ocular surface deformed in response to contact from the posterior surface of the lens; however the sliding friction coefficient was nonzero.} \rev{In Wu et al.~\cite{wu2024fea}, the eyelid pressure \rev{was} not imposed a priori, like in~\cite{wu2021biomechanical,zhao2023biomechanical}, but an elastic eyelid surface pressed the contact lens via an applied pressure. Consequently, the applied pressure on the anterior surface of the contact lens was spatially\rev{-}dependent rather than uniform. The} modeling effort account\rev{ed} for the effect of the post-lens tear film by applying a constant pressure between the contact lens and the ocular surface. \rev{Two hundred and forty-nine patient-specific model predictions were compared to patient-specific clinical measured refractive power changes. It was found that the maximum principle stress was most similar (similarity measured by cross-correlation and a structural similarity index) to the clinical data.}

All of the aforementioned models assumed that the eye is closed. In \cite{maki2014new}, Maki and Ross assumed \rev{that} the eye is open, and derived a model for the lens suction pressure under the assumption of axial symmetry and a rigid eye. This model \rev{showed} that the driving property of the contact lens to produce suction pressure is elastic tension or stretching. \rev{Rather than applying an external  pressure, it was assume\rev{d} that} the lens conforms to the shape of the eye \rev{(assumed to be frictionless).} Ross et al.~\cite{ross2016existence} constructed a well-posed system of ordinary differential equations that can be numerically solved for the suction pressure distribution. Maki and Ross extended the model to account for the dynamics of the post-\rev{lens} tear film in~\cite{maki2014exchange}.

In this work, we propose a novel mathematical model to predict ocular deformations due to contact lens wear when the eye is open. The novelty of the proposed model is \rev{that} the eye shape and the \rev{contact} lens suction pressure are all emergent properties of the coupled problem and are not imposed a priori. We build upon the previously published model of contact lens \rev{mechanics} by Maki and Ross~\cite{maki2014new,ross2016existence} and couple the lens \rev{mechanics} to \rev{a} linear elastic model of the ocular deformation via the suction pressure under the lens. \rev{The coupling is non-linear since the suction pressure non-linearly depends on the lens and ocular deformations}. We consider a simplified\rev{, axially symmetric} geometrical description of the eye and the contact lens\rev{. We consider} a homogeneous or a heterogeneous eye model in terms of material parameters (constant or spatially-dependent ocular Young’s modulus). We \rev{assume that the post-lens tear film is a thin, viscous fluid and} neglect the effect of the post-lens tear film \rev{flow} on the lens-eye interactions. 

In Section~\ref{sec:matmet}, we present the ocular and lens geometrical domains, the eye model, the lens model, the boundary and coupling conditions, and \rev{our} numerical \rev{method}. In Section~\ref{sec:results}, we first explore the results obtained with the homogeneous eye model for different contact lens stiffnesses, different ocular and lens shapes, and a constant or spatially varying lens thickness; then, we explore the effect of a heterogeneous eye model (\rev{one with a} spatially-depend\rev{ent} Young's modulus) on the lens-eye interactions.  In Section~\ref{sec:discussion}, we compare the results \rev{with} previously published experimental studies and mathematical model predictions, and \rev{we} discuss the \rev{limits} of the current work.

\section{Materials and methods}
\label{sec:matmet}
In this section, we present the mathematical \rev{model} and \rev{the} numerical \rev{method we} developed to predict the deformations of the eye and of the contact lens when a contact lens is placed on the ocular surface. \rev{W}e begin by describing the ocular and contact lens domains. Then, we introduce the model which characterizes the deformations of the eye and present the model of the contact lens mechanics.  Finally, we characterize the interactions between the contact lens and \rev{the} ocular surface, and we define the numerical \rev{method we} developed to solve the coupled problem.  

\subsection{Geometries of the eye and contact lens}
\begin{figure}[htb]
\centering
\scalebox{0.65}{
\begin{tikzpicture}[scale=1]
     	\draw[->,thin,color=gray] (0,-2.5) -- (0,10);
	\coordinate [label=left:\textcolor{gray}{$z$}] (E) at (0,10);
	\draw[->,thin,color=gray] (-1,4.63) -- (7.5,4.63);
	\coordinate [label=right:\textcolor{gray}{$r$}] (E) at (7.5,4.63);
     \draw[color=gray,dashed] (6.15,4.63) -- (6.15,-1.5); 
	\coordinate [label=above:\textcolor{black}{$\mathfrak{R}_{eye}$}] (E) at (6.15,4.63);	 
		
	\draw[thick,dotted,color=gray] ({4*cos(60)},{4*sin(60)}) arc(80:-10:5);	
	\draw[thick] ({4*cos(60)},{4*sin(60)}) arc(80:-0.3:5);		
	\draw[thick] ({4*cos(60)},{4*sin(60)}) arc(30:90:2.3);	
    \draw[blue,fill=blue] (3.2,3.1) circle (.4ex);     
    \draw[thick,blue,dashed] (3.2,3.1) -- (3.2,4.63); 
    \coordinate [label=above:\textcolor{blue}{$\chi(\mathfrak{R}_{lens})$}] (E) at (3.9,4.55);  	
    \draw[<->,thin,color=gray] (0.3,8.045) -- (0.3,8.24);
    \coordinate [label=below:\textcolor{gray}{$\tau$}] (E) at (0.3,8.04);

\coordinate [label=below:\textcolor{black}{$\Omega$}] (E) at (2,1);  	
\coordinate [label=below:\textcolor{black}{undeformed eye}] (E) at (1.7,0.5);                  
	\draw[-, thick] (0,-1.5) -- (0,4.6);
 	\draw[-, thick] (0,-1.5) -- (6.15,-1.5);     
        \coordinate [label=below:\textcolor{black}{$\Gamma_{in,h}$}] (E) at (2,-1.5);  
        \coordinate [label=left:\textcolor{black}{$\Gamma_{in,v}$}] (E) at (0,2);  		
        \coordinate [label={[rotate=-20]above:\textcolor{black}{cornea}}] (E) at (0.65,3.95);  
         \coordinate [label={[rotate=-82]above:\textcolor{black}{sclera}}] (E) at (5.5,-0.7);  	
         \coordinate [label={[rotate=-25]below:\textcolor{black}{limbus}}] (E) at (2.3,3.4);         
                     	
        \coordinate [label=below:\textcolor{black}{$\Gamma_{out}$}] (E) at (2,4.5);  
        \coordinate [label=below:\textcolor{black}{$h(r)$}] (E) at (5.6,1.9);          

\draw[thick,dotted,color=gray] ({4*cos(60)},{4*sin(60)}) -- ({-4*cos(60)},{4*sin(60)});	
  
	\begin{scope}[yscale=1,xscale=-1]
	\draw[thick, dotted,color=gray] ({4*cos(60)},{4*sin(60)}) arc(30:90:2.3);	
        \draw[thick,dotted,color=gray] ({4*cos(60)},{4*sin(60)}) arc(80:70:6);	              
	\end{scope}	

\begin{scope}[yshift=30]

	\begin{scope}[yscale=1,xscale=-1]       
	\draw[thick,dotted,color=gray] (0,7) arc(90:35:3.5);	
 	\draw[thick,dotted,color=gray] (0,7.2) arc(90:29:3.3);        
	\end{scope}	

	\draw[thick] (0,7) arc(90:35:3.5);	
 	\draw[thick] (0,7.2) arc(90:29:3.3);
    \draw[blue,fill=blue] (2.867,5.5) circle (.4ex); 

\draw[->,thin,color=gray] (-1,7.2) -- (4,7.2);
\coordinate [label=right:\textcolor{gray}{$r$}] (E) at (4,7.15);
\draw[thick,dashed,color=gray] (2.867, 7.2) -- (2.867,5.5);     
\coordinate [label=below:\textcolor{black}{$\mathfrak{R}_{lens}$}] (E) at (3.3,7.2);	
\coordinate [label=left:\textcolor{black}{$g(r)$}] (E) at (2,6.2); 
\coordinate [label=right:\textcolor{black}{undeformed contact lens}] (E) at (0,7.8);  
\end{scope}

\path[-{Stealth[length=4mm, width=2.5mm]}, fill=blue, color=blue, every node/.style={font=\sffamily\small}]
    (3.2,3.1) edge[bend right] node [left] {} (13.65,3.3);
\coordinate [label=above:\textcolor{blue}{$R_{eye}$}] (E) at (8,1.6);  
\coordinate [label=above:\textcolor{blue}{$R_{lens}$}] (E) at (8,3.2); 
\coordinate [label=above:\textcolor{blue}{$\chi$}] (E) at (2.3,4.6); 

\path[color=blue, every node/.style={font=\sffamily\small}]
    (2.867,6.55) edge[bend right] node [left] {} (13.5,3.185);

    \draw[->, thick,color=red] (.75,6.25) -- (.75,5);
    \coordinate [label=above:\textcolor{red}{$P_{out}$}] (E) at (.75,6.3); 

\path[-{Stealth[length=4mm, width=2.5mm]}, fill=blue, color=blue, every node/.style={font=\sffamily\small}]
    (2.867,6.55) edge[bend right] node [left] {} (3.15,3.16);

    
\begin{scope}[xshift=300]
     	\draw[->,thin,color=gray] (0,-2.5) -- (0,6);
	\coordinate [label=left:\textcolor{gray}{$z$}] (E) at (0,6);
	\draw[->,thin,color=gray] (-1,4.63) -- (7.5,4.63);
	\coordinate [label=right:\textcolor{gray}{$r$}] (E) at (7.5,4.63);;

\begin{scope}[yscale=1.03]
	\coordinate [label=right:\textcolor{black}{$R_{eye}(\mathfrak{R}_{eye},h(\mathfrak{R}_{eye}))$}] (E) at (6.15,4);
 	\draw[thick, dashed,color=gray] (6.15,-1.5) -- (6.15,4.5);  
        	\coordinate [label=above:\textcolor{blue}{$R_{eye}\left(\chi(\mathfrak{R}_{lens}), h\left(\chi(\mathfrak{R}_{lens})\right)\right)=R_{lens}(\mathfrak{R}_{lens})$}] (E) at (3.2+0.5,4.5);	
	\draw[thick, dashed,color=blue] (3.2,4.5) -- (3.2,3.35);

  \begin{scope}[yscale=1.07]
	\draw[thick] ({4*cos(60)},{4*sin(60)}) arc(80:0.8:5);	
  \end{scope} 
  \begin{scope}[yscale=0.98, yshift=9.5]  
	\draw[thick] ({4*cos(60)},{4*sin(60)}) arc(50:92:2.9);	
  \end{scope}  
\begin{scope}[yshift=6]
    \begin{scope}[yscale=1.7, yshift=-36]
	\draw[thick] ({4*cos(60)},{4*sin(60)}) arc(80:65:5);	
  \end{scope} 
  \begin{scope}[yscale=0.98, yshift=9.5]  
	\draw[thick] ({4*cos(60)},{4*sin(60)}) arc(50:92:2.9);	
  \end{scope}
  \end{scope}
 \draw[blue,fill=blue] (3.2,3.3) circle (.4ex);

\coordinate [label=below:\textcolor{black}{deformed eye}] (E) at (1.7,0.5);                  
	\draw[-, thick] (0,-1.5) -- (0,4.62);
 	\draw[-, thick] (0,-1.5) -- (6.15,-1.5);     		
    \coordinate [label=right:\textcolor{black}{deformed contact lens}] (E) at (0,5.5);  	       
                          
        \coordinate [label=below:\textcolor{black}{$H(r)$}] (E) at (5.2,2.6);            
\end{scope}

\end{scope}
 
\end{tikzpicture}}
\caption{A schematic of the reference frames of the eye and the contact lens models (left) and of the coupled deformed eye and contact lens (right). In blue, we highlight the deformation map for the end-point of the undeformed lens  domain, located at radial coordinate $\mathfrak{R}_{lens}$, and the corresponding point on the undeformed ocular domain, located at radial coordinate $\chi\left(\mathfrak{R}_{lens}\right)$, to the deformed configuration located at radial coordinate $R_{eye}(\chi\left(\mathfrak{R}_{lens}\right),h(\chi\left(\mathfrak{R}_{lens}\right)))=R_{lens}(\mathfrak{R}_{lens})$.}
\label{fig:contact}
\end{figure}
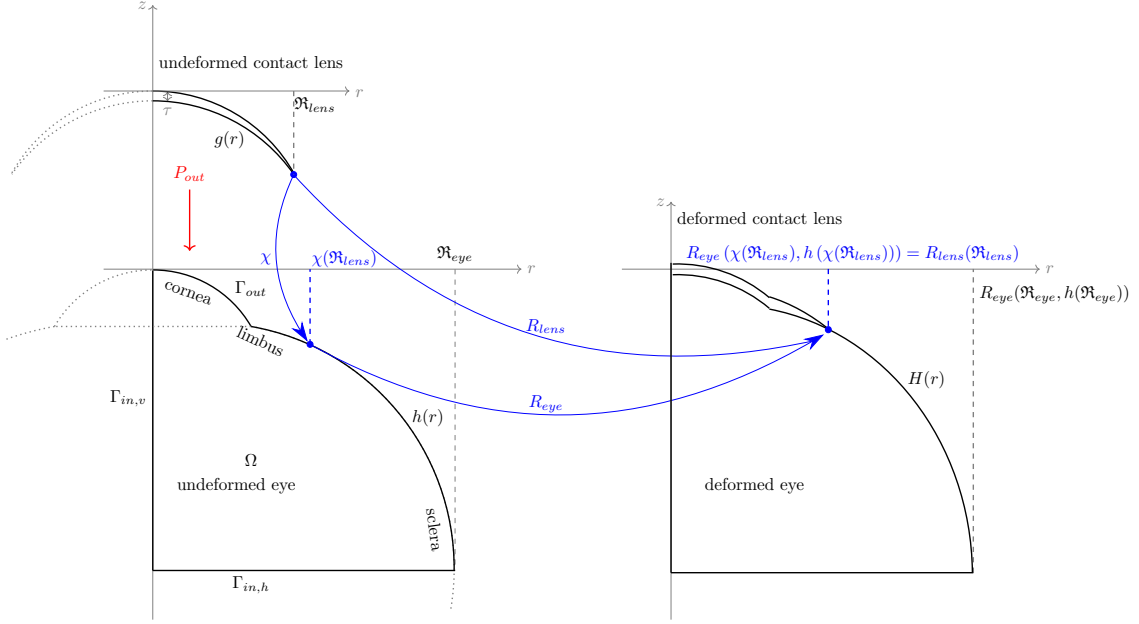

To explore the mechanical interactions between the contact lens and the ocular surface, we start with a simplified geometrical description of the eye and the contact lens. The eye and contact lens are assumed to be axisymmetric, as shown in Figure~\ref{fig:contact}.  This allows us to reduce the computational domain of the contact lens and of the eye from three dimensions (3D), $(r,\theta,z)$,  where $r$, $\theta$, and $z$ denote the radial, azimuthal, and axial coordinates, respectively, to two dimensions (2D), $(r,z)$. In particular, as shown in Figure~\ref{fig:contact}, the 2D ocular domain $\Omega$ represents the upper right quadrant of a slice of the eye in the radial direction. The 2D lens domain, as shown in Figure~\ref{fig:contact}, represents half of the cross section of the lens in the radial direction. We note that this is a simplification, as a patient's ocular topography is meridian dependent \cite{Halletal2011}, and toric contact lenses are also meridian dependent~\cite{doll2020feature}. 

\begin{figure}[htb]
    \centering
\begin{tikzpicture}
    \node at (0,0) {\includegraphics[width=.5\linewidth]{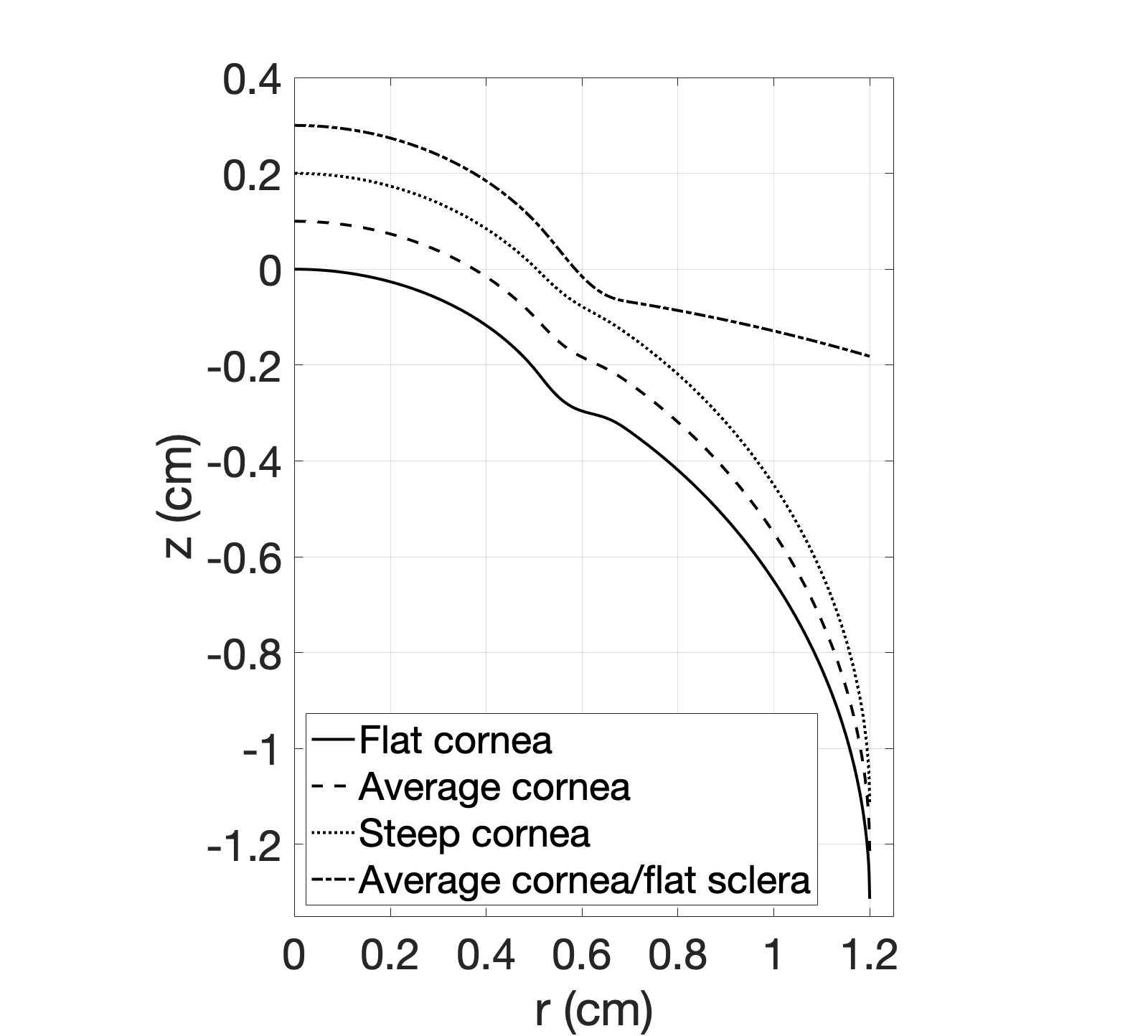}};
    \coordinate [label={left:\textcolor{black}{\bf \large A}}] (A) at (-2.8,3.5);
    \node at (8,0) {\includegraphics[width=.5\linewidth]{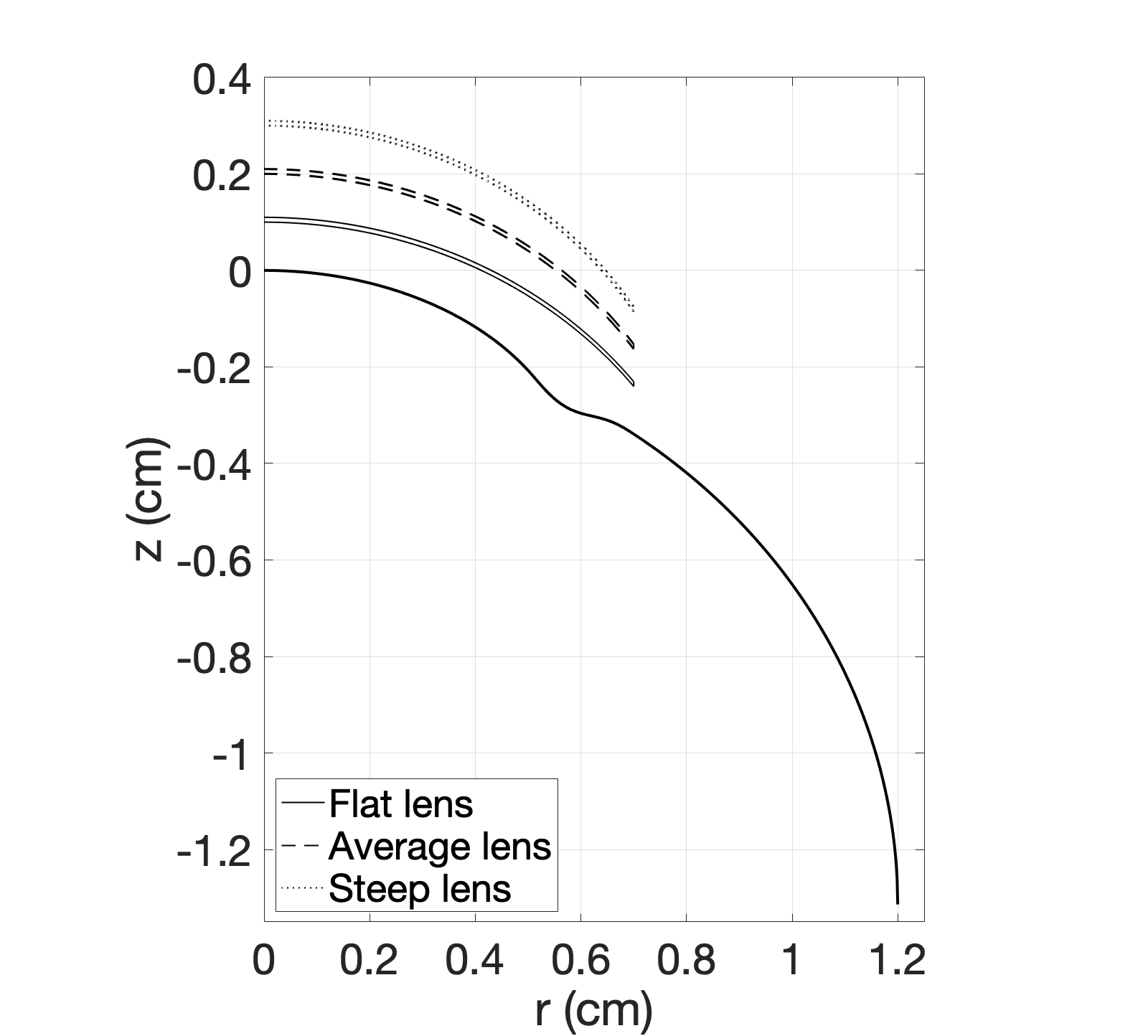}};
    \coordinate [label={left:\textcolor{black}{\bf \large B}}] (B) at (5,3.5);
\end{tikzpicture}
    \caption{({\bf A}) Different ocular surface profiles considered. ({\bf B}) Different contact lens profiles considered, placed on an average-shaped eye. All the lens profiles in ({\bf B}) and the ocular surfaces profiles in ({\bf A}), except the flat cornea, have been shifted vertically to facilitate illustrating the differences.}
    \label{fig:eye_geometry}
\end{figure}

The boundary of the eye domain $\partial \Omega$ is divided in\rev{to} the three \rev{segments}, $\Gamma_{in,v}$, $\Gamma_{in,h}$, and $\Gamma_{out}$, such that $\partial \Omega = \Gamma_{in,v} \cup \Gamma_{in,h} \cup \Gamma_{out}$ \rev{(i.e., the left, bottom, and outer (cornea, limbus and sclera) boundaries of $\Omega$, respectively, as shown in Figure~\ref{fig:contact})}. 
$\Gamma_{out}$ represents the ocular surface and is defined as follows:
\begin{equation}
\label{eq:gamma_out}
\Gamma_{out}= \displaystyle\left\{ (r, z) \in \mathbb R^2 \, : \, r\in [0,\mathfrak{R}_{eye}], \, z=h(r)   \right\},
\end{equation}
where $h(r)$ characterizes the undeformed ocular surface such that $h(0)=0$, and $\mathfrak{R}_{eye}$ is the maximum value attained by $r$ in the undeformed ocular domain. Note that $\mathfrak{R}_{eye}$ is less than or equal to the radius of the eye globe, for more details, see Section~\ref{sec:supp_eye} of the Appendix. 
$\Gamma_{in,v}$ and $\Gamma_{in,h}$ are artificial boundaries introduced by the adopted geometrical domain reduction, and are defined as follows:
\begin{align}
\Gamma_{in,h} = \displaystyle\left\{ (r, z) \in \mathbb R^2 \, : \, r\in [0,\mathfrak{R}_{eye}], \, z =h(\mathfrak{R}_{eye}) \right\} \quad \mbox{and} \quad
\Gamma_{in,v} = \displaystyle\left\{ (r, z) \in \mathbb R^2 \, : \, r=0, \, z \in \left[h(\mathfrak{R}_{eye}), 0\right]  \right\}.
\end{align}
Hence, the domain $\Omega$ can be characterized as follows:
\begin{equation}
\Omega= \displaystyle\left\{ (r, z) \in \mathbb R^2 \, : \, r\in [0,\mathfrak{R}_{eye}], \, z \in \left[h(\mathfrak{R}_{eye}), h(r)\right] \right\}.
\end{equation}
Note that since $h(0)=0$, the upper left corner of the domain $\Omega$, depicted in Figure~\ref{fig:contact}, is at the origin $(0,0)$ of the $rz$-plane and the ocular domain lies below the $r$ axis.
To characterize the function $h(r)$ which describes $\Gamma_{out}$, we divide the ocular surface into the three \rev{segments}: (i) the corneal \rev{segment}, (ii) the limbal \rev{segment}, and (iii) the scleral \rev{segment}, as shown in Figure~\ref{fig:contact}. 
We consider four different ocular shapes, which are based on biometric measurements reported in Hall et al\rev{.}~\cite{Halletal2011}. 
Three of the ocular shapes only vary in the \rev{shape of the} cornea \rev{segment} \rev{in $\Gamma_{out}$}: the cornea is either flat, average, or steep. For the fourth eye shape, the \rev{shape of the} cornea \rev{segment} is average, but the \rev{shape of the} sclera \rev{segment} is flat. Figure~\ref{fig:eye_geometry}{\bf A} displays the ocular surfaces, where all except the flat cornea \rev{shape} have been shifted vertically to facilitate illustrating differences. Section~\ref{sec:supp_eye} of the Appendix reports the expression \rev{for} $h(r)$, its derivation, and how the different ocular shapes are constructed\rev{,}.

The contact lens reference domain is characterized by a posterior curve, $z=g(r)$, and a thickness profile, $\tau(r)$, for $r\in[0,\mathfrak{R}_{lens}]$, where $\mathfrak{R}_{lens}$ is the reference radius of the lens (see Figure~\ref{fig:contact}). Hence, the anterior curve of the lens is given by $z = g(r) + \tau(r)$ for $r\in[0,\mathfrak{R}_{lens}]$. We consider three different lens shapes with constant thickness $\tau=100$~$\mu$m: a flat lens, \rev{an} average lens, and a steep lens. The lens\rev{es} are shown in Figure~\ref{fig:eye_geometry}{\bf B}, again shifted vertically to facilitate the display.  
Additionally, we consider an averaged-shape contact lens with a \rev{varying} thickness to match the lens thickness profile \rev{described} in Funkenbus\rev{c}h and Benson \cite{funkenbusch1996conformity}. \rev{The varying thickness lens case is not shown in Figure~\ref{fig:eye_geometry}{\bf B}; however, the profile of $\tau(r)$ is shown in Figure~\ref{fig:tau} in the Appendix, Section~\ref{sec:supp_lens}.
Additional d}etails \rev{of} the expression \rev{for} $g(r)$, its derivation, how the different lens shapes are constructed, and the \rev{varying} lens thickness profile are reported in the Appendix, Section~\ref{sec:supp_lens}.

\subsection{Biomechanics of the eye}\label{sec:eye_model}
We model the eye as an isotropic, linear elastic material that deforms because of its non-linear coupling with the contact lens.  The displacement vector of the eye is denoted by $\vector{u}=(u_r,0,u_z)$, where $u_{\theta}~=~0$ due to the axial symmetry assumption.  Therefore, the coordinates of the deformed eye $(R_{eye},Z_{eye})$ are given by the following:
\begin{equation}
    R_{eye}(r,z) = r + u_r(r,z) \quad \text{and} \quad Z_{eye}(r,z) = z + u_z(r,z),
\end{equation}
for $r\in[0,\mathfrak{R}_{eye}]$ and $z\in[h(\mathfrak{R}_{eye}), h(r)]$. The ocular surface of the deformed eye is as follows: 
\begin{equation}
\label{eq:H}
H\left(R_{eye}(r,h(r))\right) = Z_{eye}(r,h(r)) = h(r) + u_z(r,h(r)) \quad \mbox{for} \quad r\in[0,\mathfrak{R}_{eye}].
\end{equation}
\rev{We remind the reader of the difference in notation between $R_{eye}$ and $\mathfrak{R}_{eye}$:} $R_{eye}$ is the radial coordinate of the deformed ocular domain, while $\mathfrak{R}_{eye}$ is the maximum value attained by $r$ in the undeformed ocular domain. 

\rev{Assuming the ocular tissue as an isotropic, linear elastic material~\cite{landau1986lifshitz}, t}he displacement vector of the eye $\vector{u}$ satisfies the following force equilibrium equations:
\begin{equation}\label{eq:S}
\nabla \cdot \matrix{S} = \vector{0} \quad \mbox{ in } \quad \Omega \subset \mathbb R^2,
\end{equation}
where $\matrix{S} = \lambda_{eye} \,\mbox{tr}(\matrix{E}) \matrix{I} + 2\mu_{eye} \matrix{E}$ is the stress tensor, 
$\matrix{E}~=~\frac{1}{2}[\nabla \vector{u} + (\nabla \vector{u})^T]$
is the strain tensor, \rev{$\mbox{tr}()$ is the trace, $\matrix{I}$ is the identity matrix}, and $\lambda_{eye}$ and $\mu_{eye}$ are the spatially-dependent eye Lam\'e parameters\rev{, given in Eq~\eqref{eq:lame}}.  We characterize the elastic properties of the eye in terms of the spatially-dependent Young's modulus, $E_{eye}$, and the Poisson's ratio, $\sigma_{eye}$ ~\cite{landau1986lifshitz}.  The relationships between the Young's modulus, Poisson's ratio, and the Lam\'e constants are as follows:
\begin{equation}
\label{eq:lame}
    \mu_{eye}(r,z) = \frac{E_{eye}(r,z)}{2(1 + \sigma_{eye})} \quad \text{and} \quad \lambda_{eye}(r,z) =\frac{E_{eye}(r,z) \sigma_{eye}}{(1 + \sigma_{eye})(1-2\sigma_{eye})}.
\end{equation}

In this work, we initially assume that the ocular domain is homogeneous (i.e., we assume that Young's modulus $E_{eye}(r,z)=E_{eye}$ is constant throughout the ocular domain $\Omega$). Since, our focus is to predict the deformations of the eye near the ocular surface, we use the reported material properties of the cornea and sclera to choose the \rev{possible} range of the Young's modulus $E_{eye}$ constant, which is from $0.2$ to $1.0$~MPa, reported in~Table~\ref{tab:material}. We model the eye as a nearly incompressible material, where $\sigma_{eye}=0.49$. More information on how we chose the value of the constant $E_{eye}$ can be found in the Appendix, Section \ref{sec:supp_eye_para}. \rev{Moreover, the choices of material parameters for a contact lens are reported in Table~\ref{tab:material} for direct comparison, and they are discussed in Section~\ref{sec:lens_model}.}

\begin{figure}[htb]
    \centering
    \includegraphics[width=0.6\linewidth]{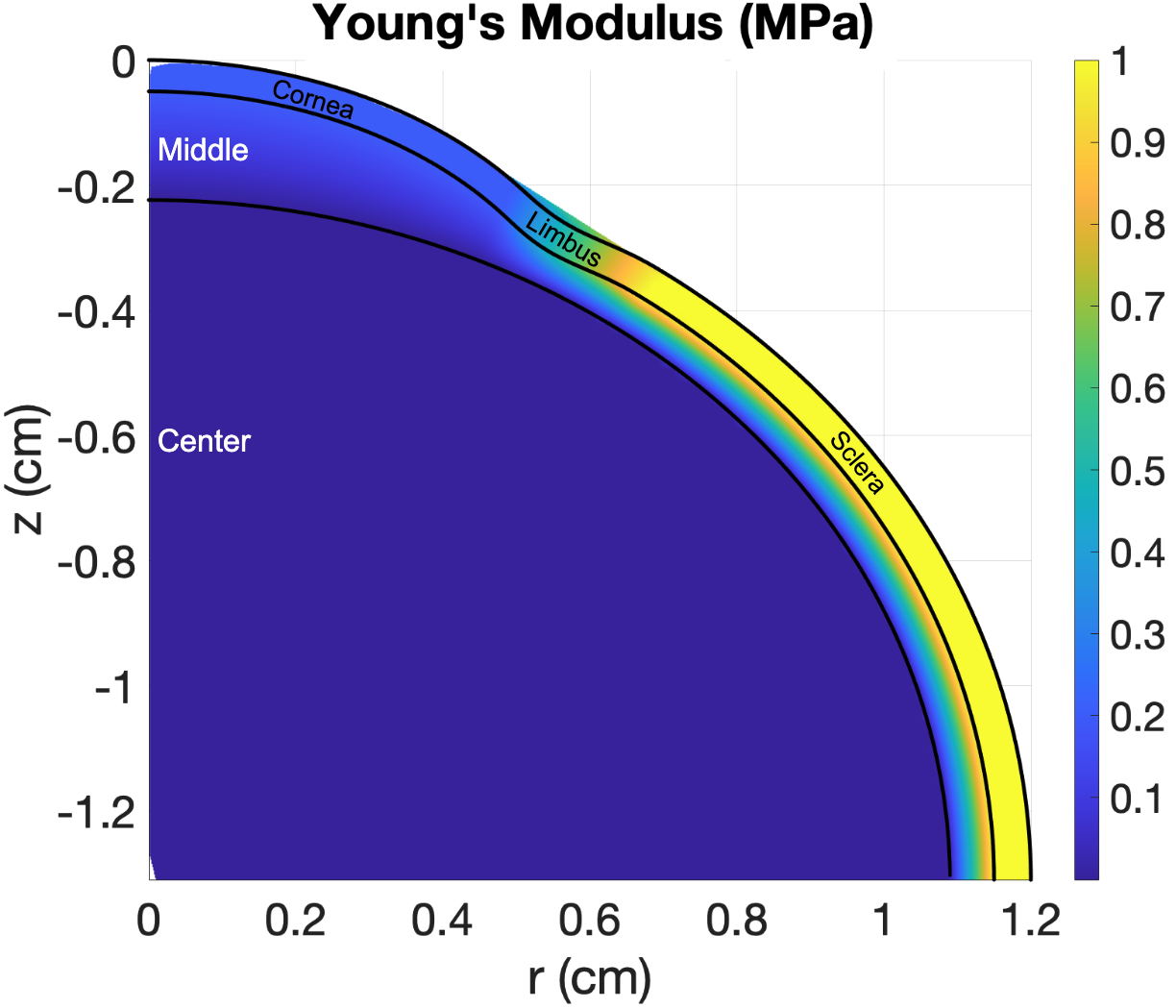}
    \caption{Spatially\rev{-}dependent Young's modulus values in the different ocular regions considered.}
    \label{fig:simplified_varied_ym}
\end{figure}

Later, to model the complex human eye anatomy, we consider a \rev{heterogeneous} model that accounts for the different responses to stresses of the cornea, sclera, and of the inside of the eye (composed mainly of gel\rev{-}like fluids). To do so, we consider a spatially-dependent Young's modulus, as shown in Figure~\ref{fig:simplified_varied_ym}. \rev{We consider an outer region, which represents the cornea, limbus, and sclera tissues, and a center region which represents the vitreous humor chamber. The middle region (i.e., in between the outer and center regions) represents the aqueous humor chamber, iris, eye lens, retinal, and choroidal tissues.} The values of $E_{eye}(r,z)$ in each region are determined according to the anatomical structure of the eye and to previous literature. We assume that the Young’s modulus in the sclera is five times the Young’s modulus in the cornea of $0.2$ MPa~\cite{Wooetal1972}. 
Inside the eye, which is represented by the center region in Figure~\ref{fig:simplified_varied_ym}, we set the Young’s modulus to be $1.17 \times 10^{-6}$~MPa~\cite{Tram2021} to mimic the vitreous humor \rev{chamber}. We assume linear transitions in the limbus region and in the middle region shown in Figure~\ref{fig:simplified_varied_ym}. More details on the construction of the spatially-dependent Young's modulus are reported in the Appendix, Section~\ref{sec:supp_YM}. 

\begin{table}[htb]
    \centering
        \caption{Material parameters of the eye and the soft contact lens, and the corresponding sources.}
    \begin{tabular}{ lll}
    \toprule
       Parameter  & Value & Reference\\
    \midrule
      Constant eye Young's modulus, $E_{eye}$  & 0.2--1.0 MPa  & \cite{BryantandMcDonnell1996,Wooetal1972} \\
      Eye Poisson's ratio, $\sigma_{eye}$ & 0.49 & \cite{Yehetal2000,Pandolfi2006} \\
      Soft contact lens Young's modulus, $E_{lens}$ & 0.1--2 MPa & \cite{tranoudis2004tensile,horst2012measuring,young2010effect,Quinceetal2021} \\
      Contact lens Poisson's ratio, $\sigma_{lens}$ & 0.49 & \cite{BhamraTighe2017} \\
      \hline
    \end{tabular}

    \label{tab:material}
\end{table}

\subsection{Mechanics of the contact lens} \label{sec:lens_model}
We use a mathematical model developed by Maki and Ross~\cite{maki2014new,ross2016existence} to describe the contact lens mechanics. \rev{Maki and Ross~\cite{maki2014new,ross2016existence} reduced the complexity of the model by assuming that} the contact lens is a linear elastic \rev{thin} shell, and thus the lens displacements are independent of the axial $z$ coordinate.  We denote the lens displacements by $\vector{\eta}=(\eta_r,0,\eta_z),$ where $\eta_{\theta}=0$ due to symmetry. 
Additionally, \rev{Maki and Ross~\cite{maki2014new,ross2016existence} derived} the model under the assumption that the deformed contact lens completely conforms to the ocular surface\rev{. In~\cite{maki2014new,ross2016existence}, the eye was rigid so the ocular surface was assumed to be fixed, while in this work, the ocular surface can deform and is} described by $z=H\left(R_{eye}(r,h(r))\right)$ in Eq~\eqref{eq:H}.  
For a given $r \in [0,\mathfrak{R}_{lens}]$ in the lens reference frame, the deformed radial coordinate of the contact lens is $R_{lens}(r) = r + \eta_r(r)$.  We remind the reader that $\mathfrak{R}_{lens}$ is the radius of the undeformed contact lens. 

\rev{Since the lens deforms to conform to the eye, while the eye deforms due to the lens suction pressure, w}e anticipate the radial displacement of the contact lens $\eta_r(r)$ to be {\it different} from the radial displacement of the ocular surface $u_r(r,h(r))$ for a fixed value of $r \in [0,\mathfrak{R}_{lens}]$. Hence, the deformed radial coordinates of the contact lens and the eye are different (i.e., $R_{lens}(r)= r + \eta_r(r) \neq r + u_r(r,h(r)) = R_{eye}(r, h(r))$, for a fixed value of $r \in [0,\mathfrak{R}_{lens}]$).
Therefore, for each $r \in [0,\mathfrak{R}_{lens}]$, there will exist $\chi(r) \in [0,\mathfrak{R}_{eye}]$ such that 
\begin{equation}
\label{eqn:PHYSICAL_R}
R_{eye}\left(\chi, h\left(\chi\right)\right)=\chi + u_r\left(\chi,h\left(\chi\right)\right) = r + \eta_r(r) = R_{lens}(r).
\end{equation} 
The function $\chi$ is a differentiable function with the domain $[0,\mathfrak{R}_{lens}]$ and the range $\left[0,\chi(\mathfrak{R}_{lens})\right]$, \rev{that can be expressed as the composition of the function $R_{lens}$ with the inverse of the function $R_{eye}$ (i.e, $\chi = R^{-1}_{eye} \circ R_{lens}$)} as depicted in blue in Figure~\ref{fig:contact}.  Therefore, the deformed contact lens lays above the deformed eye surface for $r$ in the undeformed eye domain $0 \leq r \leq \chi(\mathfrak{R}_{lens})$. 

Now that we have \rev{introduced} the function $\chi$, we can express the part of the deformed ocular surface $H$, Eq~\eqref{eq:H}, that lies under the deformed lens as follows:
\begin{equation}
\label{eq:H2}
H(R_{lens}(r)) = H\left(R_{eye}\left(\chi\rev{(r)}, h\left(\chi\rev{(r)}\right)\right)\right) = h(\chi\rev{(r)}) + u_z(\chi\rev{(r)},h(\chi\rev{(r)}))  \quad \mbox{for} \quad r\in[0,\mathfrak{R}_{lens}].
\end{equation}
Hence, we can express the contact lens axial displacement $\eta_z(r)$ as follows, under the assumption that the deformed contact lens conforms to the deformed eye surface:
\begin{equation}
\eta_z(r) = H\left(R_{lens}(r)\right) - g(r) \quad \mbox{for} \quad r\in[0,\mathfrak{R}_{lens}],
\end{equation} 
where $g$ is reference shape of the lens posterior curve (i.e., before being placed on the eye), which is sketched in Figure~\ref{fig:contact} and described in the Appendix, Section~\ref{sec:supp_lens}. 

To find the unknown function, $\eta_r(r)$, or $R_{lens}(r) = r + \eta_r(r)$,  \rev{Maki and Ross~\cite{maki2014new,ross2016existence}} derive\rev{d} a system of non\rev{-}linear ordinary differential equations which represent the Euler-Lagrange equations associated with minimizing the energy functional of the contact lens system as follows:
\begin{eqnarray}
T^\prime(r) &=& \left(\sigma_{lens} T(r) + (1-\sigma_{lens}^2)(R_{lens}(r)-r)\right) \frac{\sqrt{1 + g^{\prime}(r)^2}}{r\sqrt{1 + H^{\prime}(R_{lens})^2}} - \frac{\tau^{\prime}(r)}{\tau(r)}T(r), \label{eq:odeT}\\
R_{lens}^\prime(r) &=& \left(T(r) + (1+\sigma_{lens})r-\sigma_{lens} R_{lens}(r)\right)\frac{\sqrt{1 + g^{\prime}(r)^2}}{r\sqrt{1 + H^{\prime}(R_{lens})^2}}, \label{eq:odeRl}
\end{eqnarray}
for $r\in[0,\mathfrak{R}_{lens}]$, where $\tau(r)$ is the thickness of the contact lens, and $\sigma_{lens}$ is the Poisson's ratio of the contact lens. \rev{$T(r)$ is the contact lens radial tension scaled by the thickness of the lens $\tau(r)$, the radial coordinate $r$, and the lens material parameters; for more details, see Ross and Maki~\cite{ross2016existence}.} The boundary conditions of Eqs~\eqref{eq:odeT}--\eqref{eq:odeRl} are as follows:
\begin{eqnarray}
    T(0)=R_{lens}(0)&=&0, \label{eq:bcT1} \\
    T(\mathfrak{R}_{lens}) &=& 0. \label{eq:bcT2}
\end{eqnarray}
The detailed derivation of Eqs~\eqref{eq:odeT}--\eqref{eq:bcT2} can be found in \cite{maki2014new,ross2016existence}. \rev{In general, the system of ordinary differential equations in Eqs~\eqref{eq:odeT}--\eqref{eq:bcT2} predicts the deformed radial coordinate of the lens that requires the least amount of energy to conform the contact lens to the ocular surface.} Note that the non\rev{-}linear coupling between the eye and contact lens models in Eqs~\eqref{eq:odeT}--\eqref{eq:odeRl} is by the deformed ocular surface, $H$, which is given by Eq~(\ref{eq:H2}). Therefore, we need to know $H$ to solve the system of equations in Eqs~\eqref{eq:odeT}--\rev{\eqref{eq:bcT2}}.

Given the deformed contact lens coordinate, $R_{lens}(r)$, and the deformed eye shape, $H$, the suction pressure, $p$, under the lens can be determined from an ordinary differential equation, as detailed in~\cite{maki2014new,ross2016existence}. In particular, for a given \rev{$R_{lens}(r)$} and $r \in \rev{(}0,\mathfrak{R}_{lens}]$, 
\begin{equation}\label{eq:p}
\begin{aligned}
\rev{
p(R_{lens}) = - \frac{1}{R_{lens}} \Deriv{}{}{R_{lens}} \left[ \Frac{R_{lens} H^\prime(R_{lens})}{\sqrt{1 + H^\prime(R_{lens})^2}}
\left.\left(
\dfrac{E_{lens}}{1-\sigma_{lens}^2} \dfrac{\tau(r)}{r} T(r)
\right) \right|_{r + \eta_{r}(r)=R_{lens}}
\right],
}
\end{aligned}
\end{equation}
where $E_{lens}$ and $\sigma_{lens}$ are the Young's modulus and Poisson's ratio of the contact lens, respectively.
\rev{At $R_{lens}=0$, Eq~\eqref{eq:p} has a singularity; however, an analytical expression for $p(0)$ is found by taking the limit as $R_{lens}$ approaches zero. Note that when $p>0$, the contact lens pushes down on the ocular surface, and when $p<0$, the lens pulls up on the ocular surface.}

We assume the contact lens to be nearly incompressible; thus, $\sigma_{lens}=0.49$. The Young's modulus of a soft contact lens can significantly vary depending on the manufacturer and on the material used by the manufacturer. Following previously published literature~\cite{tranoudis2004tensile,horst2012measuring,young2010effect,Quinceetal2021}, we consider the soft contact lens Young's modulus \rev{$E_{lens}$} range reported in Table~\ref{tab:material}. 

\subsection{Boundary conditions for the eye model}
Now that we have characterized the eye and contact lens mechanics, we present the boundary conditions for the eye model needed to close the model system.  The eye model and the lens model are coupled via the boundary condition on the external surface of the eye $\Gamma_{out}$, which is defined in Eq~\eqref{eq:gamma_out} and depicted in Figure~\ref{fig:contact}.
We impose continuity of traction on $\Gamma_{out}$ as follows:  
\begin{equation}
\label{pressure_bc}
\matrix{S} \ \vector{n} = - P_{out} \ \vector{n} \quad \mbox{on }\Gamma_{out},
\end{equation}
where $\vector{n}$ is the outward normal unit vector to $\Gamma_{out}$, and $P_{out}$ is the external pressure acting on the eye. Note that in Eq~\eqref{pressure_bc}, we are imposing that the component of the traction normal to the ocular surface is equal to $-P_{out}$ and the tangential component of the traction is zero.

On $\Gamma_{out}$, the external pressure $P_{out}$ varies along the radial coordinate $r\in[0,\mathfrak{R}_{eye}]$ as follows:
\begin{equation}\label{Pe}
P_{out}(r) = 
\left\{
\begin{array}{cl}
p\left( \rev{R_{lens}\left(\chi^{-1}(r)\right)} \right) & \mbox{for }0 \leq r \leq \chi(\mathfrak{R}_{lens}),\\
0 & \mbox{for }  \chi(\mathfrak{R}_{lens}) < r \leq \mathfrak{R}_{eye}.
\end{array}
\right. 
\end{equation}
For a given eye radial coordinate $r \in [0, \chi(\mathfrak{R}_{lens})]$, the corresponding deformed radial coordinate of the eye surface is ``under" the deformed contact lens at $R_{lens}\left(\chi^{-1}(r)\right)$. Therefore, for those radial coordinates, we impose the mechanical pressure exerted by the contact lens $p$, given in Eq~\eqref{eq:p}, on the eye surface. Note that, if $p>0$, then the lens pushes on the ocular surface; hence, in Eq~\eqref{pressure_bc}, the traction $\matrix{S} \ \vector{n}<\vector{0}$. For the ocular tissue that is not ``under" the lens, we impose zero external pressure. We note that Eqs~\eqref{eq:p}--\eqref{Pe} \rev{represent the} non\rev{-}linear coupling between the eye and the contact lens models\rev{, since $p$ is non-linearly coupled to the shape of the eye $H$ in Eq~\eqref{eq:p} and $H$ is determined by $p$ in Eqs~\eqref{pressure_bc}--\eqref{Pe}}.

On $\Gamma_{in,v}$ and $\Gamma_{in,h}$, we impose symmetric boundary conditions on the displacement $\vector{u}$ and its derivatives as follows:
\begin{equation}\label{bc2}
\begin{aligned}
u_r = 0 & \mbox{ and } \deriv{}{u_z}{r}=0 & \mbox{ on } & \Gamma_{in,v},\\
u_z = 0 & \mbox{ and } \deriv{}{u_r}{z}=0 & \mbox{ on } & \Gamma_{in,h}.
\end{aligned}
\end{equation}

\subsection{Eye model weak formulation}
\label{sec:eye_weak}
The equilibrium \rev{equations} in the eye, Eq~\eqref{eq:S}, together with \rev{their} boundary conditions, Eqs~\eqref{pressure_bc}--\eqref{bc2}, are solved in the weak formulation detailed below. Given the space 
\begin{equation*}
V=\left\{ \vector{w} \in [H^1(\Omega)]^3 : w_r|_{r=0}=0,  w_{\theta}=0, w_z|_{z=h(\mathfrak{R}_{eye})}=0 \right\},
\end{equation*}
the weak formulation is to find $\vector{u}\in V$ that satisfies the boundary condition\rev{s} Eq\rev{s}~\eqref{pressure_bc}\rev{--\eqref{bc2}} and
\begin{equation}\label{weak}
\int_{\Omega} \left(\nabla \cdot \matrix{S}(\vector{u}) \right) \cdot \vector{w} \ d\Omega = 0 \quad \forall \vector{w}\in V.
\end{equation}
Using the definition of $\matrix{S}$ provided in Section~\ref{sec:eye_model} and the boundary conditions Eq\rev{s}~\eqref{pressure_bc}--\eqref{bc2}, the weak formulation Eq~\eqref{weak} can be expressed as:
\begin{equation}\label{eqn:weak_reformation}
\begin{aligned}
 \int_{\Omega} \lambda_{eye} \left( \nabla \cdot \vector{u}\right) \left(\nabla \cdot \vector{w} \right)  \ d\Omega  & + \  2  \int_{\Omega} \mu_{eye}\matrix{E}(\vector{u}) : \matrix{E}(\vector{w}) \ d\Omega \\
& + \ \int_{0}^{\mathfrak{R}_{eye}} P_{out}  \ (\vector{w} \cdot \vector{n})\Big{|}_{z=h(r)} \sqrt{1 + (h'(r))^2} \ r  \ dr = 0,
\end{aligned}
\end{equation}
where the tensor operation $\matrix{A} : \matrix{B} = \Sigma_{i,j} A_{ij}B_{ij}$. Note that $\lambda_{eye}$ and $\mu_{eye}$ are both proportional to the eye Young's modulus $E_{eye}$ via Eq~\eqref{eq:lame} and that $P_{out}$ is proportional to the lens Young's modulus $E_{lens}$ via Eq~\eqref{eq:p}. Therefore, if we divide Eq~\eqref{eqn:weak_reformation} by $E_{eye}$, then the governing equation for the displacement vector of the eye only depends on the ratio $\mathrm{E}=E_{lens}/E_{eye}$. We remind the reader that we are considering both a constant and a spatially\rev{-}dependent $E_{eye}$ in this work. 

\subsection{Numerical methods}
\label{sec:numerical_method}
The spatial discretization of the ocular domain $\Omega$ is handled via a triangular mesh refined on the ocular surface boundary $\Gamma_{out}$. We refine the mesh on $\Gamma_{out}$, where the suction pressure of the lens $p$ is applied on the eye, to better capture the effect of the variations of $p$ underneath the lens, and to \rev{better} capture the \rev{deformations} in the limbus region. A finer mesh is not needed in the part of the domain $\Omega$ which represents the inside of the eye, where the effect of the lens on the ocular domain is expected to be negligible. We solve the weak formulation in Eq~\eqref{eqn:weak_reformation} using $\mathbb{P}^2$ finite elements for the ocular displacement $\vector{u}$ implemented via the finite element library FreeFem++~\cite{hecht2012new}. After the displacement is computed, we solve weakly for the stresses (i.e., the components of the stress tensor $\matrix{S}$ in Eq~\eqref{eq:S}). Using $\mathbb{P}^2$ finite elements for $\vector{u}$ guarantees a linear/first-order approximation for the stresses.

In Ross et al.~\cite{ross2016existence}, the singular initial value problem given by Eqs~\eqref{eq:odeT}--\eqref{eq:bcT1} and governing the contact lens mechanics was shown to be well-posed. Specifically, it was shown that for every $\Sigma>0$, where $\Sigma=\frac{dR_{lens}}{dr}(0)$, there is a unique solution, $R_{lens}(r,\Sigma)$ and $T(r,\Sigma)$, of Eqs~\eqref{eq:odeT}--\eqref{eq:bcT1} defined on $r \in [0,\infty)$ that continuously depends on $\Sigma$ with $R\rev{_{lens}}(0,\Sigma)=T(0,\Sigma)=0$, $\frac{dR_{lens}}{dr}(0)=\Sigma$, and $\frac{dT}{dr}(0,\Sigma)=(1+\sigma_{lens})(\Sigma-1)$. Consequently, to approximate the boundary value problem, Eqs~\eqref{eq:odeT}--\eqref{eq:bcT2}, we implement a shooting method. The solution to the boundary value problem corresponds to the root of the following: 
\begin{equation}\label{eq:lens_root}
F(\Sigma)=T(\mathfrak{R}_{lens},\Sigma),
\end{equation}
where $T$ is the solution to the initial value problem given by Eqs~\eqref{eq:odeT}--\eqref{eq:bcT1} and $\frac{dR_{lens}}{dr}(0)=\Sigma$. \rev{Note that the root of Eq~\eqref{eq:lens_root} satisfies the lens radial tension boundary condition in Eq~\eqref{eq:bcT2}.} The bisection method is implemented to approximate the root of $F$. The initial value problem is approximated using a four-step Adams-Bashforth explicit method, with the radial grid spacing set to $0.001$~cm. To initialize, the classic Runge-Kutta method was \rev{used}. \rev{After solving Eqs~\eqref{eq:odeT}--\eqref{eq:bcT2}, we compute the suction pressure under lens $p$ via Eq~\eqref{eq:p} using a second order centered finite difference approximation.}

\rev{Our numerical algorithm to solve the coupled} eye and lens \rev{problems} \rev{is to alternate solving the eye problem and the lens problem, and to iterate} until convergence is reached, as detailed below. Let $n\geq0$ be the iteration index, and let $\phi^n$ be the approximation of a quantity $\phi$ at the n$^{th}$ iteration. As the initial conditions of the algorithm, for $n=0$, let $H^0=h$ and $p^0=0$; hence, $P_{out}^0=0$ and $\vector{\eta}^0=\vector{u}^0=\vector{0}$. Then, while $0<n\leq100$\rev{,}
\begin{enumerate}[1)]
\item Find the lens displacement $\vector{\eta}^n$ and suction pressure under the lens $p^n$ by solving Eqs~\eqref{eq:odeT}--\eqref{eq:p} using a smoothed $H^{n-1}$;
\item Find the external pressure $P_{out}^n$ using $p^n$ in Eq~\eqref{Pe};
\item Solve the weak formulation using $P_{out}^n$ in Eq~\eqref{eqn:weak_reformation} and find the ocular displacement $\vector{u}^n$;
\item Find the deformed ocular surface $H^n$ using $\vector{u}^n$ in Eq~\eqref{eq:H}; and
\item If 
$
\operatorname{max} \left(\dfrac{||\vector{u}^n-\vector{u}^{n-1} ||_{L^{\infty}(\Omega)}}{||\vector{u}^{n-1} ||_{L^{\infty}(\Omega)}}, \; \dfrac{||p^n-p^{n-1} ||_{{\infty}}}{||p^{n-1} ||_{\infty}} \right) < \varepsilon
$, then the algorithm has converge\rev{d}, otherwise let $n=n+1$ and go to Step 1.
\end{enumerate}
\rev{Note that the coupling of the lens and eye models is fully non-linear since the deformed shape of the eye $H$ is needed to determine the suction pressure under the lens $p$ in Eq~\eqref{eq:p}, and vice versa, the lens suction pressure is needed to determine the ocular deformation $\vector{u}$ in Eq~\eqref{eqn:weak_reformation}. In staggering the two numerical methods, we are linearizing the problem and such non-linearity. In Step 1 of the algorithm, since the lens suction pressure, defined by Eq~\eqref{eq:p}, depends only on the derivatives of the deformed ocular shape $H$, we approximate the shape of the eye from the previous iteration, $H^{n-1}$, using a cubic spline before taking its derivatives and conforming the contact lens to it.}  

The results presented in this paper are obtained for a value of $\varepsilon=10^{-6}$. The coupling algorithm is implemented in \textit{\rev{P}ython} and communicates with the FreeFem++ script by passing the suction pressure and running FreeFem++ in Step~3, and by passing the deformed ocular surface shape back to the main \rev{\textit{Python}} frame after Step~4. Depending on the eye and lens material parameters, for the values considered and listed in Table~\ref{tab:material}, the algorithm converge\rev{s} in 4--10 iterations.

\begin{figure}[h!]
\begin{center}
\begin{tikzpicture}
    \node at (0.5,0) {\includegraphics[width=.6\linewidth]{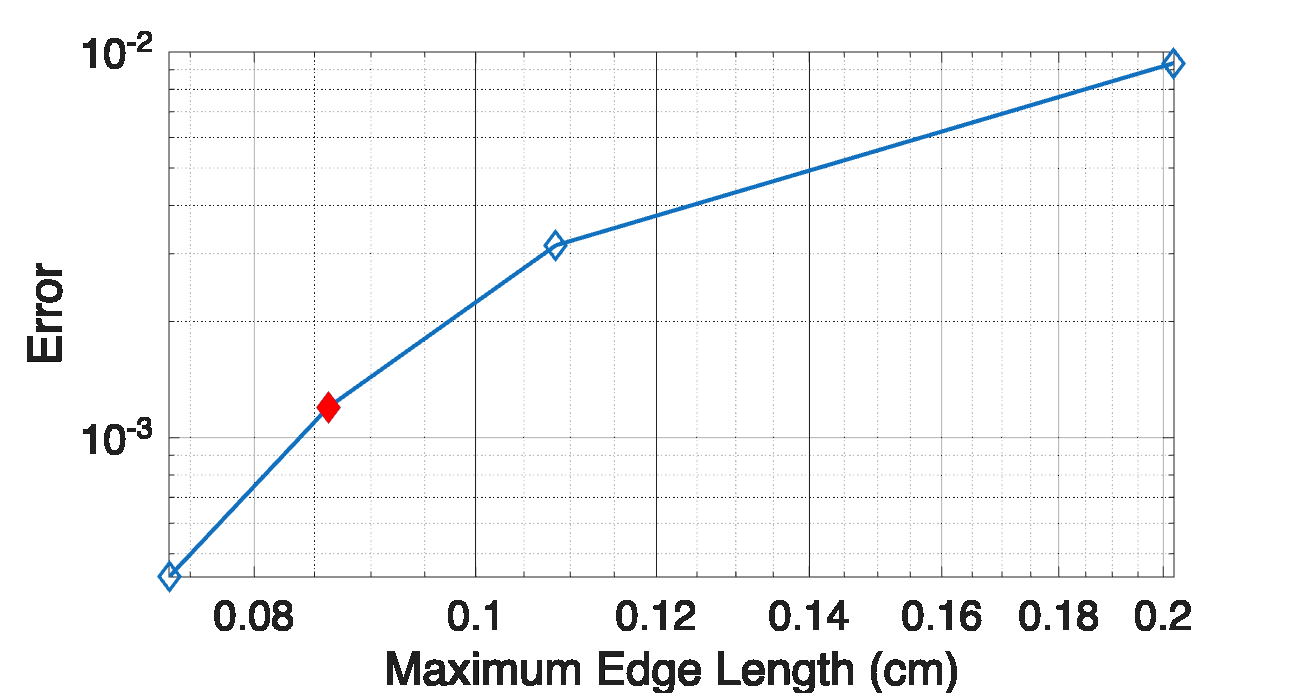}};
    \coordinate [label={left:\textcolor{black}{\bf \large A}}] (A) at (-4,3.2);    
    \node at (8.6,0) {\includegraphics[width=.5\linewidth]{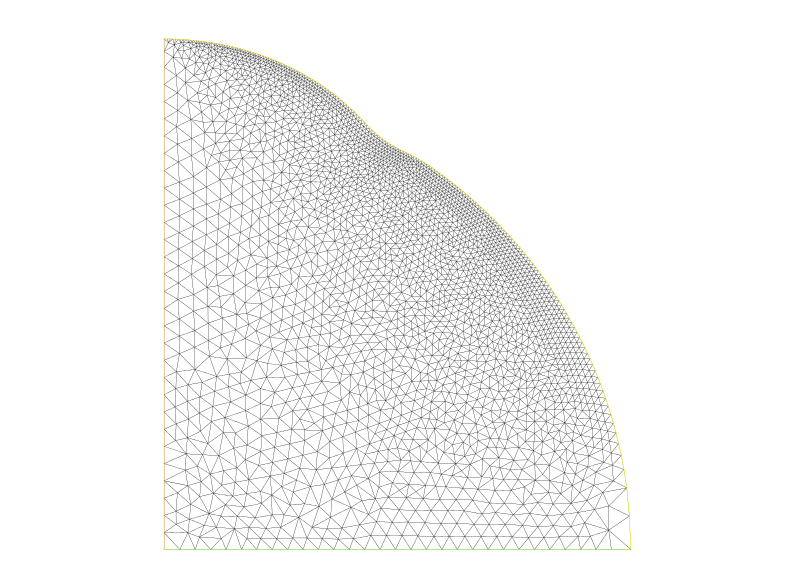}};
    \coordinate [label={left:\textcolor{black}{\bf \large B}}] (B) at (6,3.2);
\end{tikzpicture}
    \caption{(\textbf{A}) \rev{Results of the mesh convergence study performed.} (\textbf{B}) Mesh of the reference ocular domain $\Omega$ with 2806 vertices\rev{, corresponding to the red diamond mesh reported in (\textbf{A})}.}
    \label{fig:grid_refinement}
\end{center}    
\end{figure}

We perform\rev{ed} a grid refinement study of the finite element model used to approximate the displacements of the eye. We consider\rev{ed} a $100$~$\mu$m thick, average-shaped contact lens (Appendix, Section~\ref{sec:supp_lens}, Column 2 of Table~\ref{tab:lens_geometry}) inserted on an average-shaped eye (Appendix, Section~\ref{sec:supp_eye}, Column 2 of Table~\ref{tab:OcularParameter}) with the following material parameters: $E_{eye} = 0.2 $ MPa, $E_{lens} = 0.1$ MPa, and $\sigma_{eye} =\sigma_{lens} = 0.49$. The number of vertices in the mesh \rev{ranged} from 344 to 7697 vertices.  Figure~\ref{fig:grid_refinement}\textbf{A} shows the \rev{error of the grid refinement study as a function of the mesh size, measured as the maximum edge length of the mesh triangles. The error was computed as the maximum absolute difference between the norm of the displacement $||\vector{u}||$ on $\Gamma_{out}$ of a mesh and of the finest mesh considered (7697 vertices)}.  As the number of vertices increased, the solutions became closer to each other.  Given that the approximate solution on the grid \rev{with 2806 vertices} differed by \rev{approximately} $0.001$~$\mu$m to the solution with the finest grid considered \rev{(see red diamond in Figure~\ref{fig:grid_refinement}\textbf{A})}, we used a mesh grid with 2806 vertices for the remainder of the paper, see Figure~\ref{fig:grid_refinement}\textbf{B}. \rev{For a plot of the magnitude of the displacement on $\Gamma_{out}$ for the different meshes considered, see Figure~\ref{fig:grid_refinement2} in the Appendix, Section~\ref{sec:more_grid_ref}.}

\section{Results}
\label{sec:results}
First, we present our model results \rev{for} constant material parameters for the eye (i.e., \rev{for} a homogeneous eye model) in Section~\ref{sec:const_eye_res}. Then, in Section~\ref{sec:var_eye_res}, we show how a spatially-dependent eye Young's modulus (i.e., a heterogeneous eye model) affects the coupling between the eye and the lens.

\subsection{Homogeneous eye model}
\label{sec:const_eye_res}
The eye material parameters were set to be $E_{eye} = 0.2$ MPa and $\sigma_{eye} = 0.49$.  The contact lens was assumed to have a constant thickness $\tau=100$~$\mu$m with material parameters $E_{lens}=0.1$ MPa and $\sigma_{lens} = 0.49$.  In what follows, unless otherwise stated, the material parameters \rev{have} these values.
First, we present our model results on an average-shaped eye (Column~2 of Table~\ref{tab:OcularParameter} in the Appendix, Section~\ref{sec:supp_eye}) and an average-shaped contact lens (Column~2 of Table~\ref{tab:lens_geometry} in the Appendix, Section~\ref{sec:supp_lens}). 

Figure~\ref{fig:standard}({\bf A}) shows \rev{the} predicted contact lens deformations, ({\bf B}) shows \rev{the} contact lens suction pressure, ({\bf C}) shows \rev{the} deformations \rev{of the ocular surface}, and ({\bf D}--{\bf F}) show \rev{the} plots of the eye deformations in the reference coordinates $(r,z)$. 
In Figure~\ref{fig:standard}{\bf A}, the lens radial displacement, $\eta_r$, is the solid line; the lens vertical displacement, $\eta_z$, is the dashed line; the ocular normal displacement, $\eta_n=\vector{\eta} \cdot\vector{n}$, where $\vector{n}$ is the outward normal direction to the \rev{posterior surface of the lens}, is the dotted line; and the norm of the lens displacement vector, $|| \vector{\eta} ||$, is the dashed-dotted line. Recall that \rev{we formulated the model so that} the contact lens conform\rev{s} to the ocular surface.  Thus, the reference point $(r,g(r))$ on \rev{the} posterior surface of the contact lens is mapped to $(R_{lens}(r),H(R_{lens}(r))$ on the deformed ocular surface. The vertical contact lens displacements $\eta_z$ are negative except at the edge of the contact lens, where the contact lens conforms to the sclera\rev{,} thus resulting in positive displacements. Similarly, the radial displacements $\eta_r$ are negative \rev{at} the end of the cornea and \rev{in} the limbal region\rev{,} thus resulting in lens compression. Near the edge of the lens, the lens is radially stretched. Figure~\ref{fig:standard}{\bf A} shows that the maximum lens displacement \rev{occurs} in the limbal region, as shown by $||\vector{\eta}||$. The normal lens displacement $\eta_n$ varies in the radial direction $r$ similarly to the vertical \rev{d}isplacement $\eta_z$. The resulting contact lens suction pressure $p$ is shown in Figure~\ref{fig:standard}{\bf B}. The contact lens pushes down (positive~$p$) on the ocular surface in the center (\rev{the} corneal region) and \rev{at the} edge of the lens. \rev{In contrast}, in the limbal region, the contact lens pulls up (negative~$p$) on the ocular surface. \rev{The discontinuity in the suction pressure around $r=0.7$~cm is due to the assumption that the external pressure experienced by the eye is zero outside of the lens, see Eq~\eqref{Pe}.}
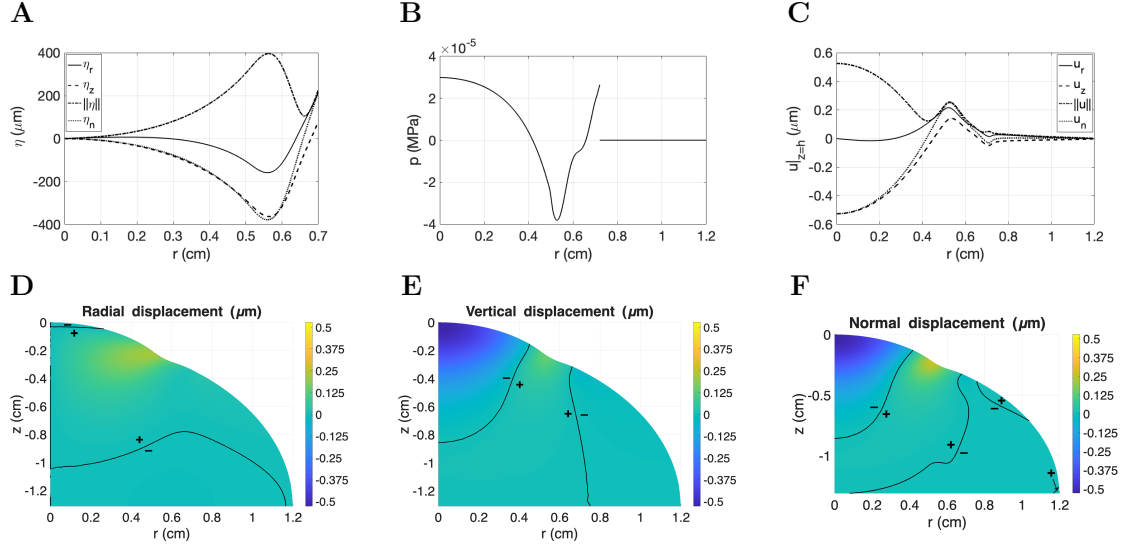
\begin{figure}[t]
\begin{center}
\scalebox{0.9}{
\begin{tikzpicture}
\node at (-.2,0) {\includegraphics[width=0.33\linewidth]{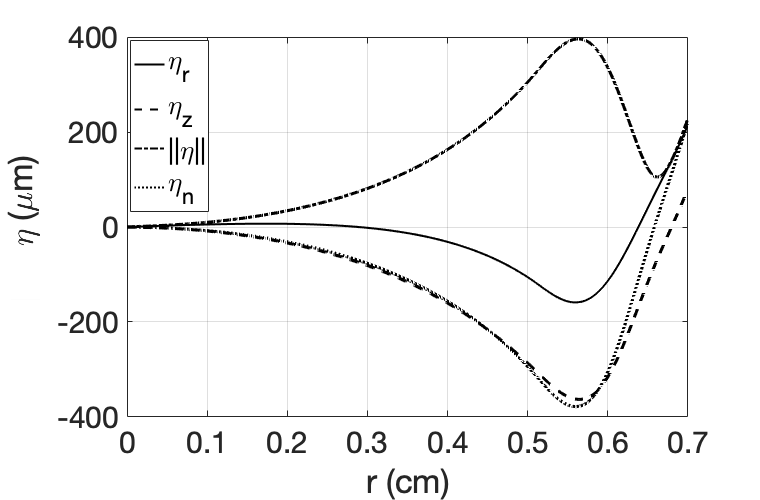}};
\coordinate [label={left:\textcolor{black}{\bf \large A}}] (A) at (-2.2,2);
\node at (5.5,0) {\includegraphics[width=0.33\textwidth]{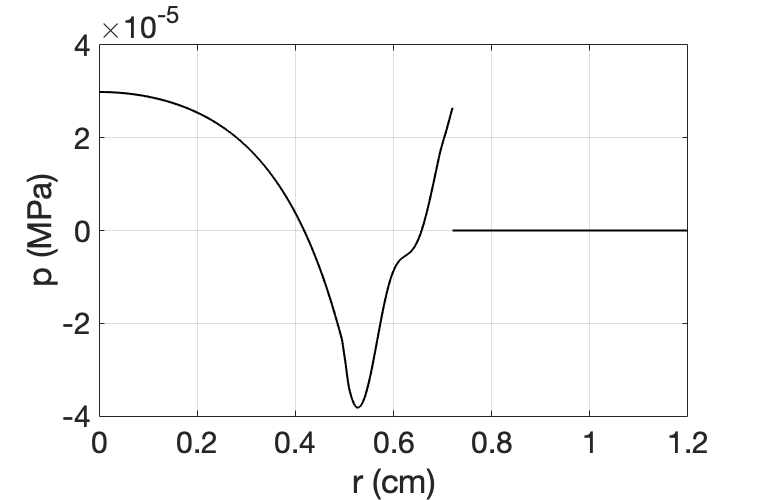}};
\coordinate [label={left:\textcolor{black}{\bf \large B}}] (B) at (3.5,2);
\node at (11.2,0) {\includegraphics[width=0.33\textwidth]{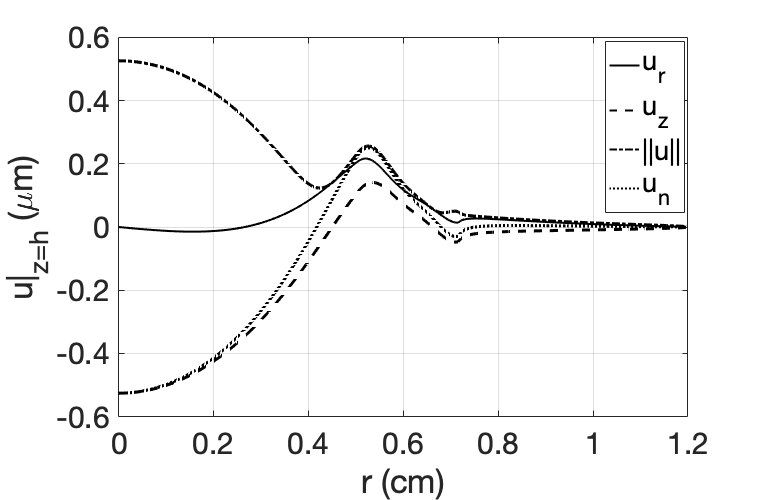}};
\coordinate [label={left:\textcolor{black}{\bf \large C}}] (C) at (9.2,2);
\node at (-.2,-4) {\includegraphics[width=0.33\linewidth]{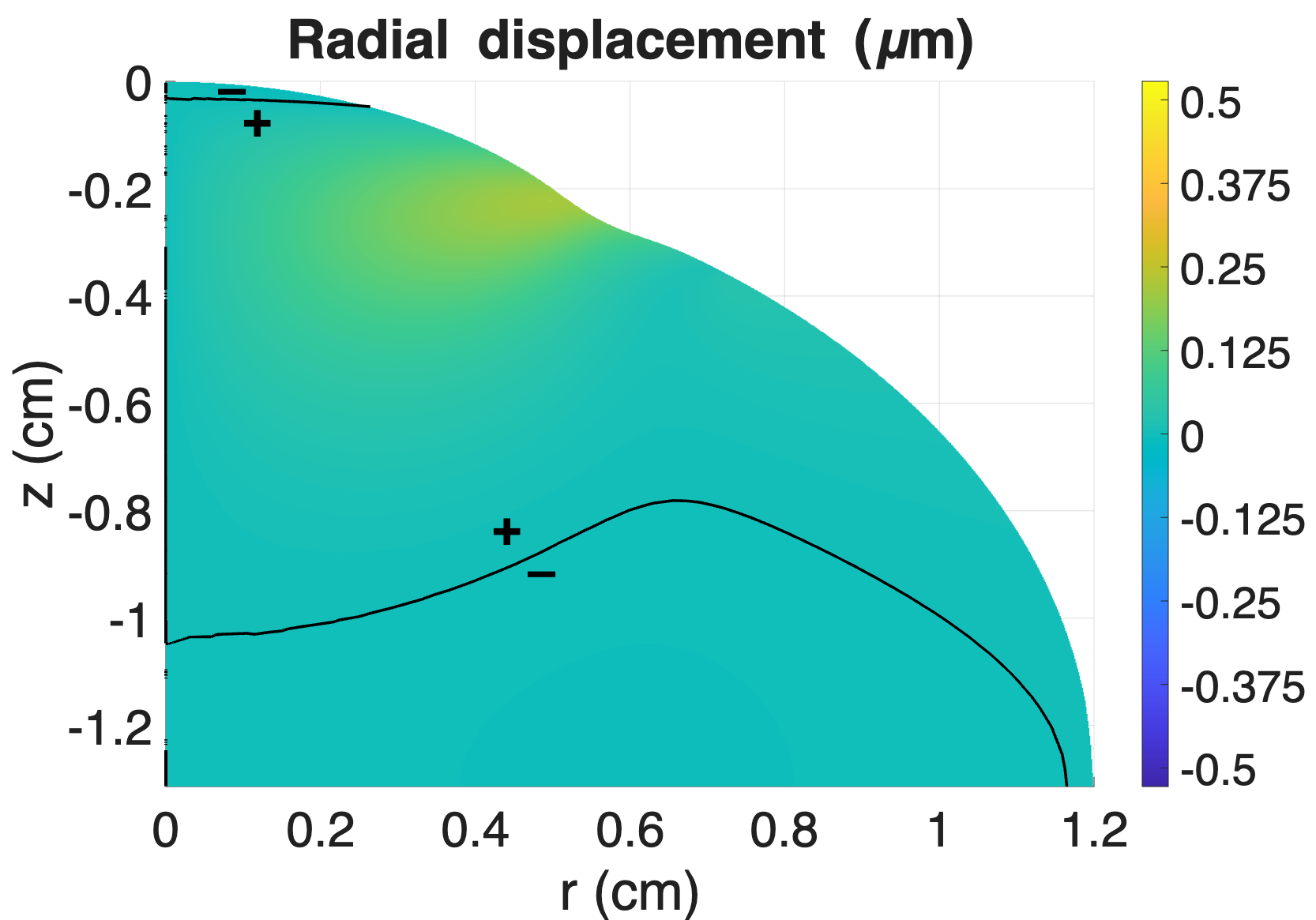}};
\coordinate [label={left:\textcolor{black}{\bf \large D}}] (D) at (-2.2,-2);
\node at (5.5,-4) {\includegraphics[width=0.33\linewidth]{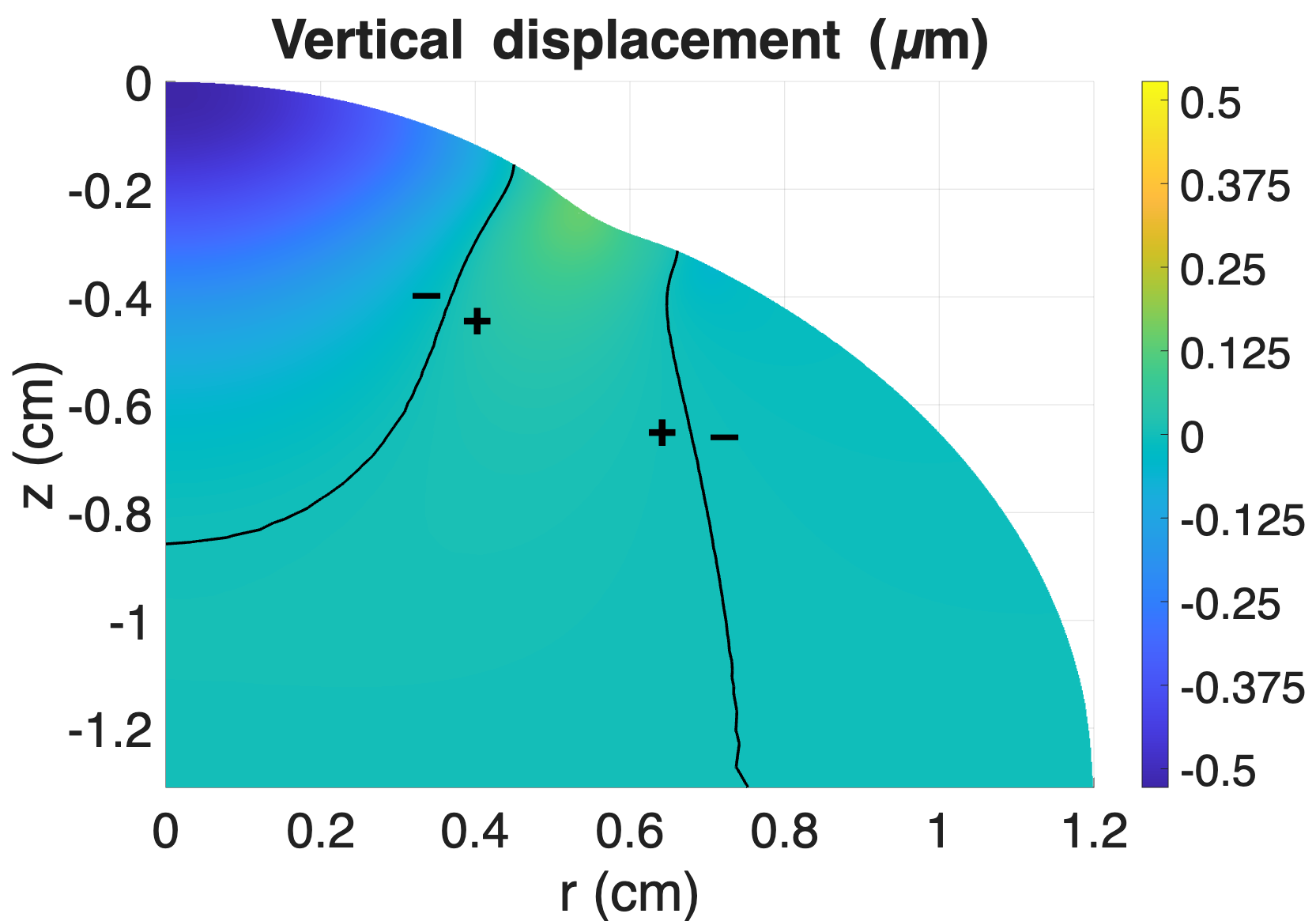}};
\coordinate [label={left:\textcolor{black}{\bf \large E}}] (E) at (3.5,-2);
\node at (11.2,-4) {\includegraphics[width=0.33\linewidth]{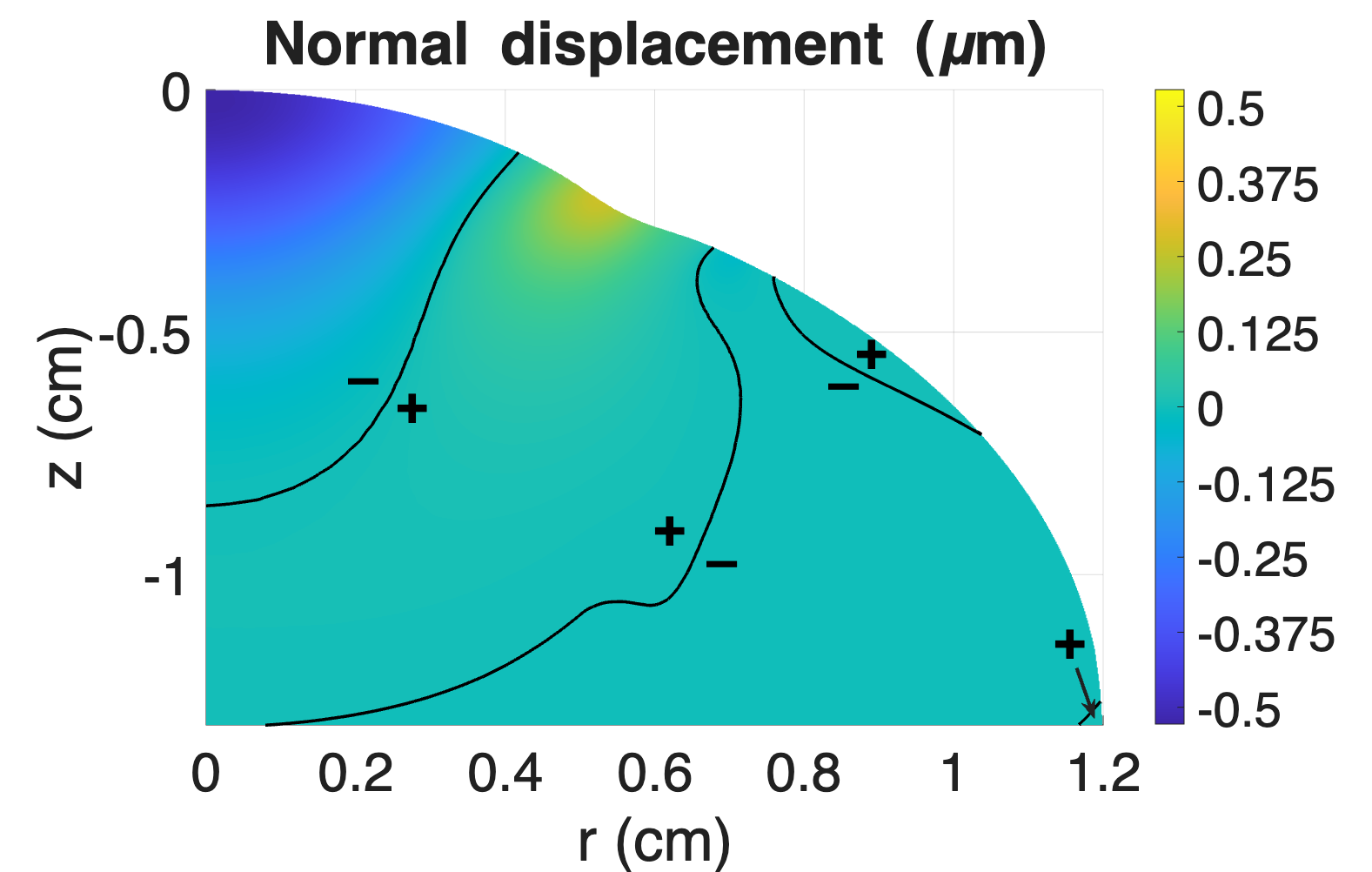}};
\coordinate [label={left:\textcolor{black}{\bf \large F}}] (F) at (9.2,-2);
\end{tikzpicture}
}
\end{center}
\caption{\textit{Homogeneous average-shaped eye and average-shaped contact lens of constant thickness.} ({\bf A}) Contact lens displacements. ({\bf B}) Contact lens suction pressure, see Eq~($\ref{eq:p}$). ({\bf C}) Ocular surface displacements. ({\bf D}) The radial eye displacements. ({\bf E}) The vertical eye displacements. ({\bf F})~Eye displacement in the outward normal direction to the ocular surface. The solid black lines in the contour plots {\bf D}--{\bf F} are the zero level curves \rev{and the positive or negative signs indicate the sign of the displacement near the level curves}.}
\label{fig:standard}
\end{figure}
 
The ocular displacements on the surface $\Gamma_{out}$ are shown in  Figure~\ref{fig:standard}{\bf C}, where the ocular radial displacement, $u_r|_{z=h}$, is the solid line; the ocular vertical displacement, $u_z|_{z=h}$, is the dashed line; the ocular normal displacement, $u_n=\vector{u}|_{z=h} \cdot\vector{n}$, where $\vector{n}$ is the outward normal to $\Gamma_{out}$, is the dotted line; and the norm of the ocular displacement vector, $|| \vector{u}|_{z=h} ||$, is the dashed-dotted line.  In the region close to the center of the eye ($0 \leq r \leq 0.4$~cm), where the lens suction pressure is positive, the ocular surface is compressed inward by at most $0.53$~$\mu$m.  Close to the edge of the lens \rev{($r=0.7$~cm)}, where the lens suction pressure is also positive, the ocular surface is compressed slightly inward (at most $0.048$~$\mu$m). Near the limbal region ($0.4<r<0.65$~cm), where the lens suction pressure is negative, since the contact lens must conform to the concave inward ocular surface of the limbus, the ocular surface extends outward by at most $0.26$~$\mu$m.   

The eye deformations in the reference domain $\Omega$ are shown in Figures~\ref{fig:standard}{\bf D}--{\bf F}. The radial and vertical displacements are plotted in Figures~\ref{fig:standard}{\bf D}--{\bf E}. The solid black lines are the level curves when the displacement is zero \rev{and the positive or negative signs indicate the sign of the displacement near the level curves}. In general, the displacements decrease in magnitude as you move inward into the eye and away from the ocular surface. Figure~\ref{fig:standard}{\bf \rev{F}} plots the outward normal (of the reference ocular surface) displacement. \rev{To compute the outward normal displacement for a point in the interior of the ocular domain $\Omega$, we consider the line that connects the point to the bottom left corner of the ocular domain and use the outward unit normal at the intersection between the line and $\Gamma_{out}$.} The magnitude of the displacement is less than $0.6$~$\mu$m, and the largest displacements occur at the top center of the eye and in the limbal region, where the eye experiences the maximum and minimum suction pressure values. 

\begin{figure}[h!]
\begin{center}
\scalebox{0.93}{
\begin{tikzpicture}
\node at (-.2,0) {\includegraphics[width=0.4\linewidth]{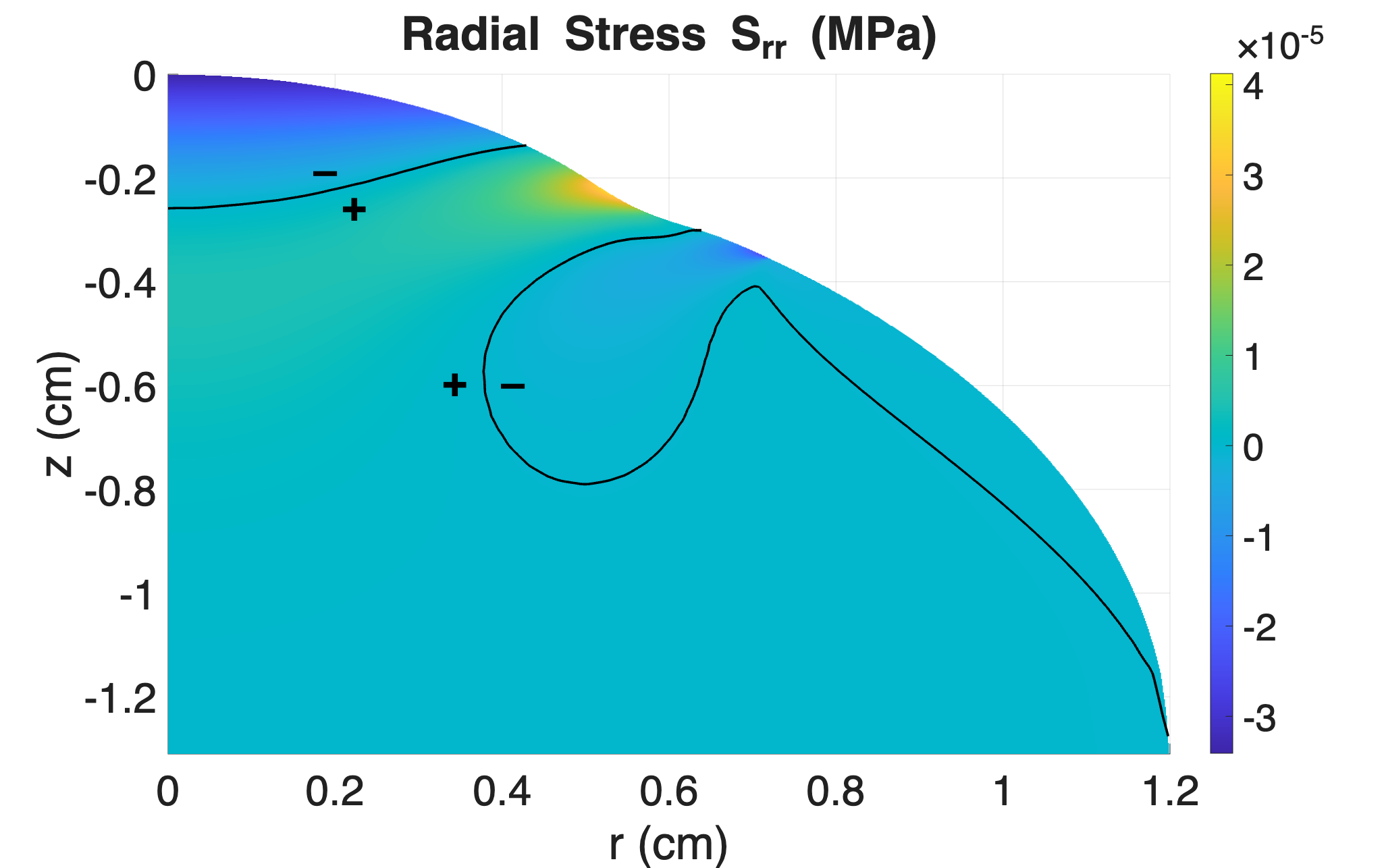}};
\coordinate [label={left:\textcolor{black}{\bf \large A}}] (A) at (-2.8,2.5);
\node at (7,0) {\includegraphics[width=0.4\textwidth]{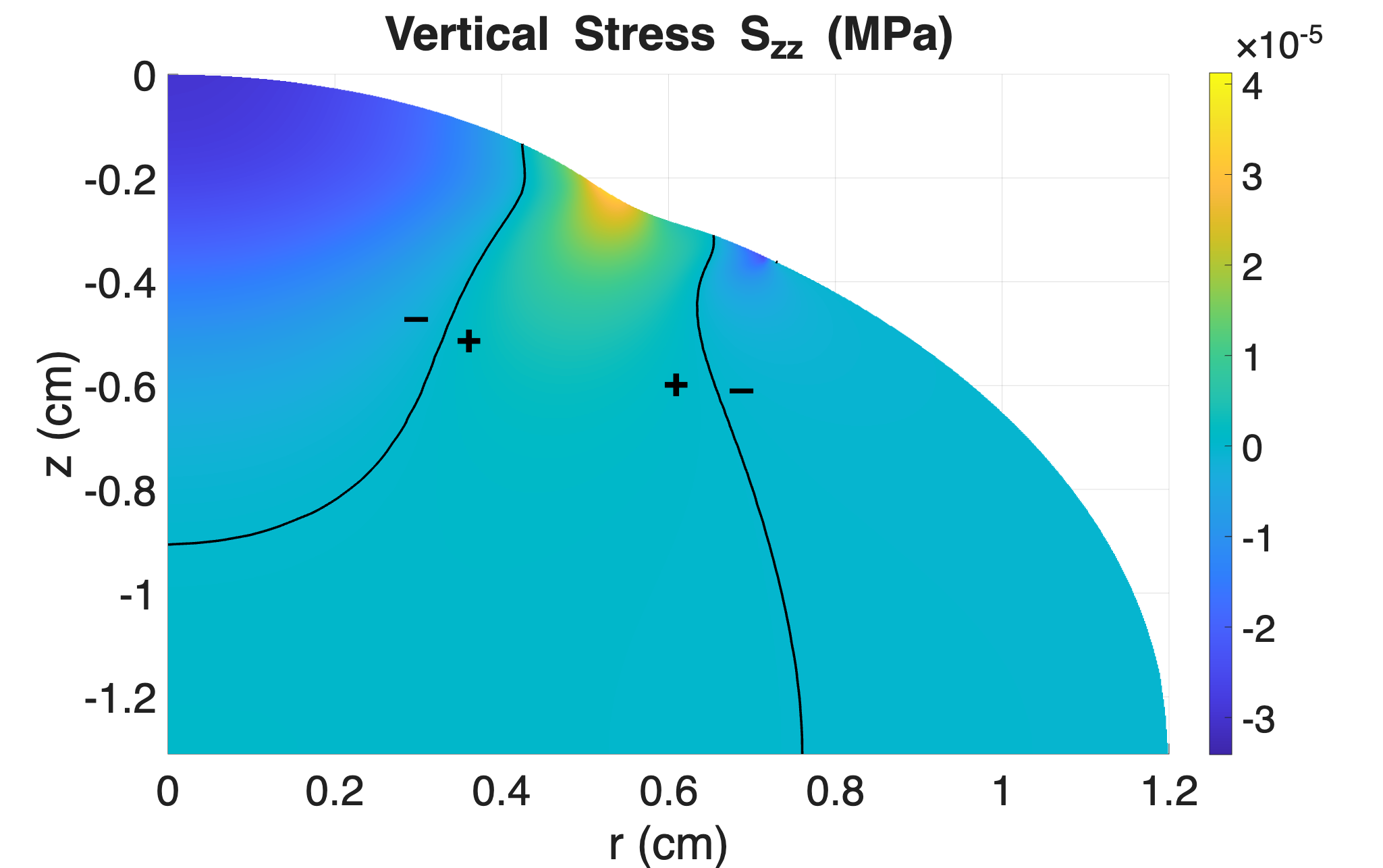}};
\coordinate [label={left:\textcolor{black}{\bf \large B}}] (B) at (4.4,2.5);
\node at (-.2,-5) {\includegraphics[width=0.4\textwidth]{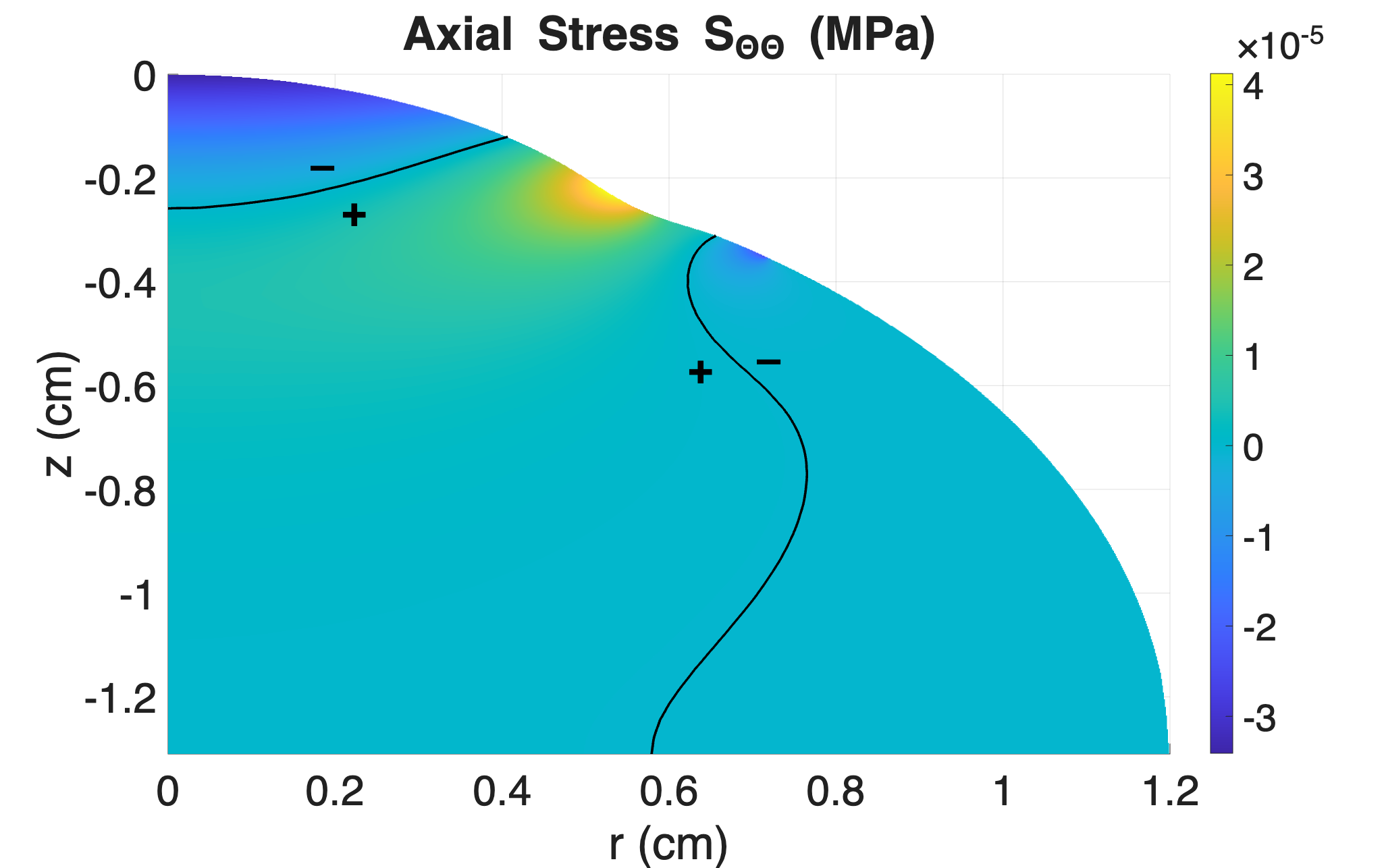}};
\coordinate [label={left:\textcolor{black}{\bf \large C}}] (C) at (-2.8,-2.5);
\node at (7,-5) {\includegraphics[width=0.4\linewidth]{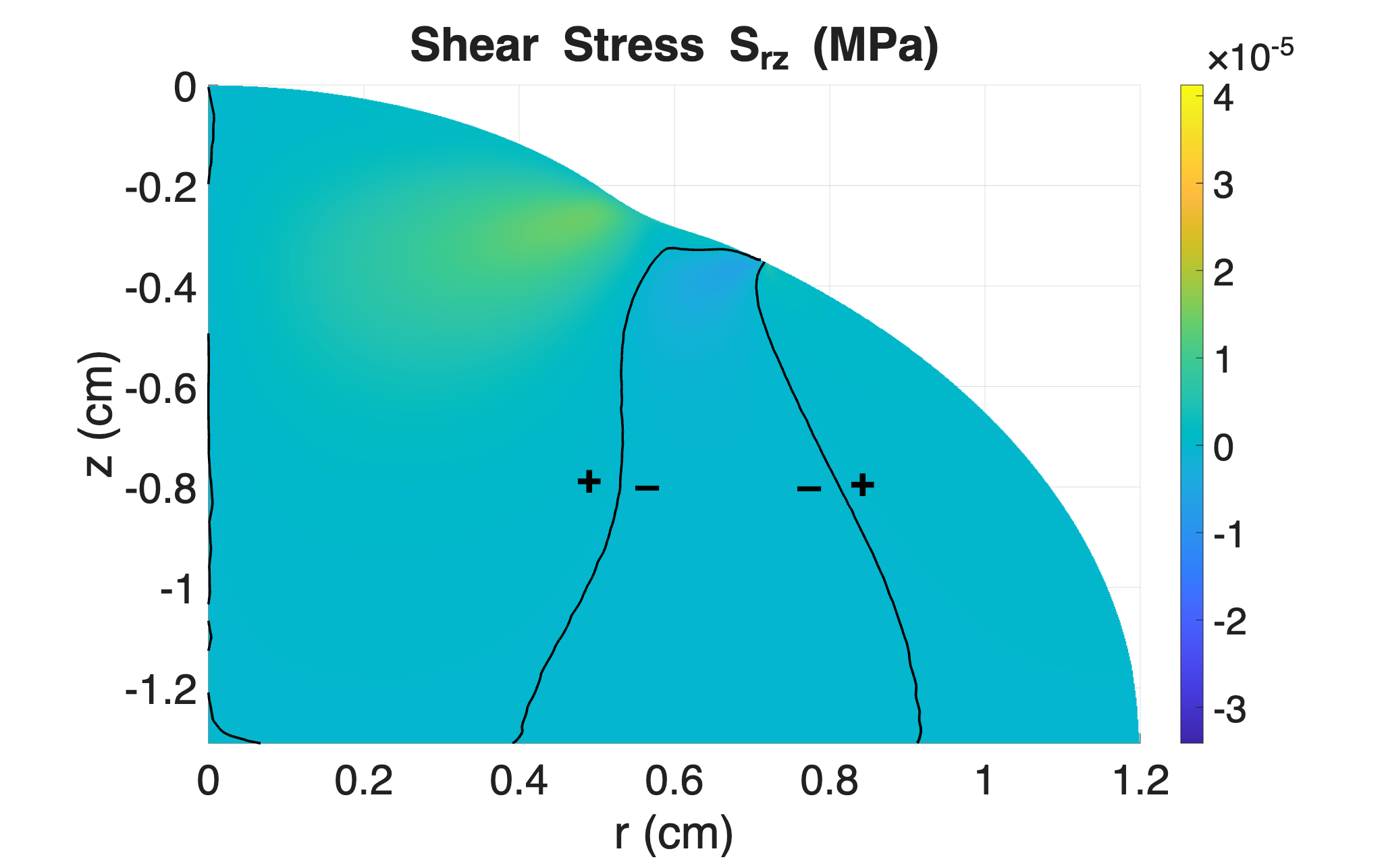}};
\coordinate [label={left:\textcolor{black}{\bf \large D}}] (D) at (4.4,-2.5);
\node at (8,-10) {\includegraphics[width=0.55\linewidth]{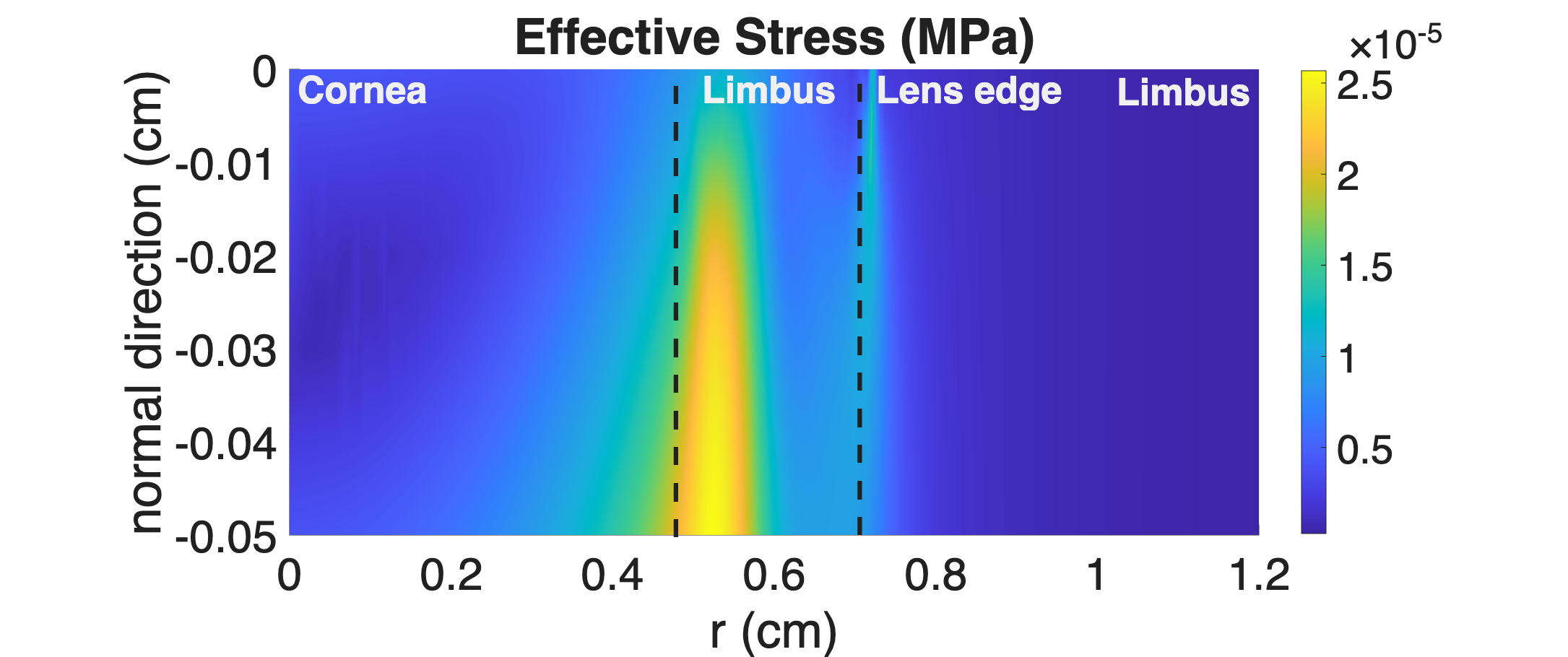}};
\node at (-.5,-10) {\includegraphics[width=0.49\linewidth]{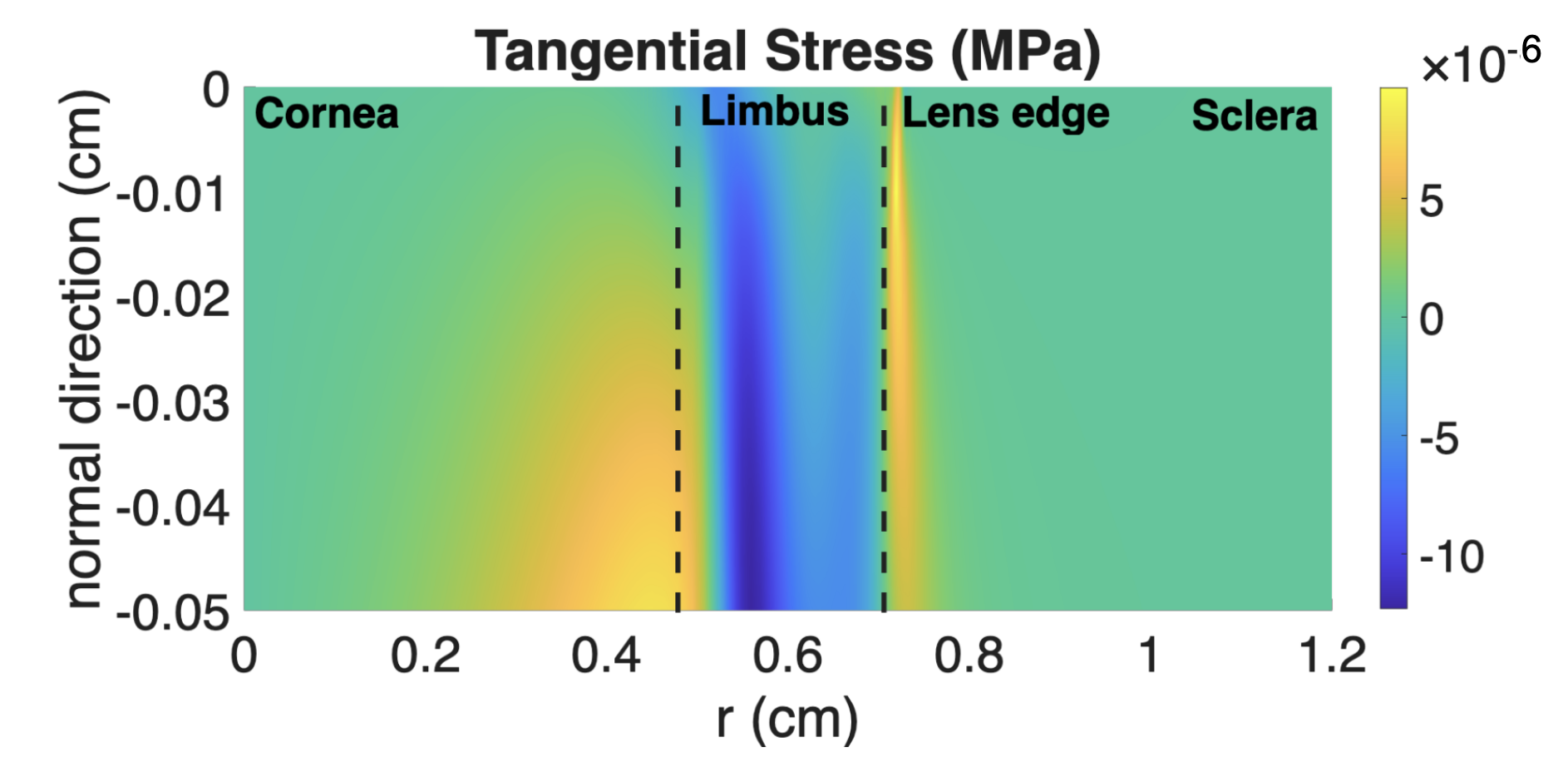}};
\coordinate [label={left:\textcolor{black}{\bf \large E}}] (E) at (-3.6,-7.8);
\coordinate [label={left:\textcolor{black}{\bf \large F}}] (F) at (5,-7.8);
\end{tikzpicture}
}
\end{center}
\caption{\textit{Homogeneous average-shaped eye and average-shaped contact lens of constant thickness.} Ocular stresses:
({\bf A}) radial stress $S_{rr}$; ({\bf B}) vertical stress $S_{zz}$; ({\bf C}) axial stress $S_{\theta \theta}$; and ({\bf D}) shear stress $S_{rz}$. ({\bf E}) The tangential stress $S_{tan}$ and ({\bf F}) the effective stress $S_{eff}$ are shown in the strip of tissue normal to the reference ocular surface $\Gamma_{out}$ (at \rev{zero} normal direction \rev{on the vertical axis}) of thickness $0.05$~cm. Note that the ranges of the color bars are the same in {\bf A}--{\bf D}, but are different in {\bf E} and {\bf F}. The solid black lines in the contour plots {\bf A}--{\bf D} are the zero level curves \rev{and the positive or negative signs indicate the sign of the stress near the level curves}.}
\label{fig:stresses}
\end{figure}
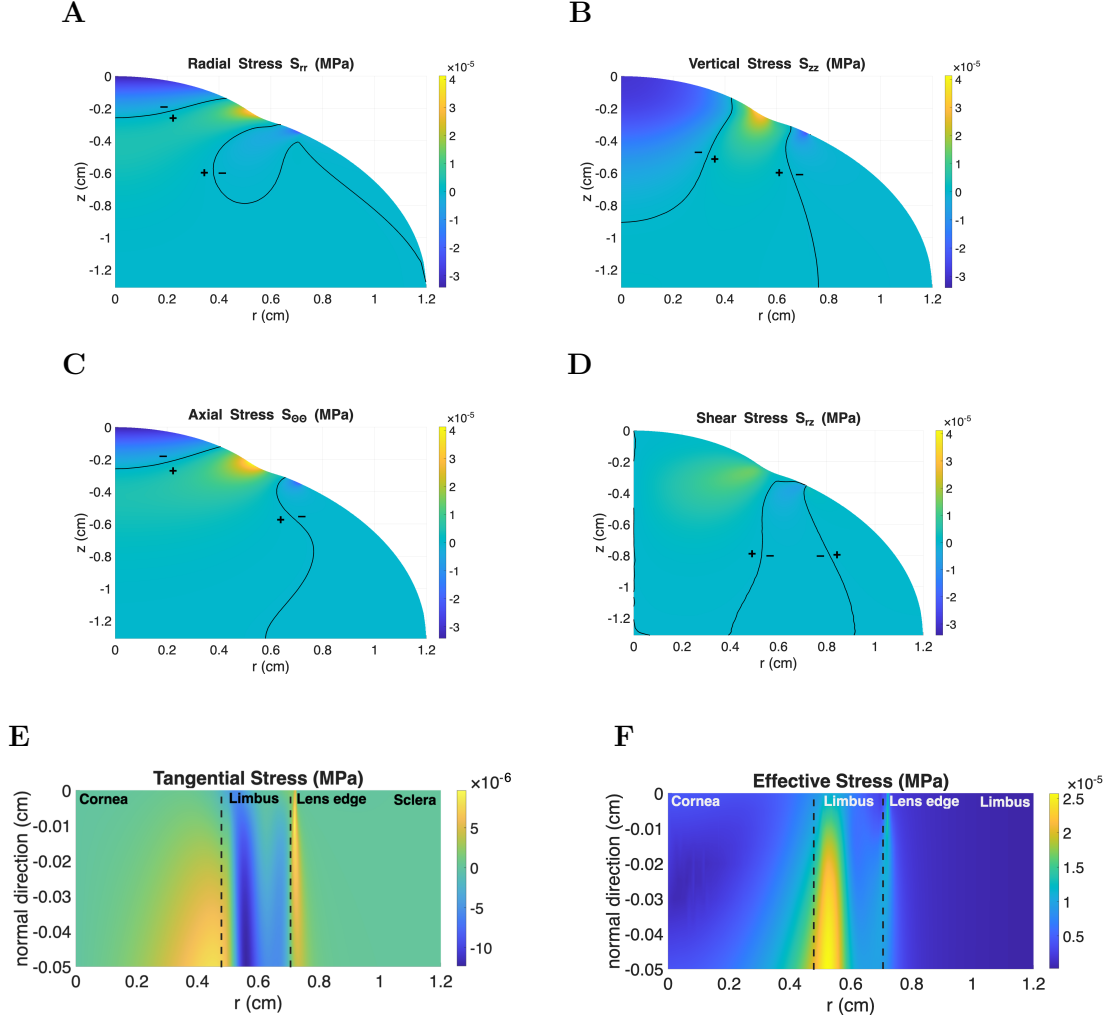

Figure~\ref{fig:stresses} displays the different components of the predicted stresses within the eye due to contact lens wear: ({\bf A}) radial stress $S_{rr}$; ({\bf B}) vertical stress $S_{zz}$; ({\bf C}) axial stress $S_{\theta \theta}$; and ({\bf D}) shear stress $S_{rz}$. The solid black lines are the level curves when the stress is zero. The magnitudes of the stresses are less than $4 \times 10^{-5}$ MPa, with the stresses \rev{of the greatest magnitude} occur\rev{ing} near or at the ocular surface. 

Figures~\ref{fig:stresses}{\bf E}--{\bf F} plot the tangential and effective stresses in \rev{the outer} strip of ocular tissue \rev{shown in Figure~\ref{fig:simplified_varied_ym}}, respectively. The strip of tissue is normal to the reference ocular surface $\Gamma_{out}$ (at \rev{zero} normal direction \rev{on the vertical axis}) and is of a thickness of $0.05$~cm. 
Tangential stress is computed as $S_{tan}=\vector{t} \cdot \matrix{S} \cdot \vector{n}$, where $\vector{n}$ and $\vector{t}$ are the normal and tangent vectors to the ocular reference surface $\Gamma_{out}$, respectively.  
The effective, or von Mises, stress is computed as follows\cite{Jones2009,carichino2012effect}:
\begin{equation}\label{Eq:Effective Stress}
    S_{eff}\left(\matrix{S}\right) = \sqrt{\dfrac{3}{2}\left(\left(\matrix{S}-\dfrac{1}{3}\tr{\matrix{S}}\matrix{I}\right)\left(\matrix{S}-\dfrac{1}{3}\tr{\matrix{S}}\matrix{I}\right)^T\right)}.
\end{equation} 
\rev{The ocular tissue is subject to a complex load due to contact lens wear and the effective
stress quantifies the risk of 
the ocular tissue to yield or fracture, thus experiencing permanent deformation due to such a load~\cite{Jones2009}.}
At the ocular surface (zero normal direction \rev{on the vertical axis}), the regions of nonzero tangential stress occur in the limbus, and these regions of nonzero tangential stress widen to include parts of the cornea as you move inside the tissue in the normal direction. Similarly, the regions of nonzero effective stress are at the edges of the limbus and widen to encompass more of the limbus as you move into the tissue in the normal direction. Differently, the regions of the nonzero effective stresses do not significantly penetrate into the corneal region.

\subsubsection{Influence of material parameters}
As we mentioned in Section~\ref{sec:eye_model}, the governing equation for the displacement vector of the eye\rev{, Eq~\eqref{eqn:weak_reformation},} only depends on the ratio $\mathrm{E}=E_{lens}/E_{eye}$. We explore the material parameter space $\mathrm{E}=E_{lens}/E_{eye}$ by fixing $E_{eye} = 0.2$ MPa and varying $E_{lens}$, values given in Table~\ref{tab:E_exploration}.  The values of $\mathrm{E}$ between $0.1$ and $10$ correspond to soft contact lenses, while $\mathrm{E}=25$ and $50$ mimic the interaction between the eye and stiffer contact lenses. As $\mathrm{E}$ increases, the contact lens is stiffer than the eye; thus, the eye will show larger deformations.  As the ocular surface deforms, its shape will change and, consequently, the surface the contact lens must conform to.  
\begin{table}[htb]
    \centering
    \caption{Young's modulus of the contact lens used to explore the $\mathrm{E}$ parameter space. $E_{eye}$ is kept constant at $0.2$~MPa.}
    \begin{tabular}{llllllllll}
        \toprule
        & \multicolumn{6}{c}{Soft contact lenses} & \multicolumn{3}{r}{Stiffer contact lenses} \\
        \toprule
        $E_{lens}$ [MPa] & 0.1 & 0.2 &  0.4 &  0.8 &  1.6 &  2.0  & \; &  5.0 & 10.0 \\   
       \midrule        
       $\mathrm{E}\rev{=E_{lens}/E_{eye}}$ [-]  & $\frac{1}{2}$ & 1 & 2 & 4 & 8 & 10  & \; & 25 & 50\\
       \bottomrule
    \end{tabular}
    
    \label{tab:E_exploration}
\end{table}

\begin{figure}[h!]
    \centering
    \scalebox{0.9}{
\begin{tikzpicture}
    \node at (0,0) {\includegraphics[width=.5\linewidth]{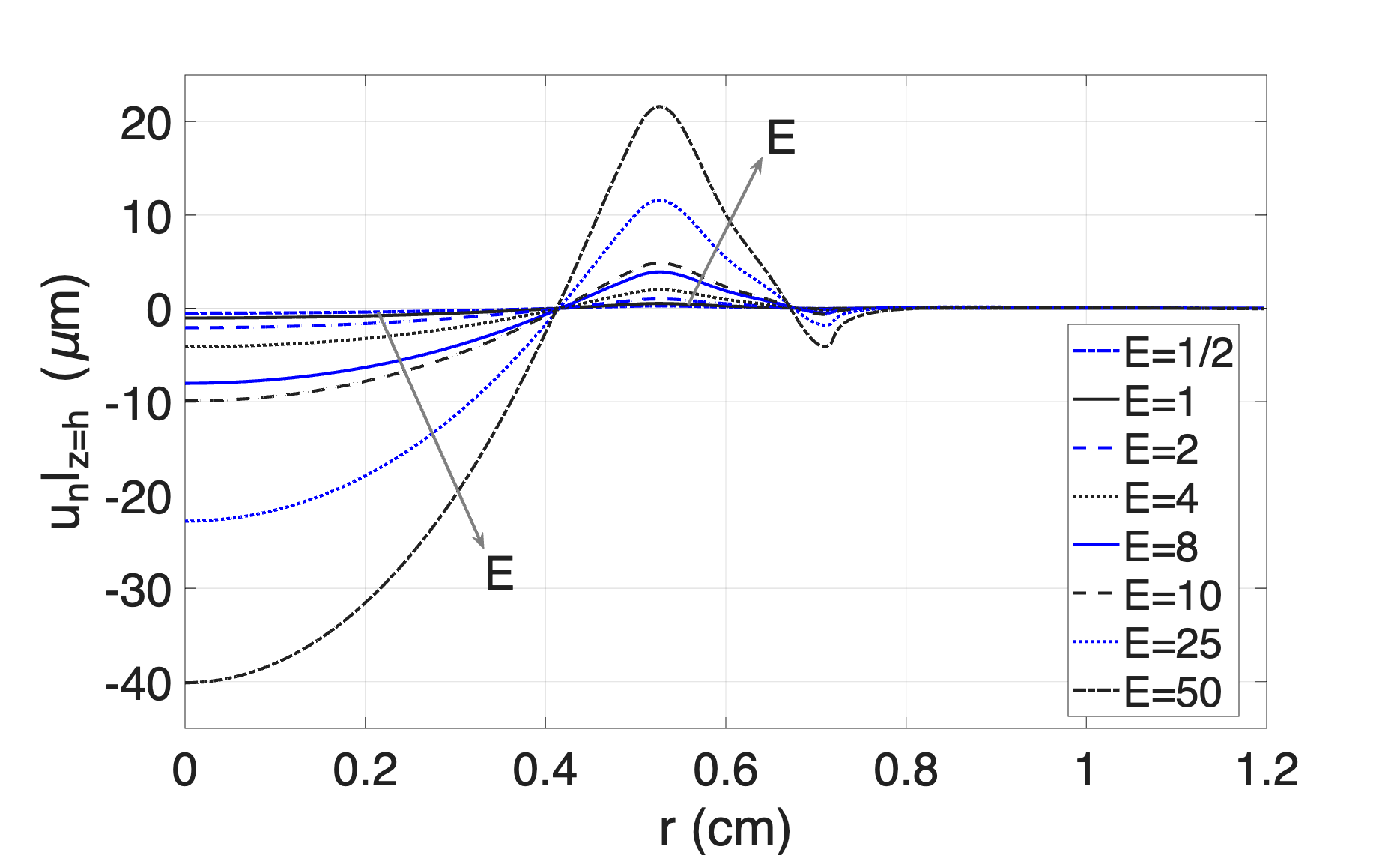}};
    \coordinate [label={left:\textcolor{black}{\bf \large A}}] (A) at (-3.5,3.15);
    \node at (9,0) {\includegraphics[width=.5\linewidth]{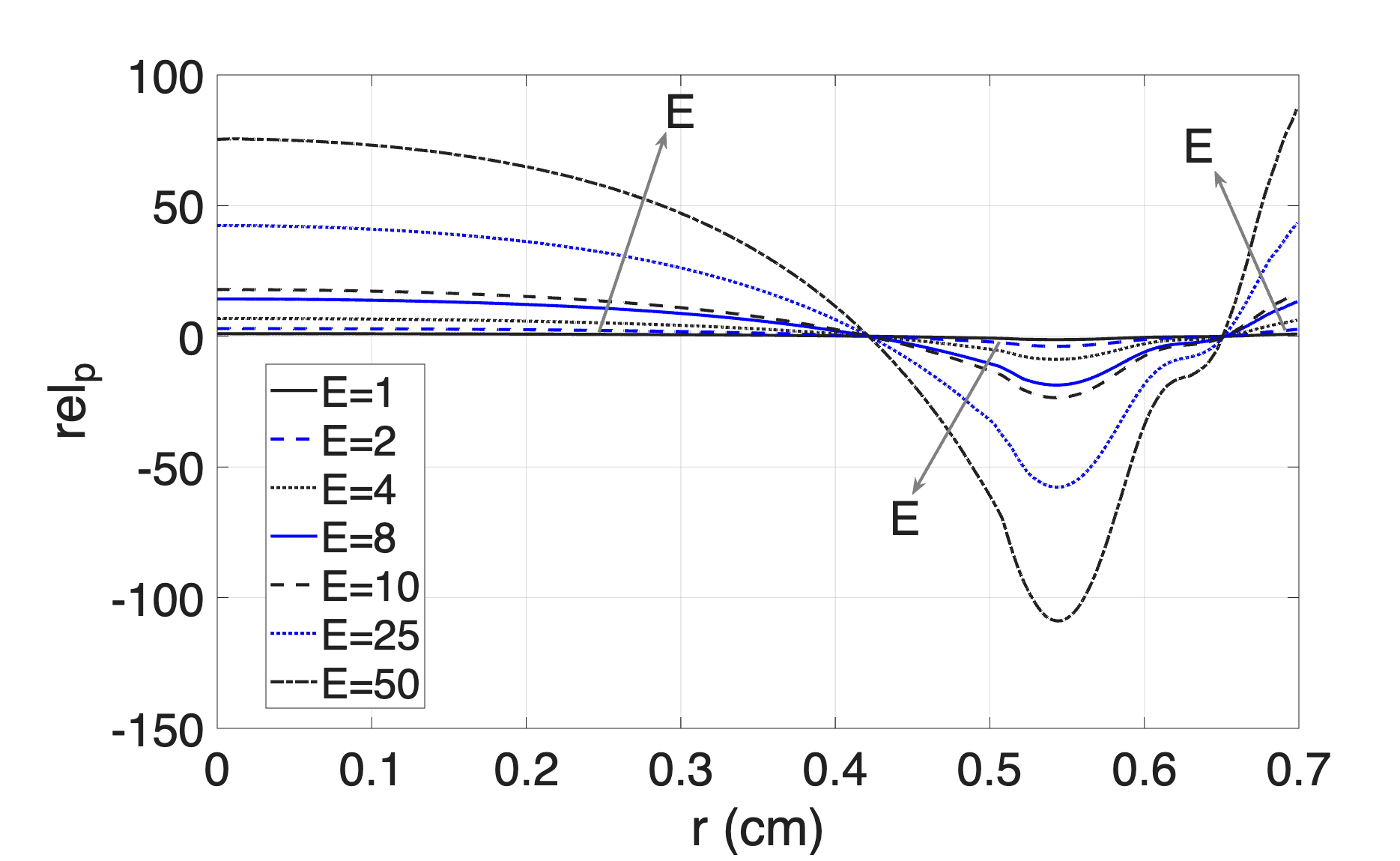}};
    \coordinate [label={left:\textcolor{black}{\bf \large B}}] (B) at (5.5,3.15);
    \node at (4,-6.75) {\includegraphics[width=.5\linewidth]{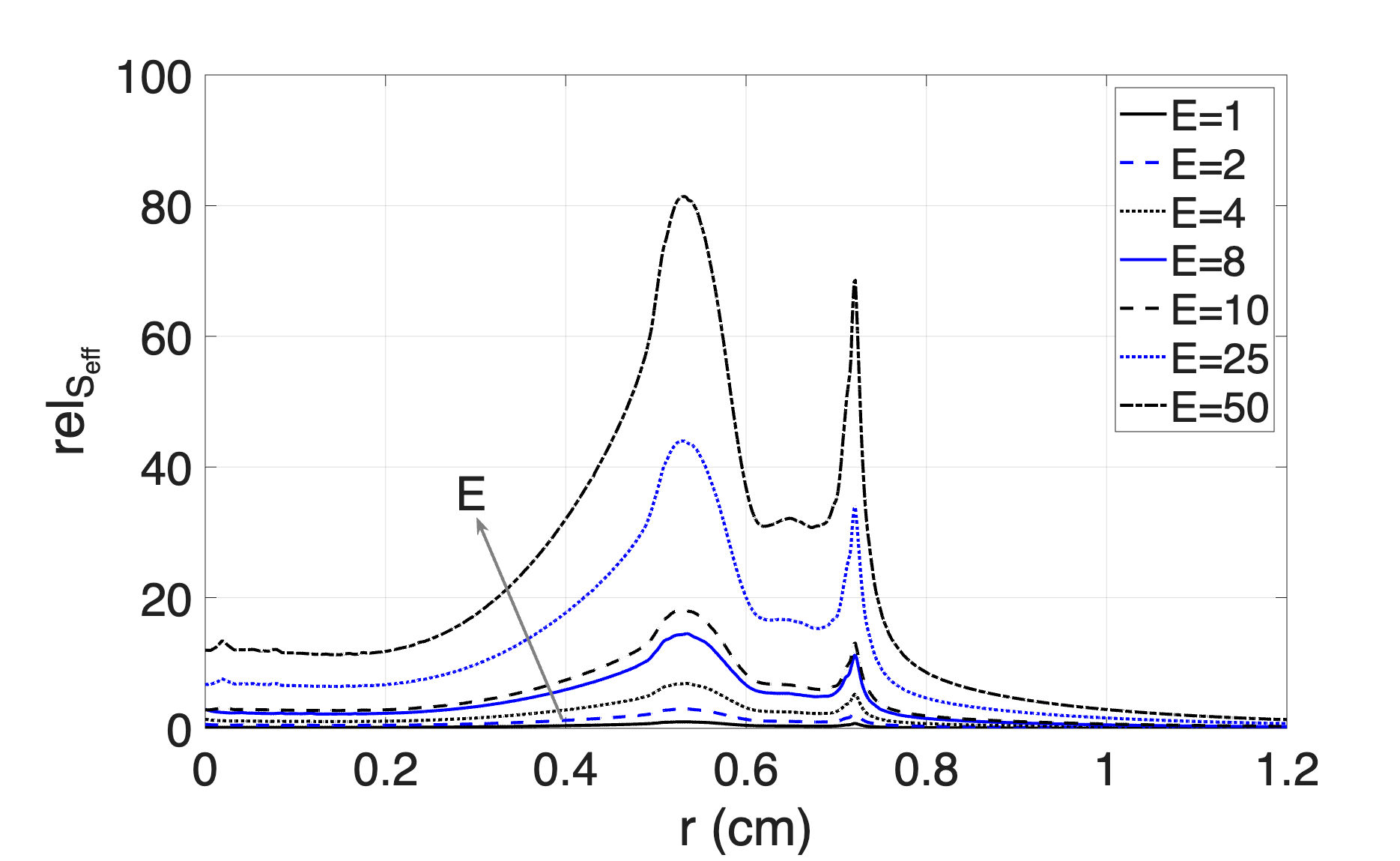}};
    \coordinate [label={left:\textcolor{black}{\bf \large C}}] (C) at (0.35,-3.6);
\end{tikzpicture}}
    \caption{\textit{Homogeneous average-shaped eye and average-shaped contact lens of constant thickness.} ({\bf A}) The ocular normal displacement, $u_n$, as function of the radial reference coordinate for different values of $\mathrm{E} = E_{lens}/E_{eye}$.  ({\bf B}) The relative difference of suction pressures, $rel_p$, as a function of the radial reference coordinate for different values of $\mathrm{E}$. ({\bf C}) The relative difference of effective stress, $rel_{S_{eff}}$, as a function of the radial reference coordinate for different values of $\mathrm{E}$. \rev{The arrows indicate the direction which $\mathrm{E}$ increases.}}
    \label{fig:material_props}
\end{figure}
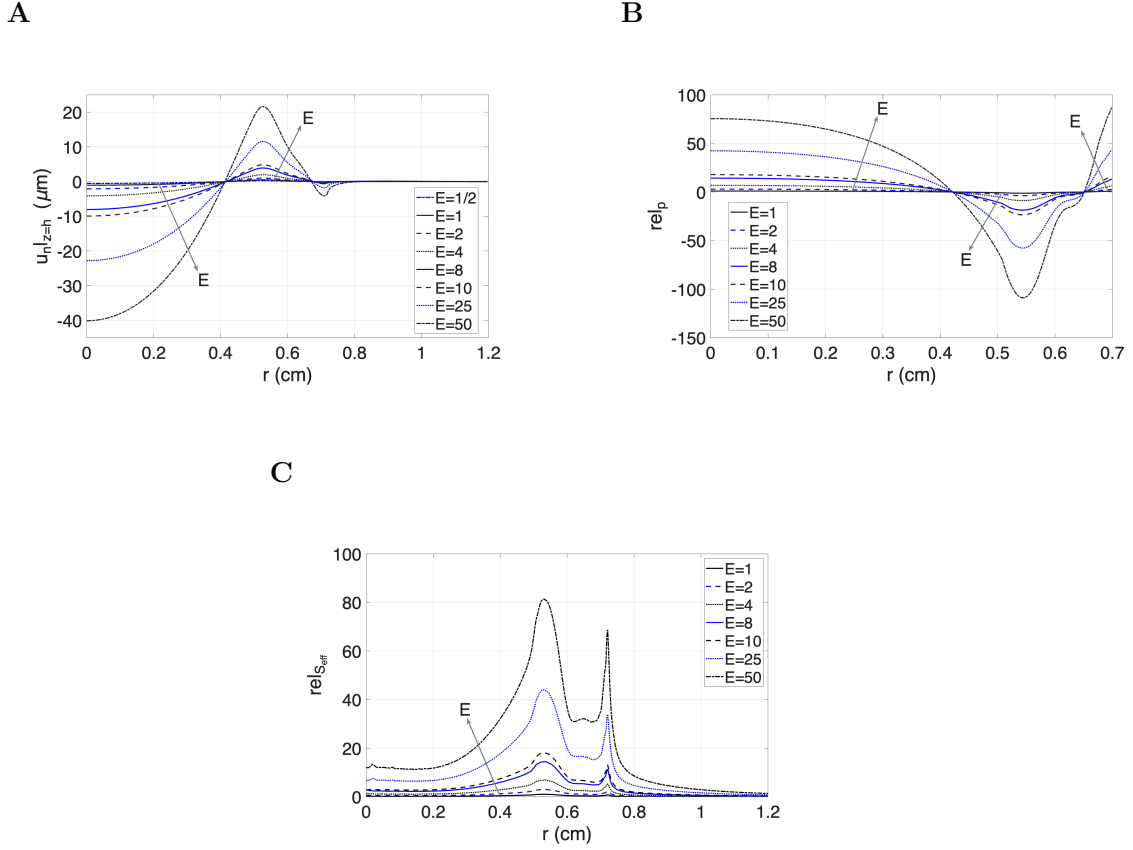

Figure~\ref{fig:material_props}{\bf A} shows the normal displacement $u_n$ of the ocular surface $\Gamma_{out}$ as $\mathrm{E}$ increases. \rev{The arrows indicate the direction which $\mathrm{E}$ increases.} The results show an inward ocular surface deformation at the center of the eye, outward at the limbus, and then inward at the end of the lens. We can see that the magnitude of the ocular normal displacements increases as $\mathrm{E}$ increases. If we double the value of $\mathrm{E}$ (i.e., a contact lens two time\rev{s} stiffer), then the maximum ocular normal displacement increases by a factor ranging from $1.78$ (increasing $\mathrm{E}$ from 25 to 50) to $1.98$ (for all smaller values of $\mathrm{E}$). 
The largest ocular normal displacements are of the order of microns. 
In the limit as $\mathrm{E}$ approaches infinity (i.e., a rigid lens), the deformed ocular surface would approach the posterior curve of the contact lens.  Differently, in the limit as $\mathrm{E}$ approaches zero (i.e., a rigid eye), the deformed ocular surface is given by the reference ocular surface.

Figure~\ref{fig:material_props}{\bf B} illustrates the effect of $\mathrm{E}$ on the contact lens suction pressure. For the different values of $\mathrm{E}$ considered, the figure illustrates the change in the suction pressure relative to the maximum absolute suction pressure when $\mathrm{E}=1/2$.
The relative difference of the suction pressure is computed as follows:
\begin{equation}
  rel_p(r; \mathrm{E}) =   \frac{p(r; \mathrm{E}) - p(r; \mathrm{E}=1/2)   }{ || p( r; \mathrm{E}=1/2)  ||_{\infty}}.
\end{equation}
\rev{Since we are accounting for the relative difference of pressure in $rel_p$, the vertical axis of Figure~\ref{fig:material_props}{\bf B} can be interpreted as a multiplication factor, positive or negative, to the maximum magnitude of $p$ predicted in the baseline scenario $\mathrm{E}=1/2$.}
We can see that the magnitude of $rel_p$ increases as $\mathrm{E}$ increases. For example, we find that when $\mathrm{E}=50$, the change in suction pressure can be at most $86$ times at the center of the eye and $109$ times at the limbus compared to $\mathrm{E}=1/2$. In the corneal region, the suction pressure change is large and positive. Consequently, the compression/flattening of the ocular surface increases as $\mathrm{E}$ increases. \rev{Recall that a positive suction pressure pushes downward on the ocular surface.} The range of $\mathrm{E}$ for currently manufactured soft contact lens materials and measurements of the ocular surface material properties is $0.5$ to $10$. Therefore, in the range of $\mathrm{E}$ from $0.5$ to $10$, as the contact lens becomes more stiffer relative to the eye, the ocular deformations can cause \rev{a change in the suction pressure of} at most $17$ \rev{times when compared to $\mathrm{E}=1/2$}.  

Figure~\ref{fig:material_props}{\bf C} illustrates the effect  of $\mathrm{E}$ on the ocular surface stresses.  The relative difference of the effective stress is computed as follows:
\begin{equation}
  rel_{S_{eff}}(r,h(r);\mathrm{E}) =   \frac{S_{eff}(r,h(r); \mathrm{E}) - S_{eff}(r,h(r); \mathrm{E}=1/2)   }{ || S_{eff}(r,h(r); \mathrm{E}=1/2)  ||_{\infty}}.
\end{equation}
To remove numerical instabilities due to the discontinuity of the suction at the end of the lens, in Figure~\ref{fig:material_props}{\bf C}, we plot the average value of $rel_{S_{eff}}$ throughout the outermost strip of tissue of thickness $0.02$~cm for each $r$ value. 
\rev{Since we are accounting for the relative difference of effective stress in $rel_{S_{eff}}$, the vertical axis of Figure~\ref{fig:material_props}{\bf C} can be interpreted as a multiplication factor, positive or negative, to the maximum value of $S_{eff}$ predicted in the baseline scenario $\mathrm{E}=1/2$.}
We find that the change in the effective ocular stress can be at most $81$ times the maximum effective stress when compared to $\mathrm{E}=1/2$. In general, the effective ocular stress increases as the contact lens becomes stiffer than the eye ($\mathrm{E}$ increases).  If we double the value of $\mathrm{E}$ (i.e., a contact lens two time\rev{s} stiffer), the maximum relative effective stress increases by a factor ranging from $1.8$ (increasing $\mathrm{E}$ from 25 to 50) to $3$ (for all smaller values of $\mathrm{E}$). The largest change in effective stress occurs at the beginning of the limbal region, where the effective stress is largest (see Figure~\ref{fig:stresses}{\bf F}). For existing commercial contact lenses and ``typical" eye tissue properties (i.e., $\mathrm{E}$ from $0.5$ to $10$), the effective stress can at most increase by 18 times its value when compared to $\mathrm{E}=1/2$.  

\subsubsection{Influence of contact lens thickness profile}
The contact lens thickness $\tau$ in the center of the lens is determined by the vision correction~\cite{douthwaite2006contact,stapleton2017impact}. The thickness profile at the edge of the lens (i.e., outside of the vision zone), is chosen to enhance on-eye performance (e.g., comfort)~\cite{douthwaite2006contact,srinivasan2015science,stapleton2017impact}. Here, we explore how a previously published contact lens thickness profile affects \rev{the} eye tissue \rev{response}. The thickness profile mimics the contact lens thickness profile studied by Funkenbusch and Benson \cite{funkenbusch1996conformity}. The center and edge lens thicknesses are $35$~$\mu$m. The lens thickness increases to $269$~$\mu$m as the radius increases from $0$ to $0.623$~cm, and then decreases to $35$~$\mu$m at the lens edge. The mathematical representation of the contact lens thickness is written in Eq~\eqref{eqn:thickness} \rev{and shown in Figure~\ref{fig:tau}} in the Appendix, Section~\ref{sec:supp_lens}.
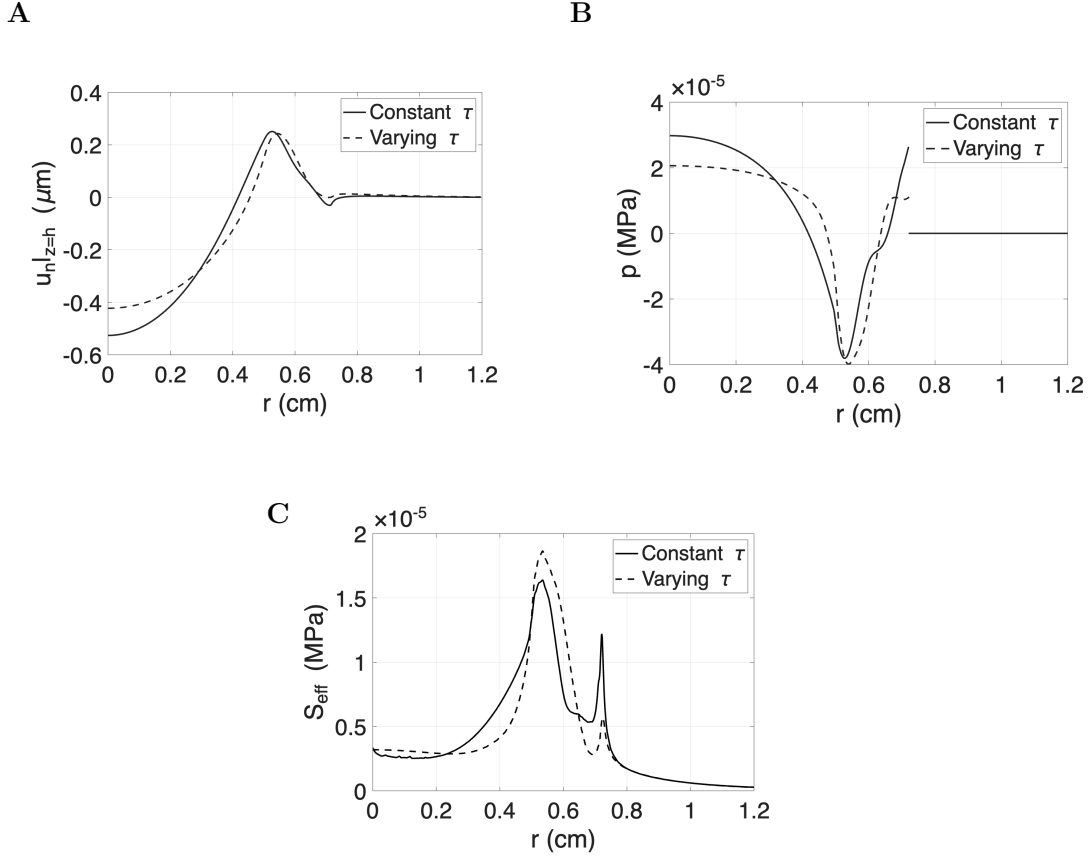
\begin{figure}[h!]
    \centering
        \scalebox{0.9}{
\begin{tikzpicture}
    \node at (0,0) {\includegraphics[width=.45\linewidth]{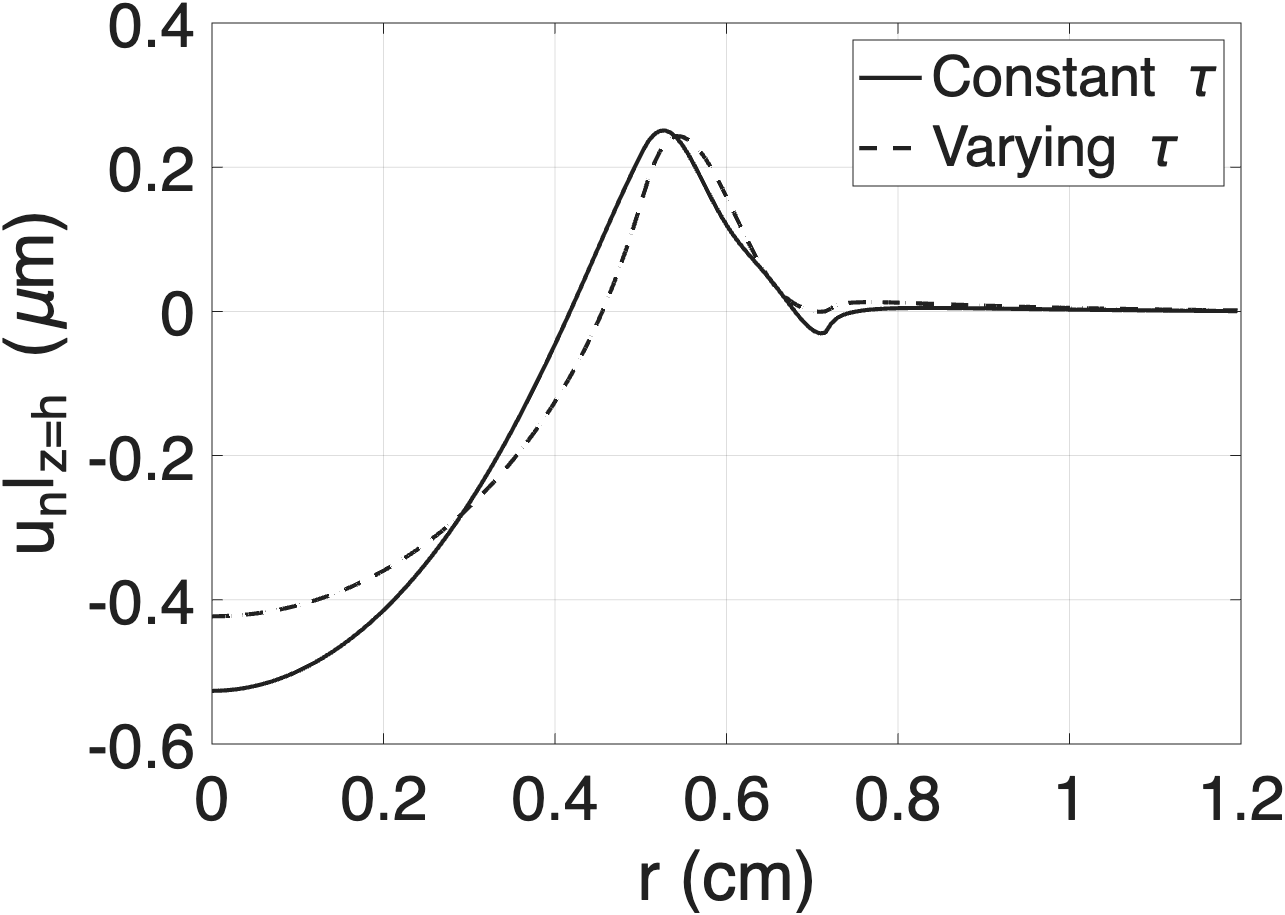}};
    \coordinate [label={left:\textcolor{black}{\bf \large A}}] (A) at (-3.3,3.5);
    \node at (8.6,0) {\includegraphics[width=.45\linewidth]{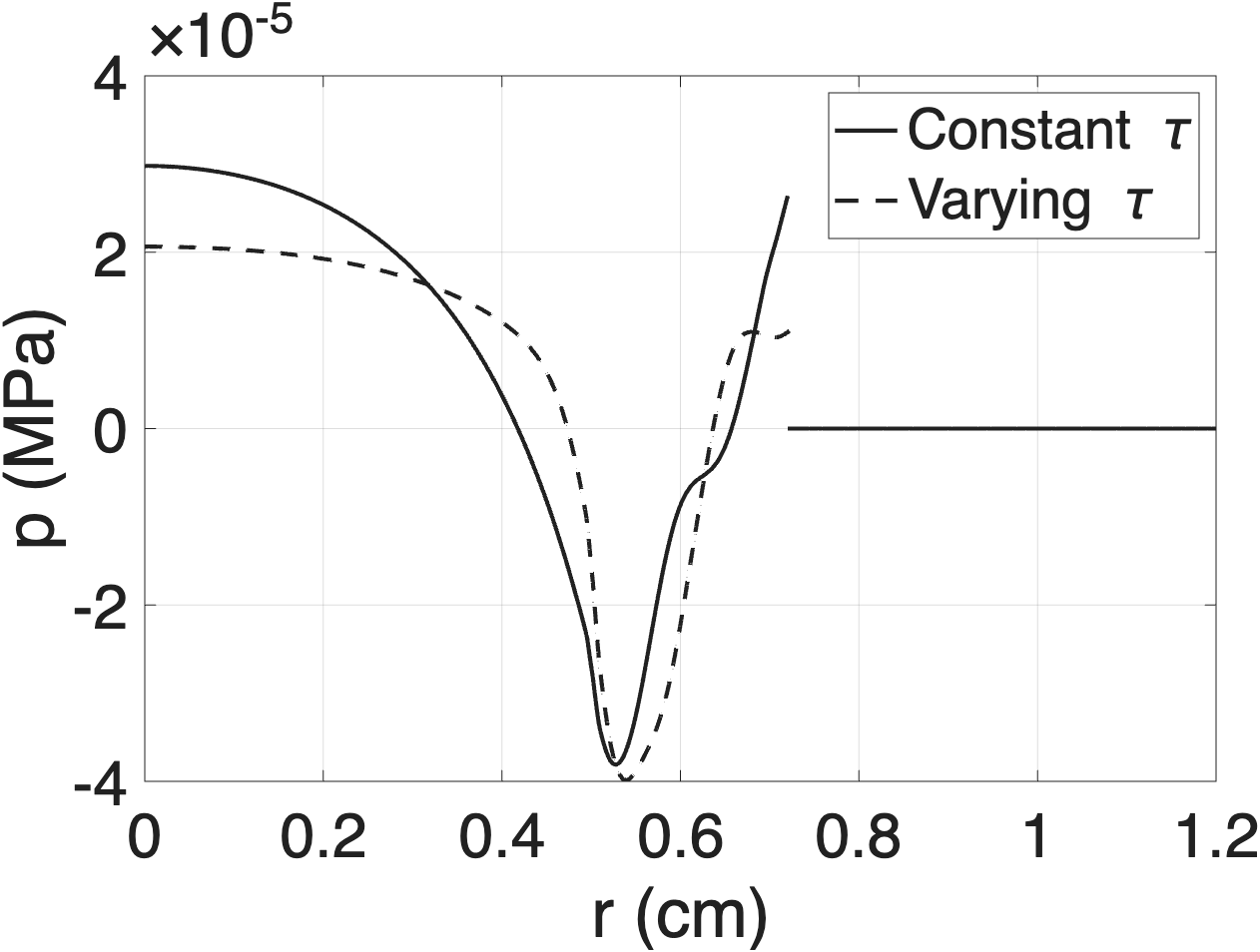}};
    \coordinate [label={left:\textcolor{black}{\bf \large B}}] (B) at (4.95,3.5);
    \node at (4,-6.3) {\includegraphics[width=.45\linewidth]{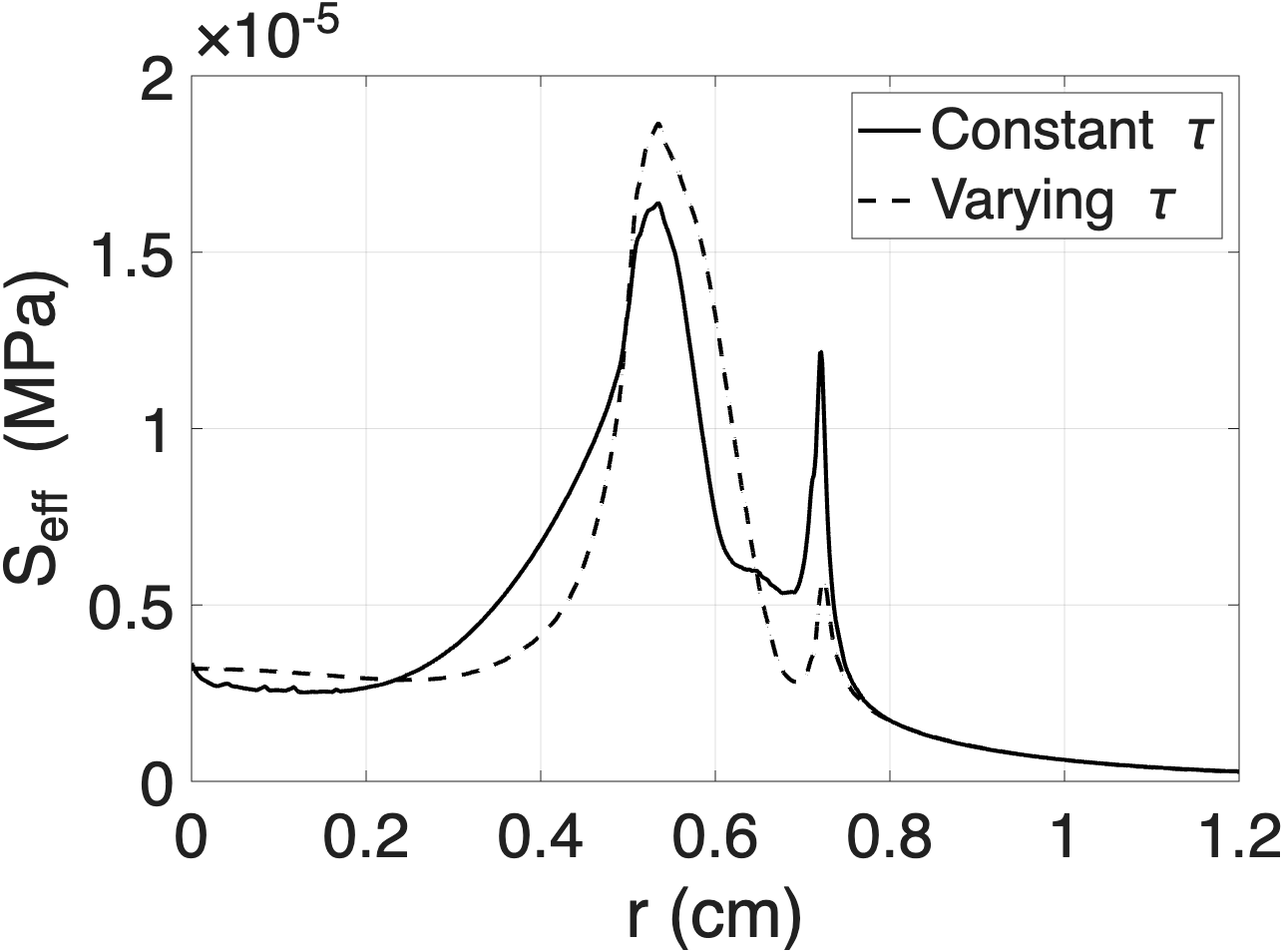}};
    \coordinate [label={left:\textcolor{black}{\bf \large C}}] (C) at (0.5,-3.85);
\end{tikzpicture}}
    \caption{\textit{Homogeneous average-shaped eye and average-shaped contact lens of constant thickness (solid line) and varying thickness (dashed line)}. \rev{({\bf A})} The ocular surface normal displacement $u_n$, \rev{({\bf B})} the suction pressures $p$, and \rev{({\bf C})} the effective stress $S_{eff}$ as a function of the radial reference coordinate $r$. }
    \label{fig:thickness}
\end{figure}

Figure~\ref{fig:thickness} displays the ({\bf A}) normal ocular deformations, ({\bf B}) lens suction pressure, and ({\bf C}) ocular surface effective stress  for a constant lens thickness $\tau =100$~$\mu$m (solid line) and a varying thickness $\tau(r)$ (Appendix Eq~\eqref{eqn:thickness} in Section~\ref{sec:supp_lens}). \rev{The solid line results reported in Figure~\ref{fig:thickness} are the same as the results reported in Figures~\ref{fig:standard} and~~\ref{fig:stresses}}. The varying thickness lens results in less deformation in the center of the eye. The results show that the suction pressure at the center of the eye is proportional to the contact lens thickness: the thinner lens at the center ($35$~$\mu$m compared to $100$~$\mu$m at $r=0$) results in a smaller\rev{, positive} suction pressure and less inward ocular deformations. 
Similarly, at the lens edge, the normal ocular displacements is smaller if the thickness of the lens varies since the suction pressure is smaller for the varying thickness lens ($35$~$\mu$m compared to $100$~$\mu$m at $r=\mathfrak{R_{lens}}=0.7cm$). In the limbal region, the varying lens thickness ranges from $37$~$\mu$m to $291$~$\mu$m.
The \rev{results} of the constant vs varying thickness lenses are similar in the limbal regions. This suggests that the curvature changes in the limbal region are the main factor that determines the \rev{lens-eye interactions}. Figure~\ref{fig:thickness}{\bf C} compares the effective stresses near the ocular surface for the constant and varying lens thickness scenarios. Similarly to Figure~\ref{fig:material_props}{\bf C}, Figure~\ref{fig:thickness}{\bf C} shows the average value of $rel_{S_{eff}}$ throughout the outermost strip of a tissue $0.02$~cm thick for each $r$ value.
The largest differences are found at the edge of the lens. 

\subsubsection{Influence of the contact lens and ocular shapes}

We return to our baseline reference material parameter values, $\mathrm{E}=1/2$, $\sigma_{eye}=\sigma_{lens} = 0.49$, and $\tau=100$~$\mu$m, and study how different combinations of contact lens and ocular shapes affects their interaction. Figure~\ref{fig:shapes} displays the lens and ocular predictions for the different eye and lens combinations (four eyes shapes and three lens shapes, as described in the Appendix, Sections~\ref{sec:supp_lens} and~\ref{sec:supp_eye}). Each row of Figure~\ref{fig:shapes} corresponds to a different eye shape. Each column plots a different output from the model. In each plot, the prediction for the flat lens (solid line), average lens (dashed line), and steep lens (dotted line) are shown. 

The contact lens mechanics are shown in the first and second columns of Figure~\ref{fig:shapes}. The lens radial displacement $\eta_r$ is plotted in first column. In general, for all the combinations considered, we find the contact lenses are slightly stretched radially for most of the corneal region ($0<r<0.4$~cm), and the lens are radially compressed in the limbal region ($0.5<r<0.7$~cm). Such compression is released at the edge of the lens ($r=0.7$~cm), where the lens is radially stretched (except for the average cornea with flat sclera and a flat lens) in most of the combinations considered. The suction pressure is plotted in second column of  Figure~\ref{fig:shapes}.  The general shape of the suction pressure profile is similar for all the eye and contact lens combinations considered in this study.  The suction pressure has a local maximum in the center of the eye, which monotonically decreases to a global minimum located in the first half of the limbal region ($0.5<r<0.6$~cm).  Then, the suction pressure either monotonically or non-monotonically increases to the a local maximum at the edge of the contact lens ($r=0.7$~cm). 

There are general relationships between differences in eye and lens shapes and the predicted contact lens mechanics (first and second columns in Figure~\ref{fig:shapes}).  Specifically, the flat-, average-, and steep-shaped contact lenses' radii of curvature are larger than the radii of curvature of the corneas considered.  Therefore, if the contact lens' radius of curvature decreases (moving from flat to average to steep) or if the cornea's radius of curvature decreases (moving from flat to average to steep), then, in the corneal region ($0<r<0.5$~cm), the difference between the eye's and the lens' radii of curvature decreases (see Figure~\ref{fig:eye_geometry}). 
Consequently, in the corneal region, the magnitude of the contact lens radial displacement (see  first column in Figure~\ref{fig:shapes}) and the  magnitude of the suction pressure (see second column in Figure~\ref{fig:shapes}) decrease from the flat-shaped lens to the steep-shaped lens. Similarly, in the corneal and limbal regions, except at the end of the lens, the magnitude of the contact lens radial displacement and the magnitude of the suction pressure decrease from flat-shaped cornea to the steep-shaped cornea with average sclera.  

On the other hand, the difference between the contact lens' radius of curvature and the sclera's radius of curvature increases (see Figure~\ref{fig:eye_geometry}) from flat- to steep-shaped lens\rev{es} (see two first columns and first three rows in Figure~\ref{fig:shapes}).  Therefore, the magnitude of the radial displacement and the magnitude of the suction pressure increase at the edge of the contact lens from flat- to steep-shaped lens. At the edge of the lens $r=0.7$~cm, the magnitude of the suction pressure and of the lens radial displacement increases from flat- to average-shaped sclera (see two first columns and the second and fourth rows of Figure~\ref{fig:shapes}). 
For any ocular shape considered, there is a trade off.  A reduction in the suction pressure at the center of the eye\rev{,} resulting from inserting a steeper contact lens\rev{,} leads to a larger suction pressure at the edge of the contact lens. For a fixed lens, the suction pressure minimum value and the value at the edge of the lens change depending on the corneal and scleral radii of curvature.

While the general shape of the suction pressure is the same for the different shaped lenses and eyes, the global and local maxima and minima of the suction pressure profiles are different.  Consequently, the magnitude of ocular surface normal displacements (third column in  Figure~\ref{fig:shapes}) and the effective stress on the ocular surface (fourth column in  Figure~\ref{fig:shapes}) are different for the different lens/eye combinations. For any eye shape considered, the eye deformations decrease in the corneal region ($0<r<0.4$~cm) and increase at the edge of the lens ($r=0.7$~cm) from the flat-shaped lens to the steep-shaped lens. 
For an average shaped sclera, the effective stress experienced by the ocular surface decreases from the flat-shaped lens to the steep-shaped lens, except near the edge of the lens ($r=0.7$~cm) (see fourth column and first three rows in Figure~\ref{fig:shapes}), where the stresses increase from the flat to steep lens. For a flat sclera, the ocular surface effective stress decreases from the flat-shaped lens to the steep-shaped lens for any value of $r$. For a fixed lens, as the cornea changes from flat-shaped to steep-shaped, the ocular surface deformation and the effective stress at in the limbus region ($0.5<r<0.7$~cm) decrease. 
The flat cornea paired with the flat-shaped contact lens produce the largest ocular effective stresses, which correspond to the largest magnitudes of suction pressures (see second column in Figure~\ref{fig:shapes}). Interestingly, the flat sclera produces smaller effective stresses in the center of the eye, even though the suction pressure is large in magnitude there. The largest ocular effective stress shifts from the limbal region on the average cornea, flat cornea, and steep cornea eye shapes closer to the beginning of the sclera ($r=0.7$~cm) on the flat sclera eye shape. 

\begin{landscape}
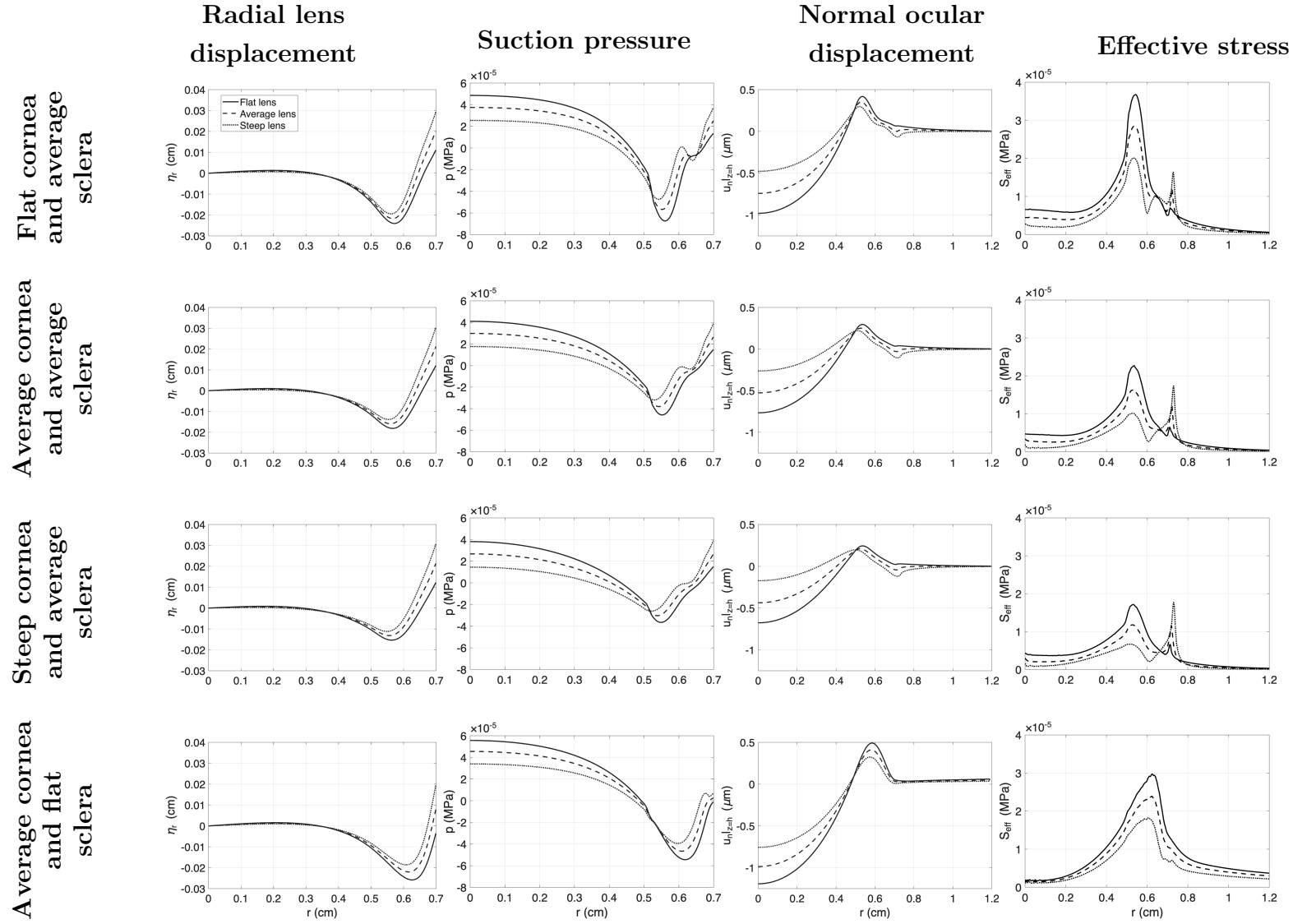
\begin{figure}[htb]
\centering
\begin{tikzpicture}
\node at (1.5,-1.5) {\includegraphics[width=0.22\linewidth]{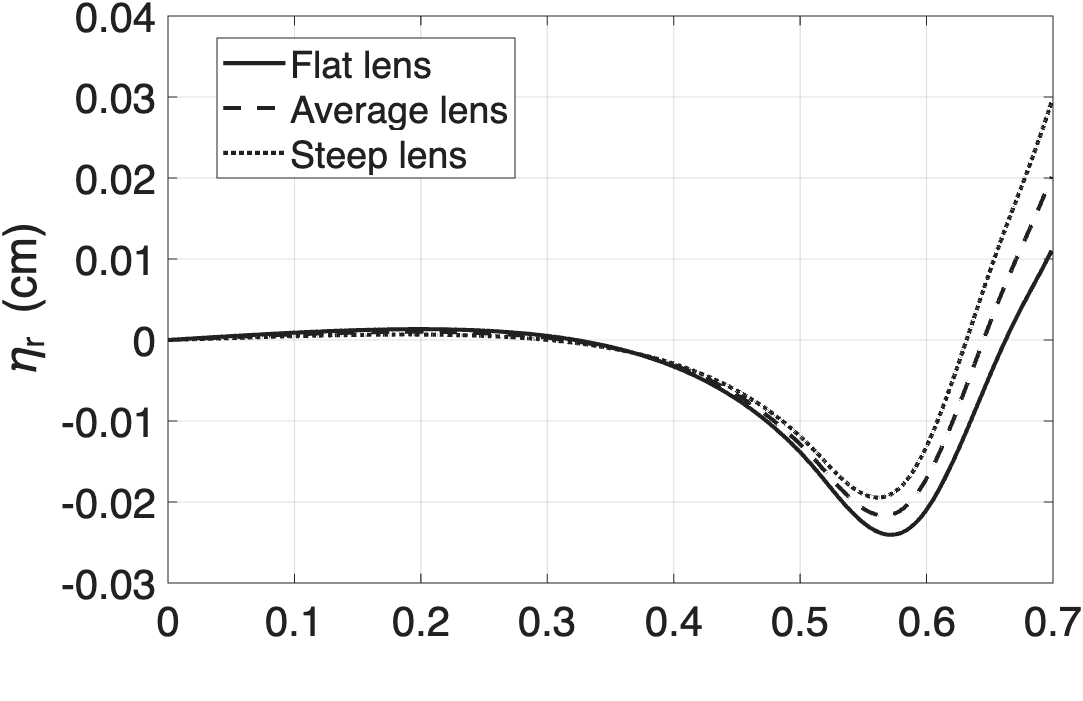}\includegraphics[width=0.22\linewidth]{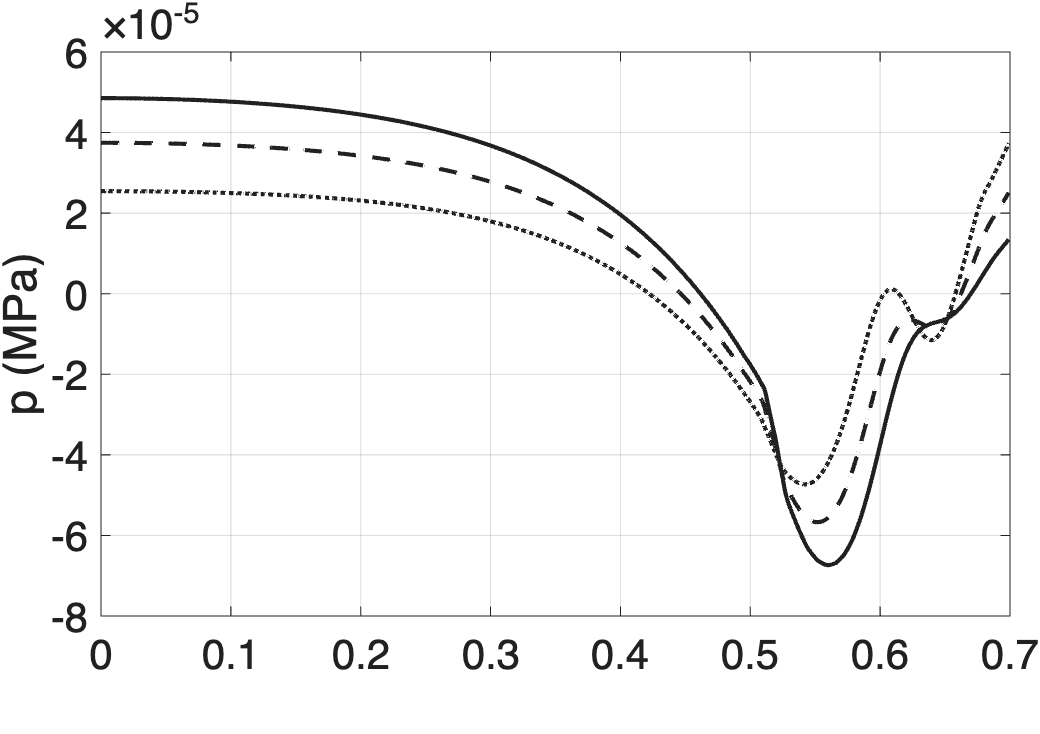}\includegraphics[width=0.22\linewidth]{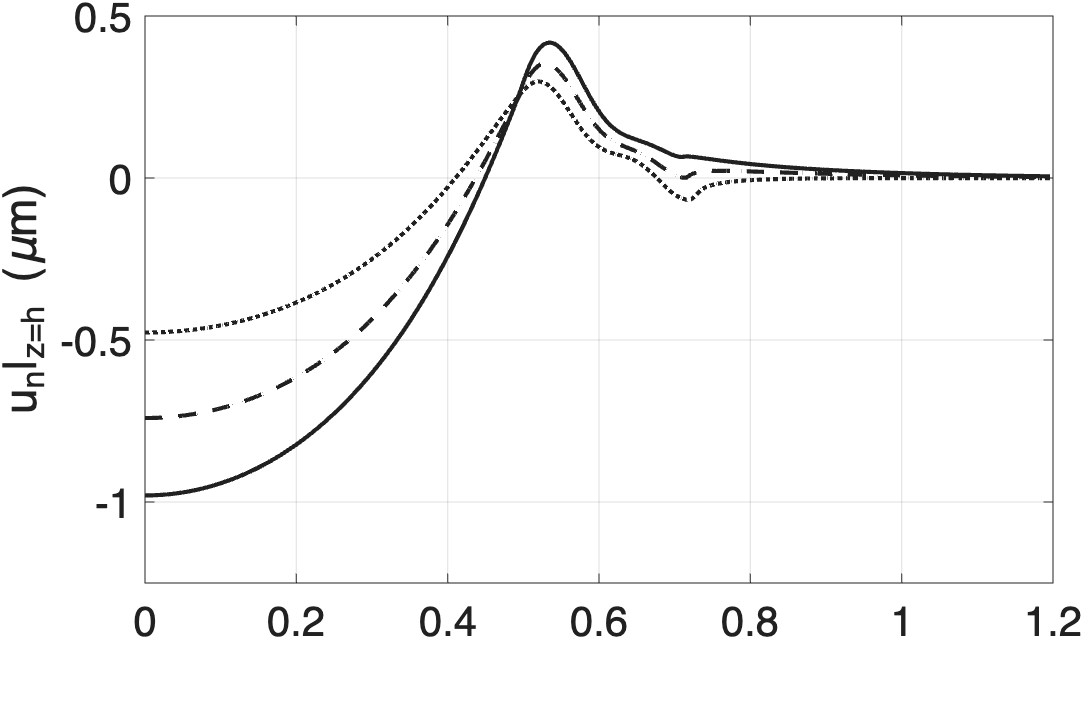}\includegraphics[width=0.22\linewidth]{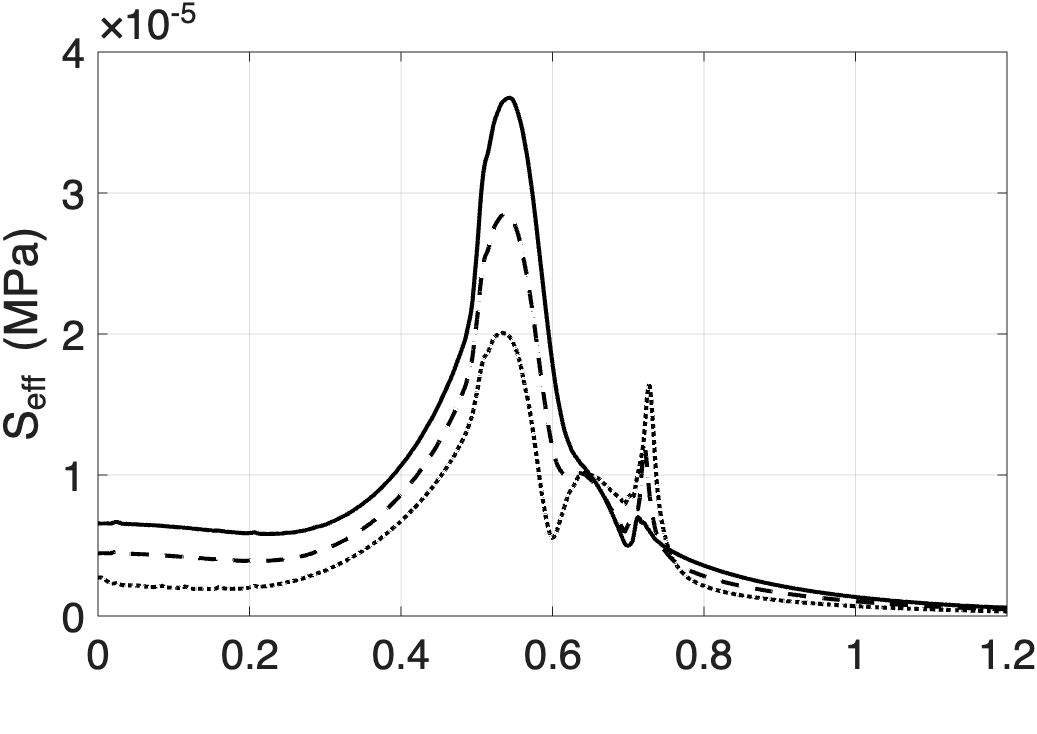}};
\node at (1.5,-5) {\includegraphics[width=0.22\linewidth]{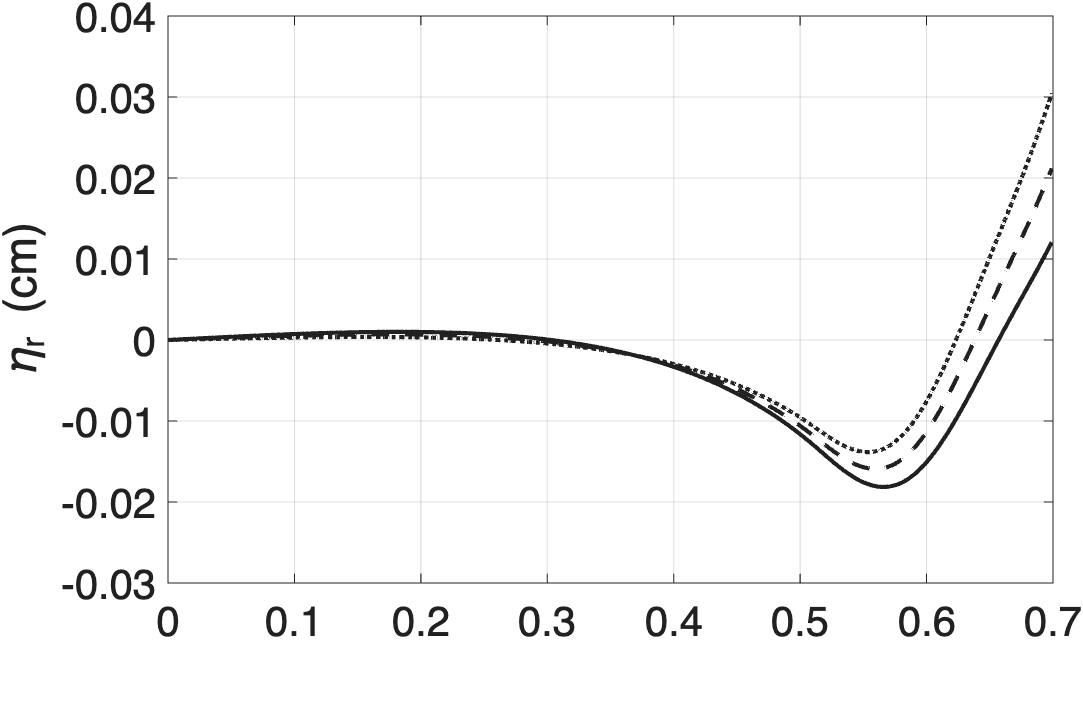}\includegraphics[width=0.22\linewidth]{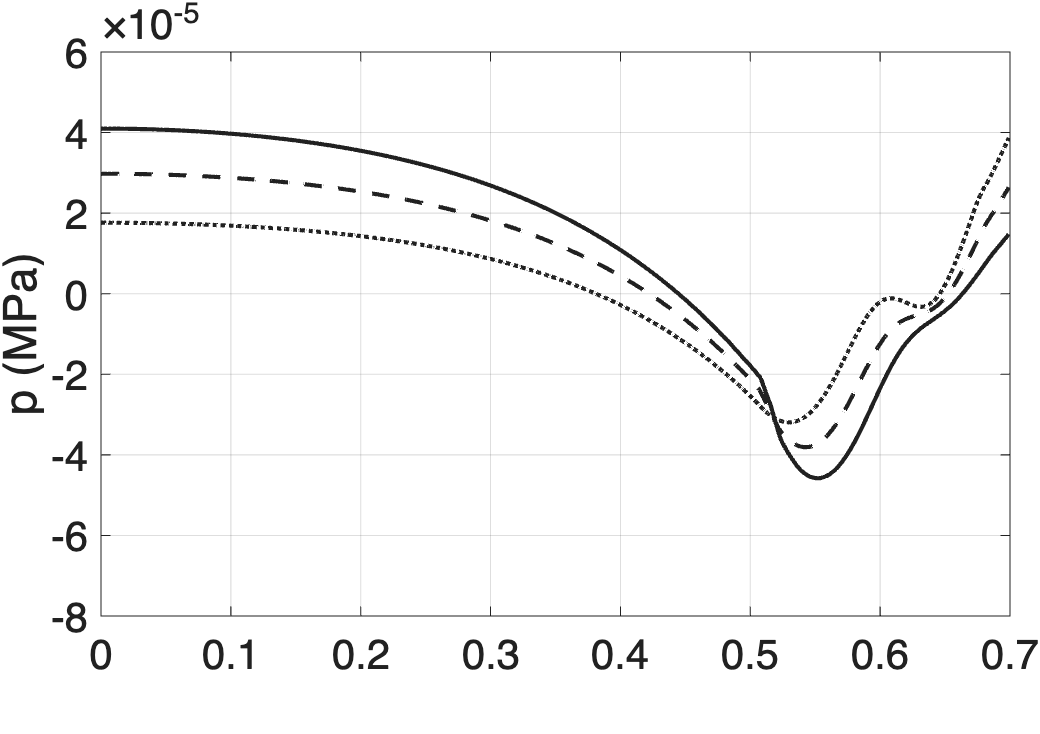}\includegraphics[width=0.22\linewidth]{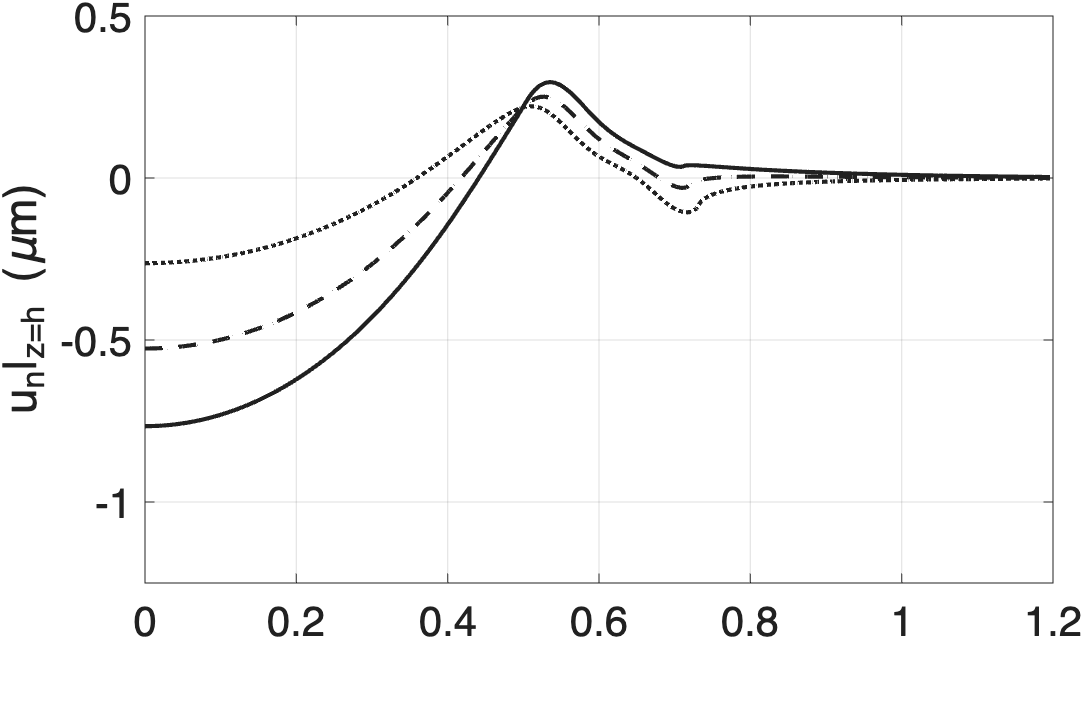}\includegraphics[width=0.22\linewidth]{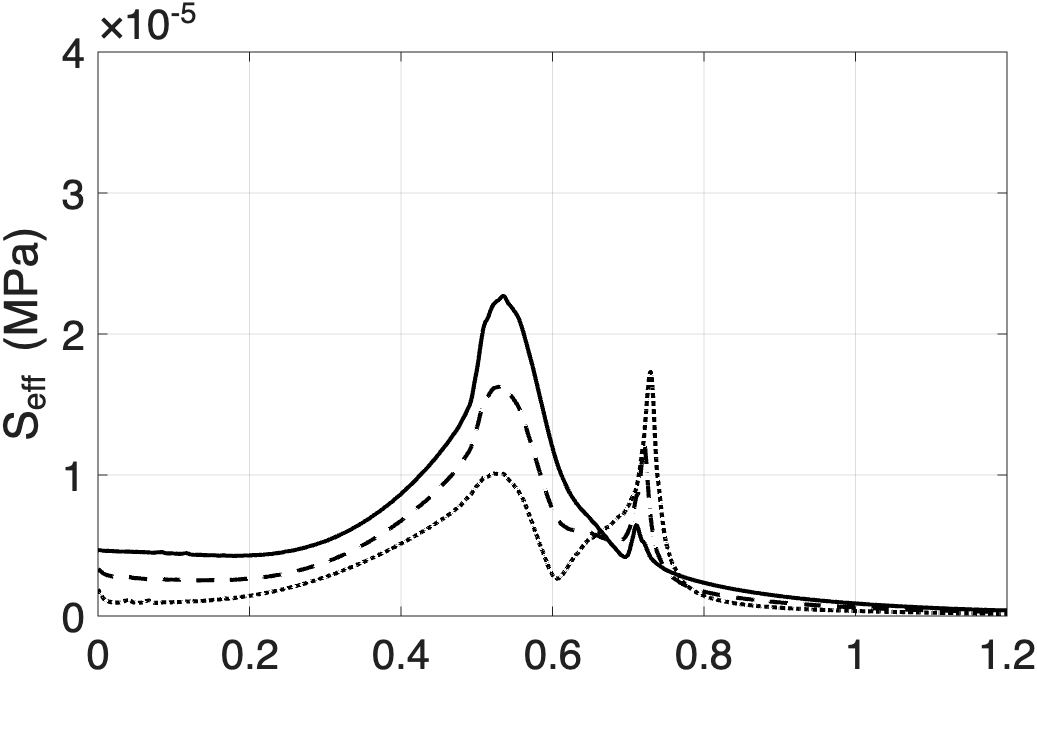}};
\node at (1.5,-8.5) {\includegraphics[width=0.22\linewidth]{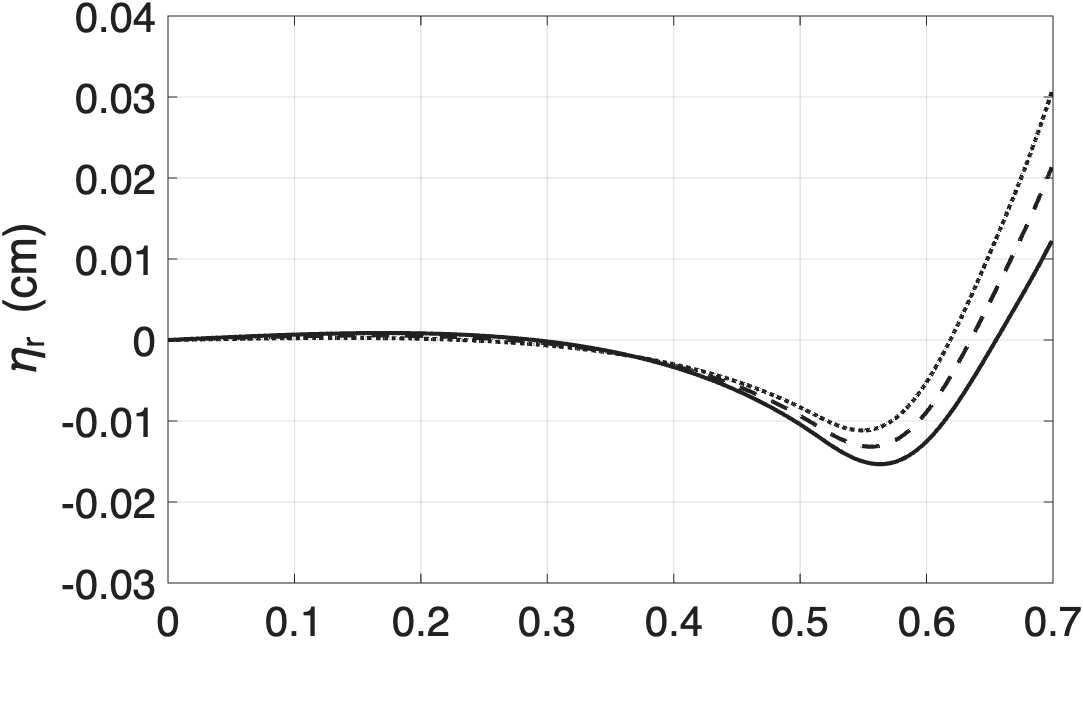}\includegraphics[width=0.22\linewidth]{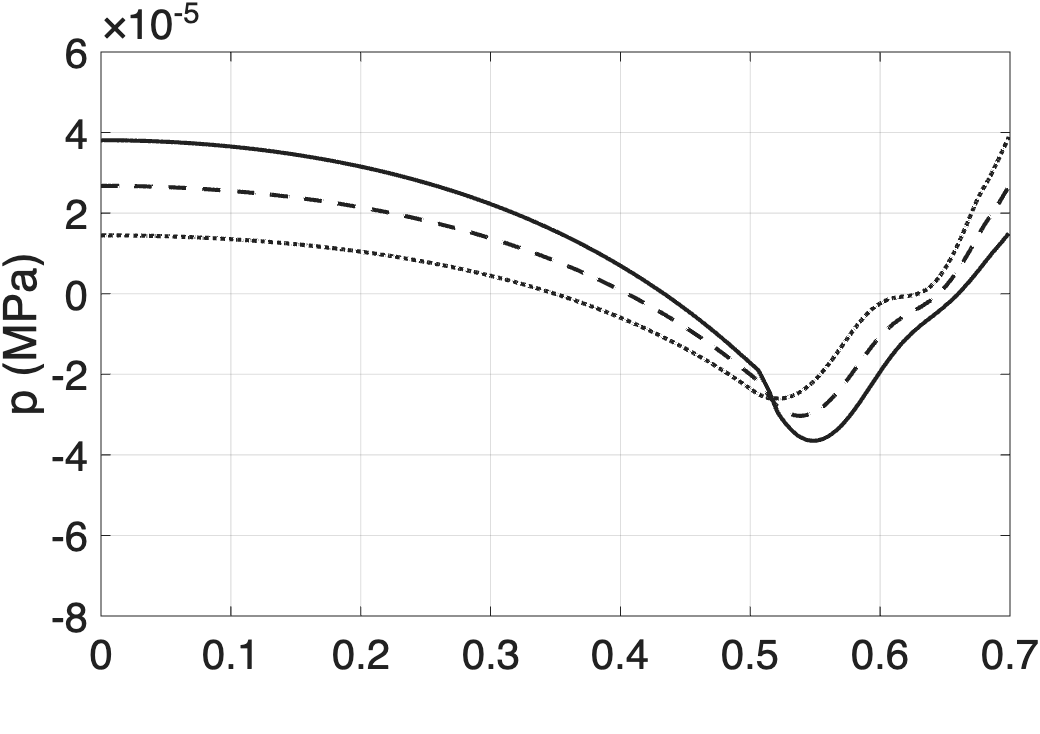}\includegraphics[width=0.22\linewidth]{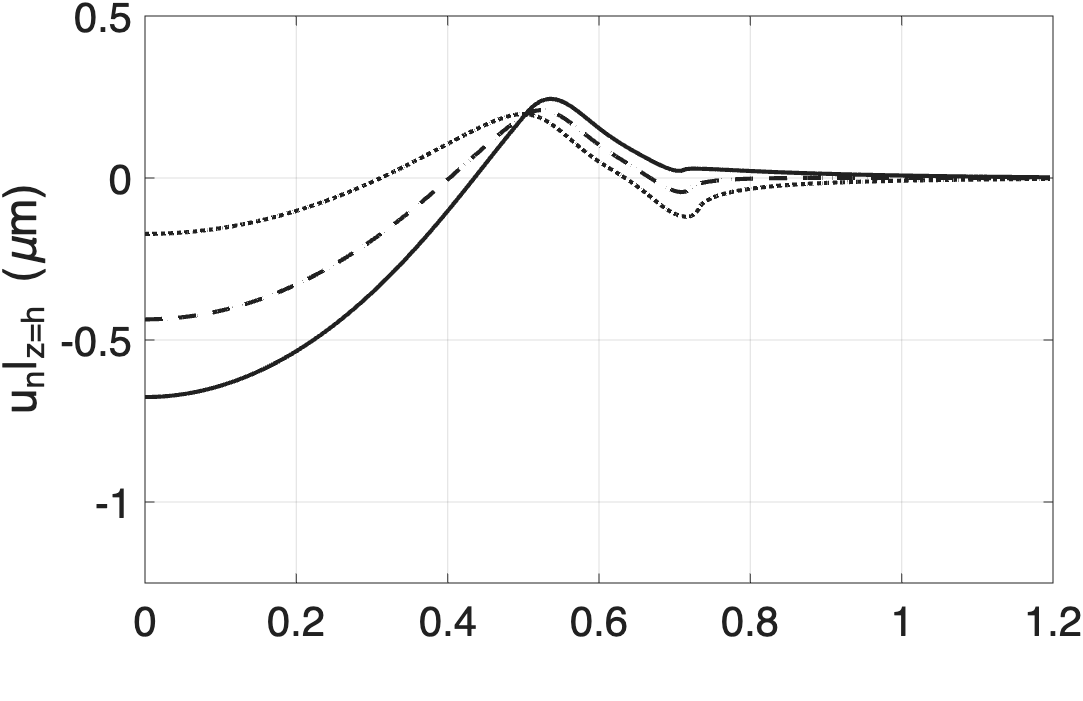}\includegraphics[width=0.22\linewidth]{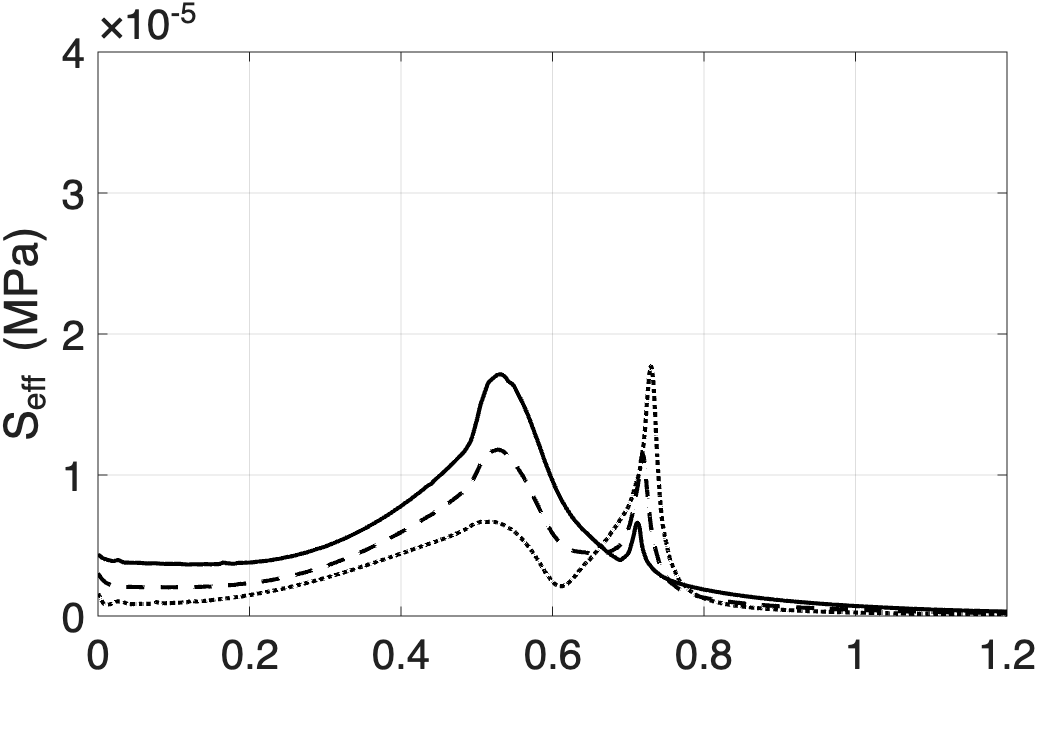}};
\node at (1.5,-12) {\includegraphics[width=0.22\linewidth]{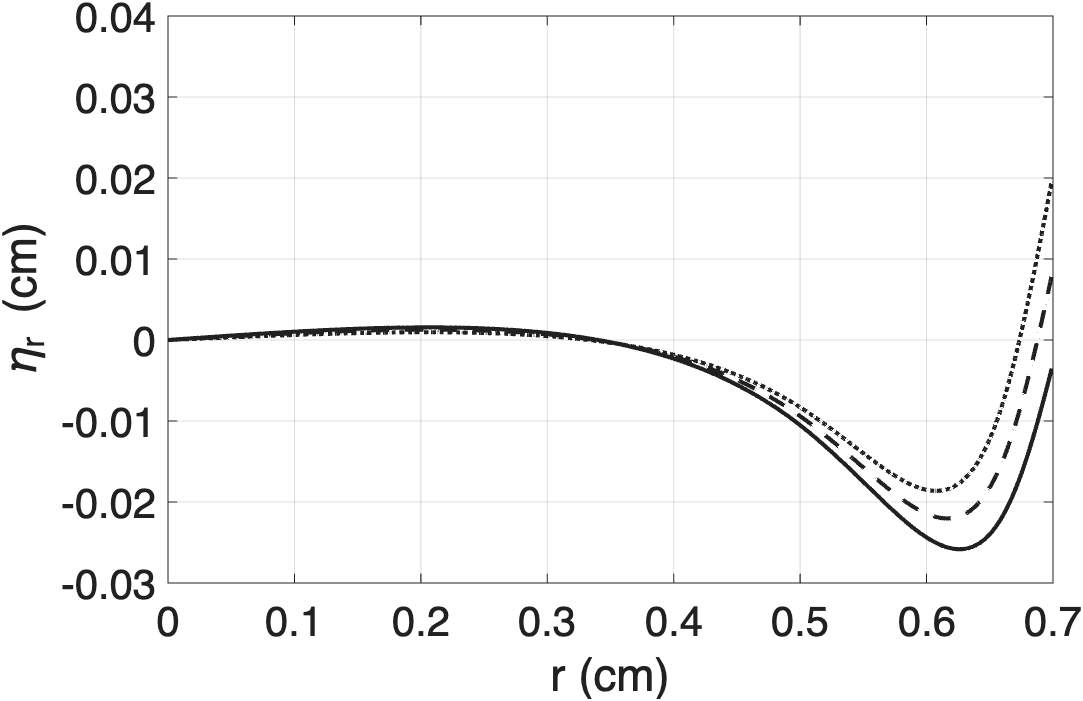}\includegraphics[width=0.22\linewidth]{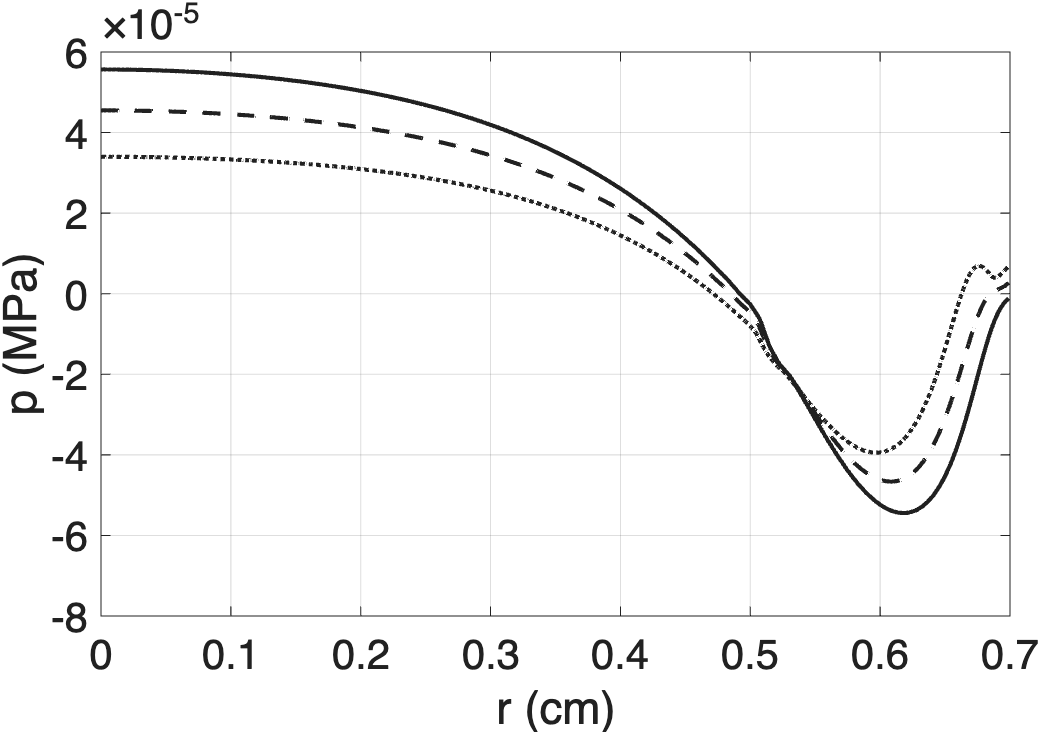}\includegraphics[width=0.22\linewidth]{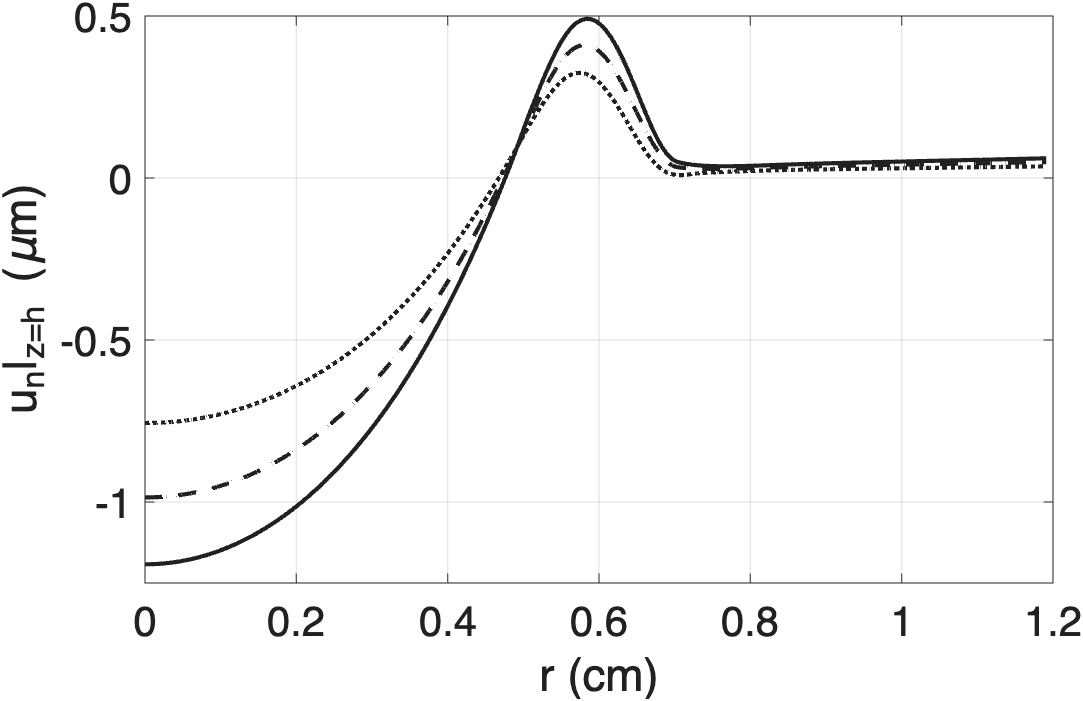}\includegraphics[width=0.22\linewidth]{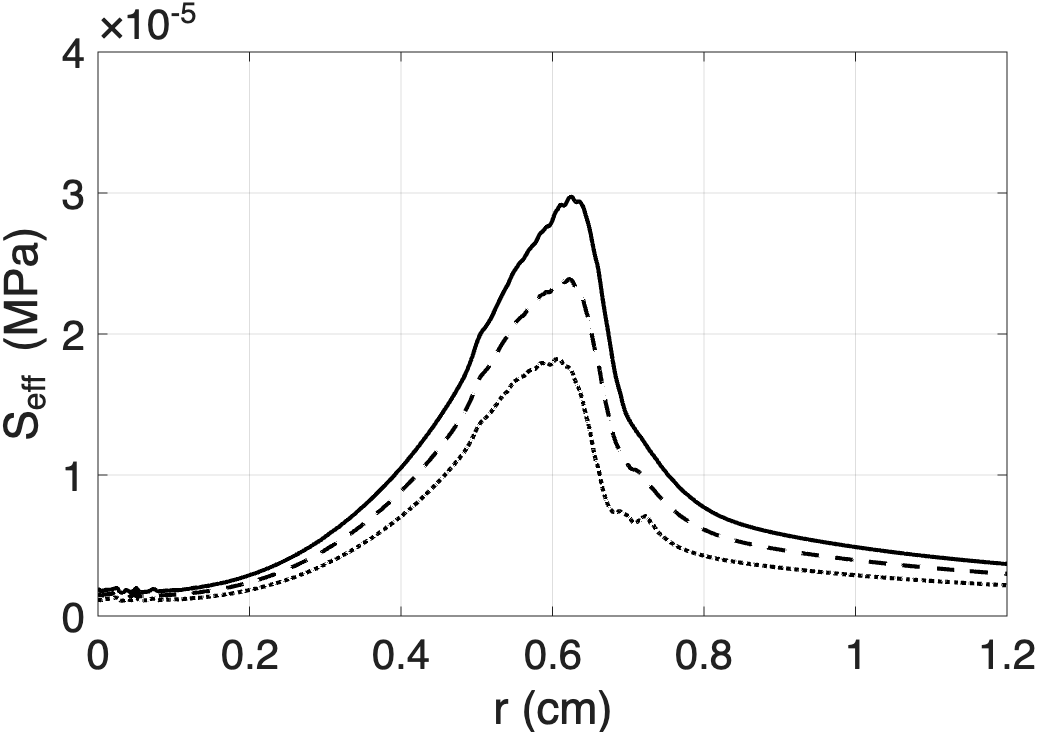}};
\coordinate[label={\textcolor{black}{\bf \large Radial lens}}] (c1a) at (-5.7,.7);
\coordinate[label={\textcolor{black}{\bf \large displacement}}] (c1b) at (-5.7,0);
\coordinate[label={\textcolor{black}{\bf \large Suction pressure}}] (c2) at (-.7,0.2);
\coordinate[label={\textcolor{black}{\bf \large Normal ocular}}] (c3a) at (4.25,.7);
\coordinate[label={\textcolor{black}{\bf \large displacement}}] (c3b) at (4.25,0);
\coordinate[label={\textcolor{black}{\bf \large Effective stress}}] (c4) at (9.1,.2);
\coordinate[label={[rotate=90]\textcolor{black}{\bf \large Flat cornea}}] (r1) at (-9.4,-1.5);
\coordinate[label={[rotate=90]\textcolor{black}{\bf \large and average}}] (r14b) at (-8.9,-1.5);
\coordinate[label={[rotate=90]\textcolor{black}{\bf \large sclera}}] (r14b) at (-8.5,-1.5);
\coordinate[label={[rotate=90]\textcolor{black}{\bf \large Average cornea}}] (r2) at (-9.4,-4.9);
\coordinate[label={[rotate=90]\textcolor{black}{\bf \large and average}}] (r14b) at (-8.9,-4.9);
\coordinate[label={[rotate=90]\textcolor{black}{\bf \large sclera}}] (r14b) at (-8.5,-4.9);
\coordinate[label={[rotate=90]\textcolor{black}{\bf \large Steep cornea}}] (r3) at (-9.4,-8.5);
\coordinate[label={[rotate=90]\textcolor{black}{\bf \large and average}}] (r14b) at (-8.9,-8.5);
\coordinate[label={[rotate=90]\textcolor{black}{\bf \large sclera}}] (r14b) at (-8.5,-8.5);
\coordinate[label={[rotate=90]\textcolor{black}{\bf \large Average cornea}}] (r4a) at (-9.4,-12);
\coordinate[label={[rotate=90]\textcolor{black}{\bf \large and flat}}] (r14b) at (-9,-12);
\coordinate[label={[rotate=90]\textcolor{black}{\bf \large sclera}}] (r14b) at (-8.5,-12);

\end{tikzpicture}
\caption{\textit{Homogeneous eye model.} Radial lens displacement (first column), suction pressure (second column), normal ocular displacement (third column), and effective stress (fourth column) profiles for different contact lens and eye combinations. Results for three lens shapes, flat lens (solid line), average lens (dashed line), and steep lens (dotted line), are shown in each plot. Each row corresponds to a different ocular shape considered. Details on the lens and eye shapes considered are in the Appendix, Sections~\ref{sec:supp_lens} and~\ref{sec:supp_eye}.
}
\label{fig:shapes}
\end{figure}
\end{landscape} 

\subsection{Heterogeneous eye model}
\label{sec:var_eye_res}
Lastly, an eye tissue with a spatially-dependent Young's modulus  is explored, which we refer to as the heterogeneous eye model. From Section~\ref{sec:eye_model}, recall that the Young's modulus is smaller at the center of the eye and increases in magnitude along rays from the center of the eye to the ocular surface. Then, along the ocular surface, the sclera is five times stiffer than the corneal region. For more details on the spatially-depend\rev{ent} Young's modulus, see the Appendix, Section~\ref{sec:supp_YM}. In what follows, we consider an average-shaped lens with a constant thickness $\tau=100$~$\mu$m and $E_{lens}=0.1$~MPa on an average-shaped eye.

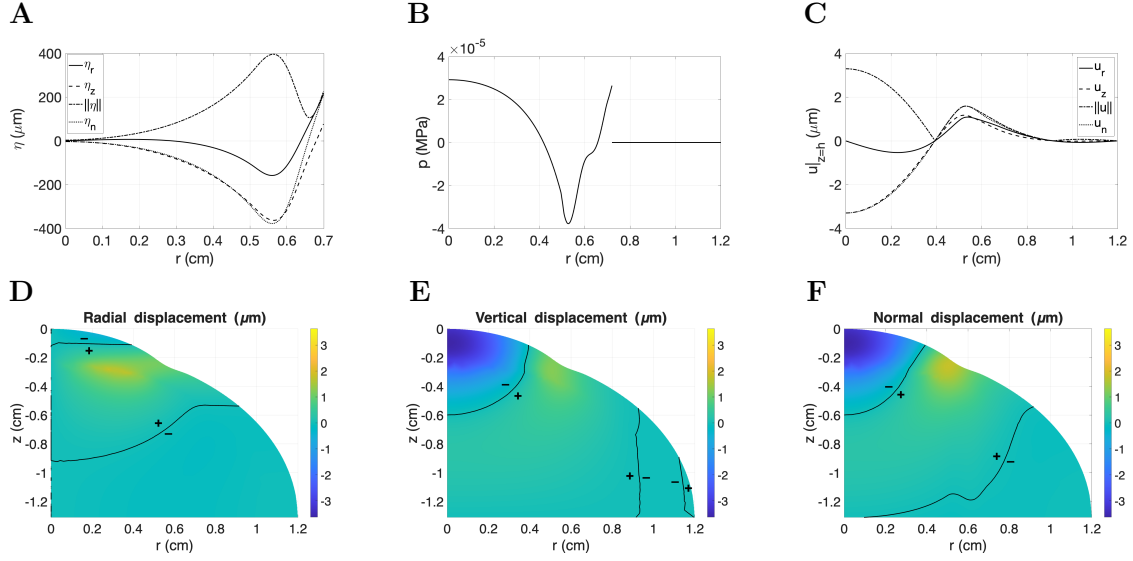
\begin{figure}[h!]
\begin{center}
\scalebox{0.92}{
\begin{tikzpicture}
\node at (-.2,0) {\includegraphics[width=0.33\linewidth]{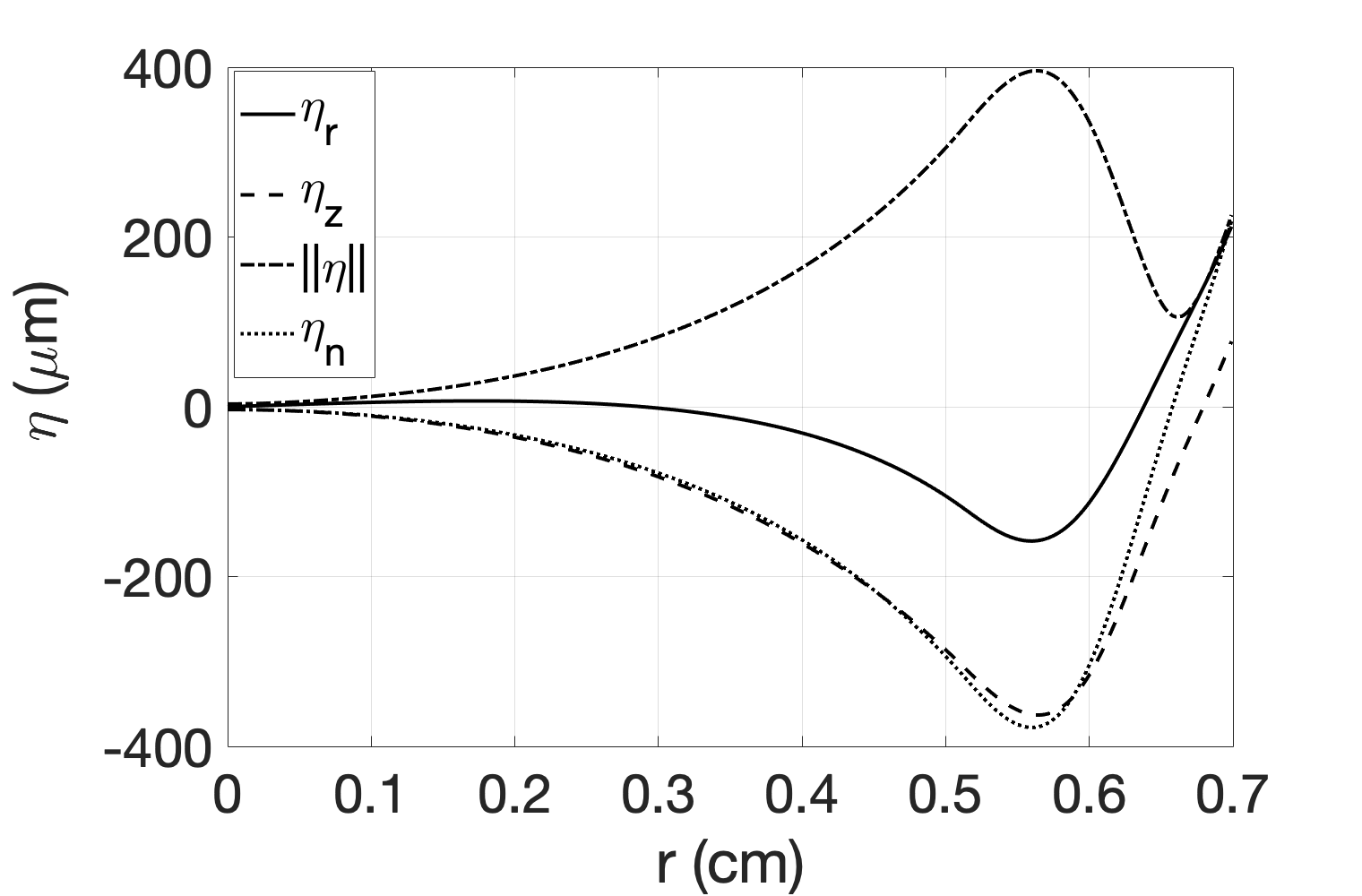}};
\coordinate [label={left:\textcolor{black}{\bf \large A}}] (A) at (-2.2,2);
\node at (5.5,0) {\includegraphics[width=0.33\textwidth]{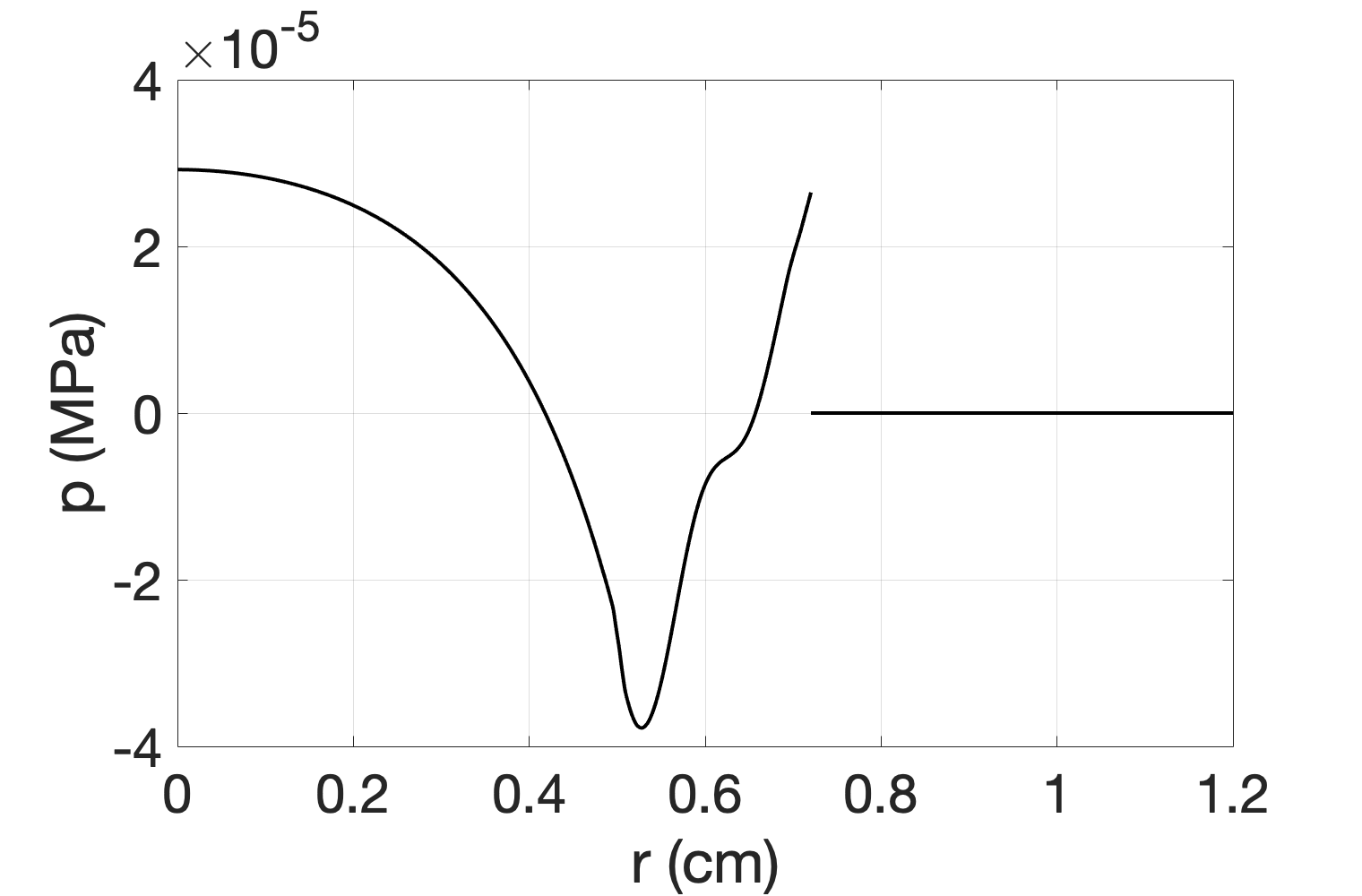}};
\coordinate [label={left:\textcolor{black}{\bf \large B}}] (B) at (3.5,2);
\node at (11.2,0) {\includegraphics[width=0.33\textwidth]{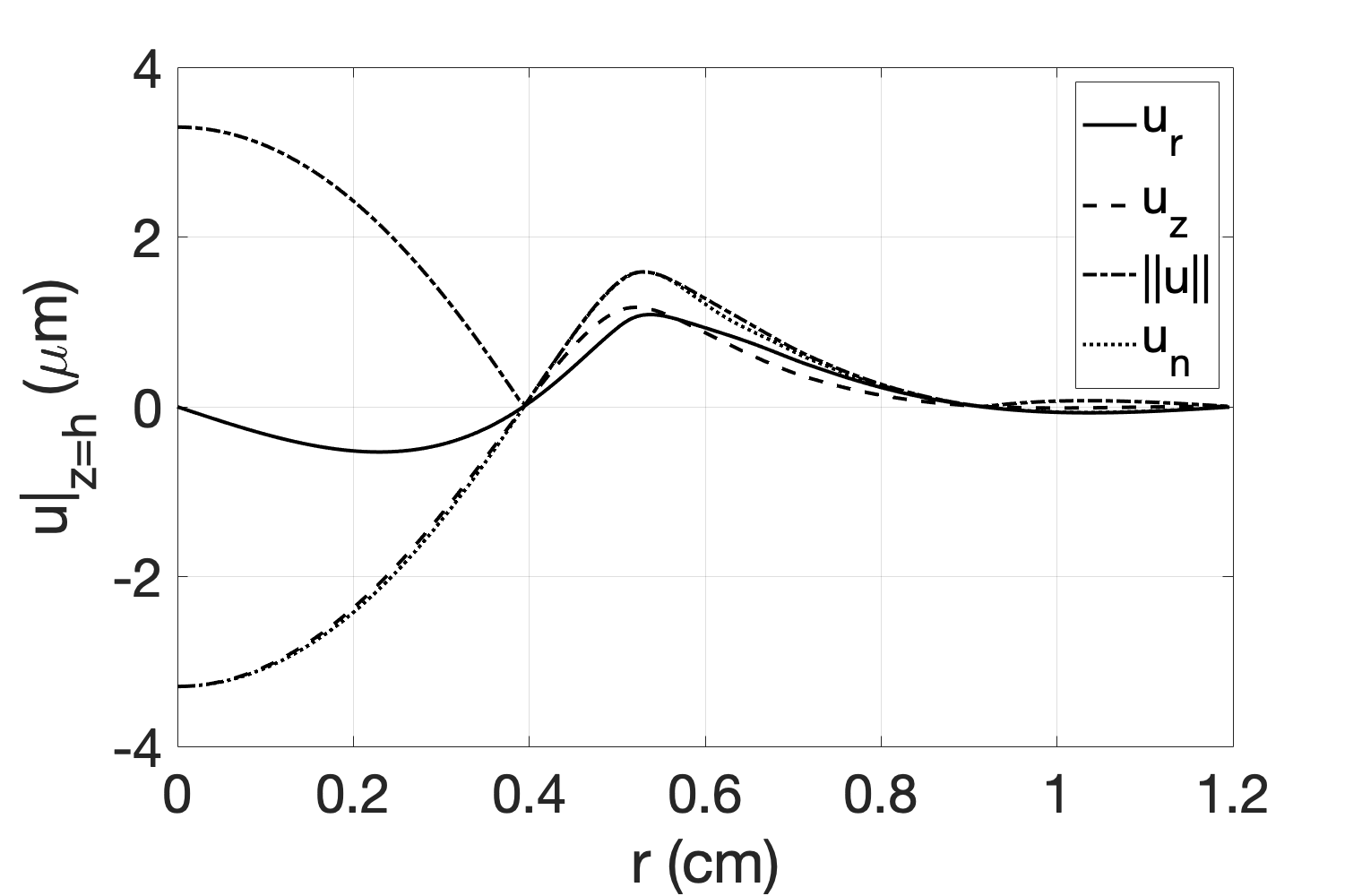}};
\coordinate [label={left:\textcolor{black}{\bf \large C}}] (C) at (9.2,2);
\node at (-.2,-4) {\includegraphics[width=0.33\linewidth]{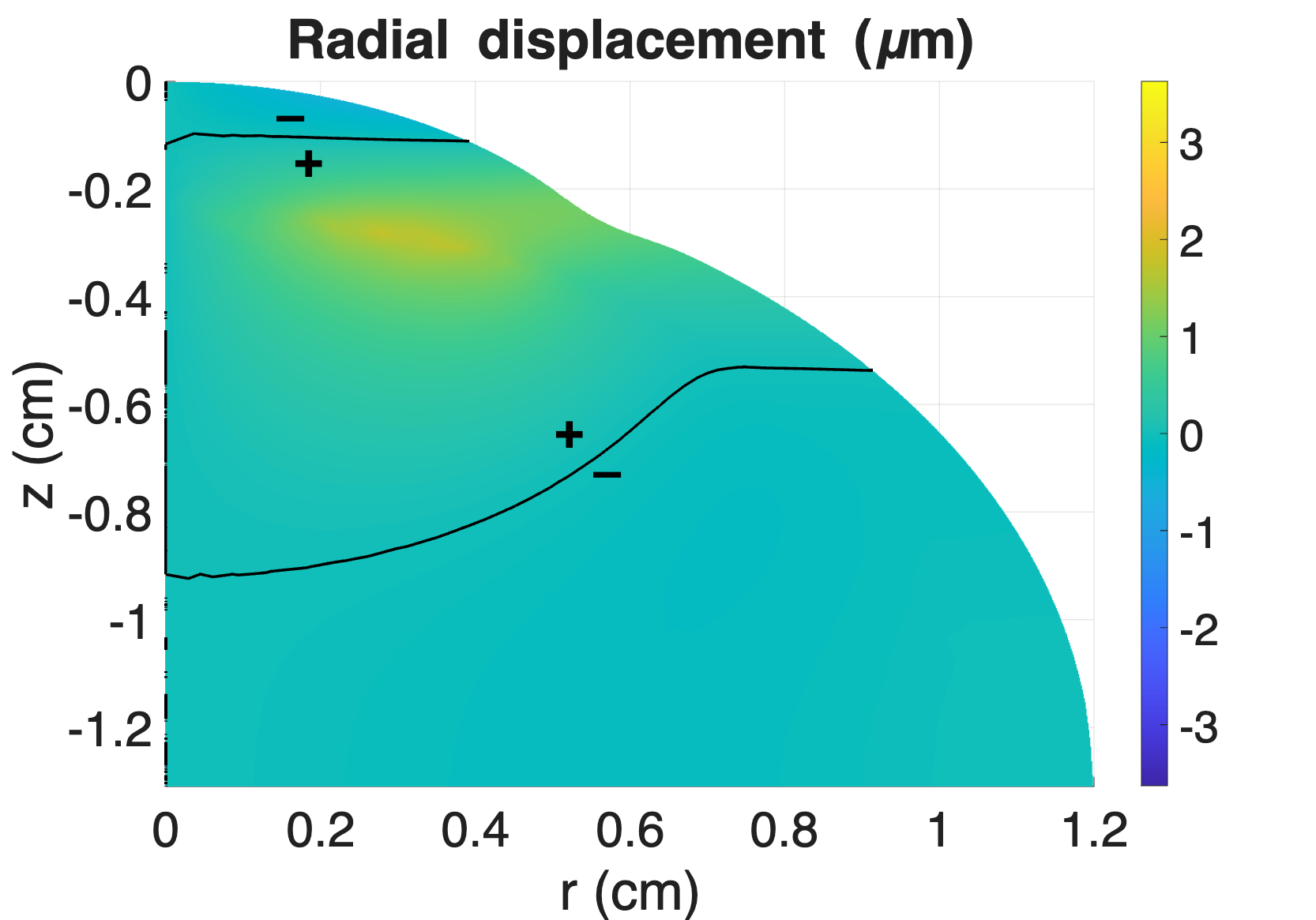}};
\coordinate [label={left:\textcolor{black}{\bf \large D}}] (D) at (-2.2,-2);
\node at (5.5,-4) {\includegraphics[width=0.33\linewidth]{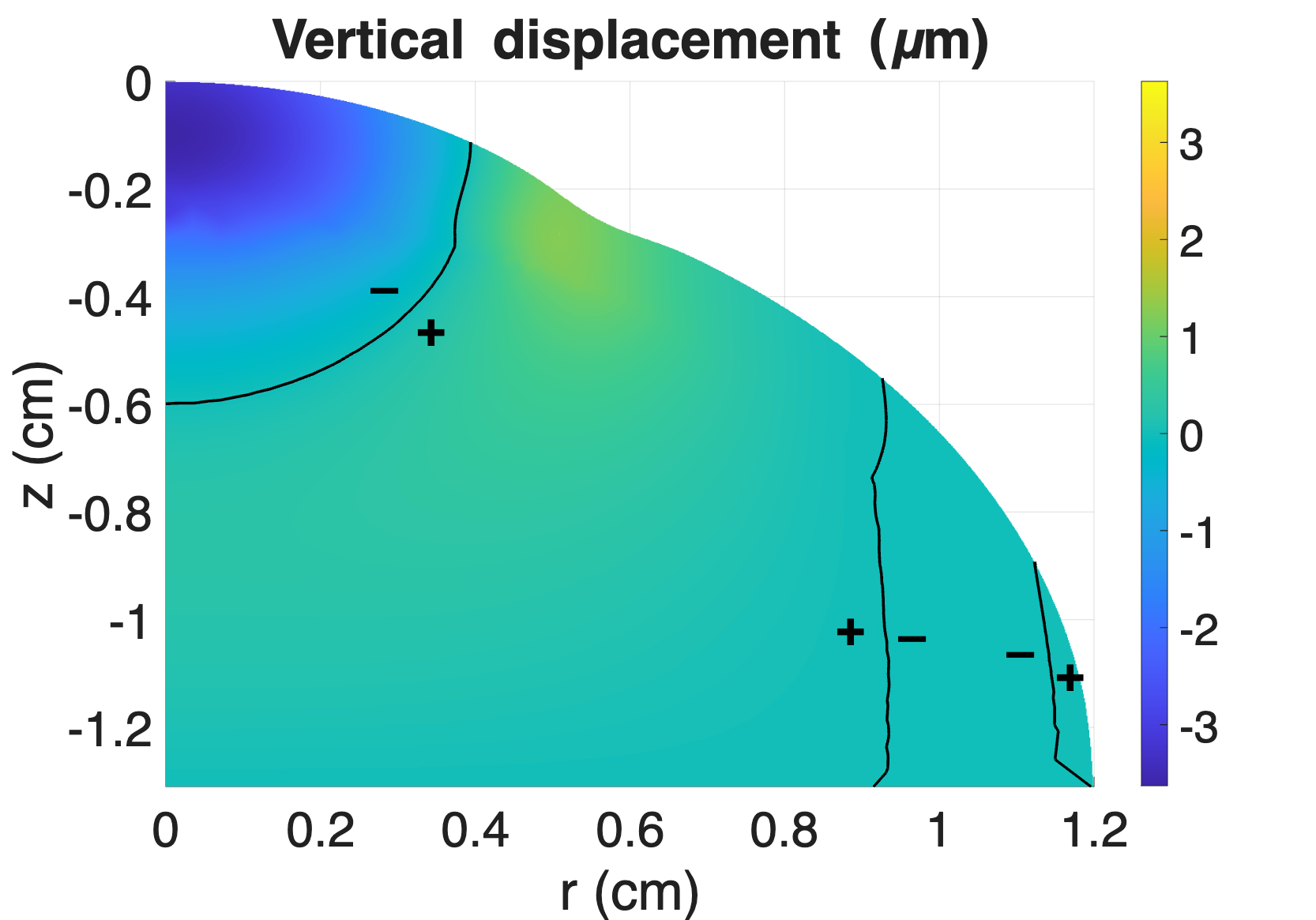}};
\coordinate [label={left:\textcolor{black}{\bf \large E}}] (E) at (3.5,-2);
\node at (11.2,-4) {\includegraphics[width=0.33\linewidth]{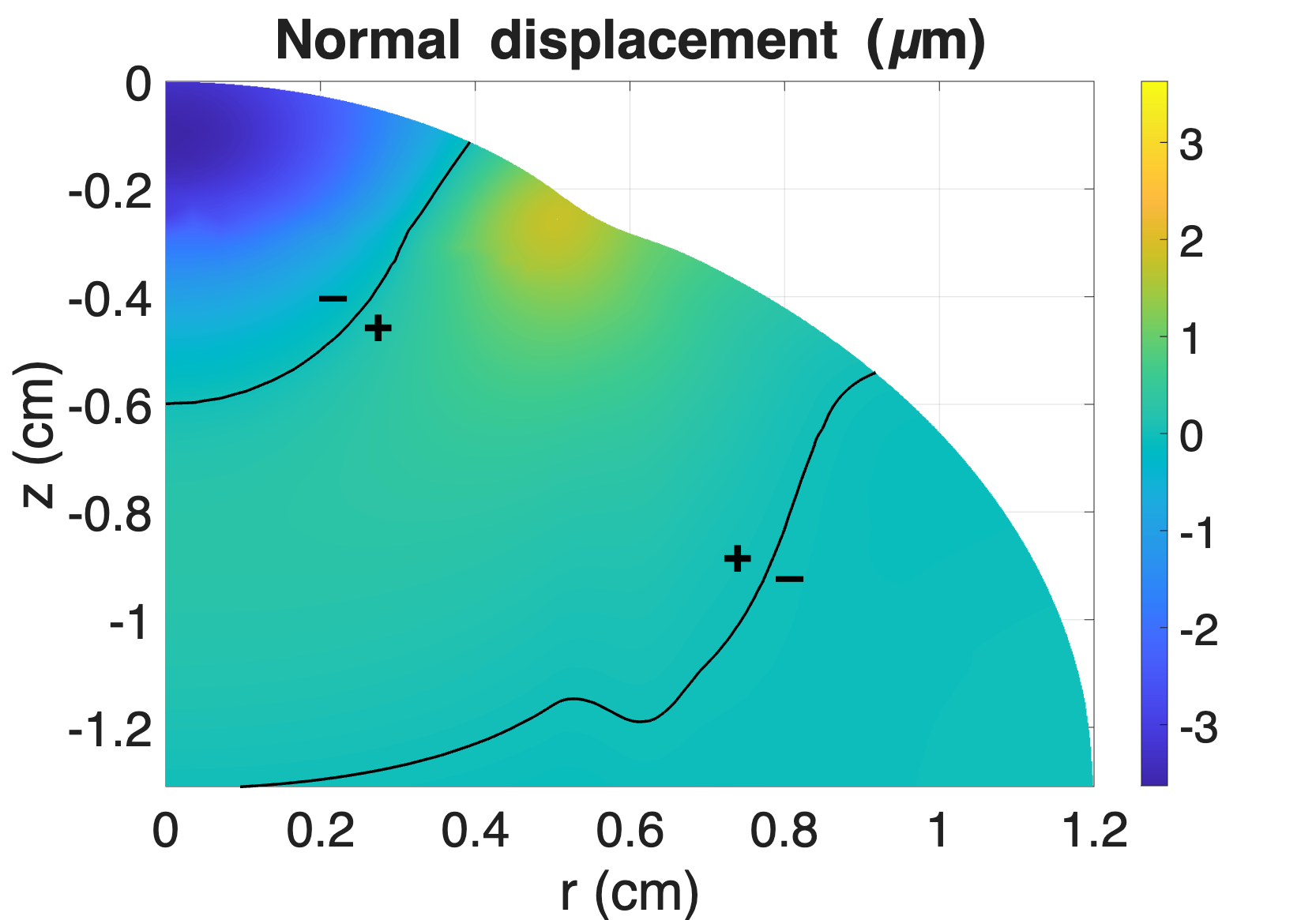}};
\coordinate [label={left:\textcolor{black}{\bf \large F}}] (F) at (9.2,-2);
\end{tikzpicture}}
\end{center}
\caption{
\textit{Heterogeneous average-shaped eye and average-shaped contact lens of constant thickness.} ({\bf A}) Contact lens displacements. ({\bf B}) Contact lens suction pressure, see Eq~($\ref{eq:p}$). ({\bf C}) Ocular surface displacements. ({\bf D}) The radial eye displacements. ({\bf E}) The vertical eye displacements. ({\bf F})~Eye displacement in the outward normal direction to the ocular surface. The solid black lines in the contour plots {\bf D}--{\bf F} are the zero level curves \rev{and the positive or negative signs indicate the sign of the displacement near the level curves}.}
\label{fig:standard_YM}
\end{figure}
Figure~\ref{fig:standard_YM} displays the predicted contact lens mechanics and eye deformations for the heterogeneous eye model. The contact lens mechanics shown in Figures\ref{fig:standard_YM}{{\bf A}--{\bf B}} are similar to the homogeneous eye model (constant Young's modulus) results presented in Figure~\ref{fig:standard}{\bf {\bf A}--{\bf B}}.
\rev{The largest relative difference in the suction pressure in the homogeneous vs the heterogeneous model is $0.02\%$ at the center of the eye $r=0$~cm.} The main changes introduced by a spatially-depend\rev{ent} Young's modulus are shown in the ocular displacement plots of Figure~\ref{fig:standard_YM}{\bf C}--{\bf F}. Qualitatively, the displacement on the ocular surface for the heterogeneous model (Figure~\ref{fig:standard_YM}{\bf C}) varies along the radial coordinate $r$ similar to the homogeneous model (Figure~\ref{fig:standard}{\bf C}), except at the edge of the lens ($r=0.7$~cm), where assuming a spatially-dependent Young's modulus in the limbus smooths out the effect of the suction pressure on the ocular surface (no local minimum or maximum of the surface displacement shown in Figure~\ref{fig:standard_YM}{\bf C} at $r=0.7$~cm compared to Figure~\ref{fig:standard}{\bf C}). The model predicts peaks of $||\vector{u}||$ $6$ times higher in the corneal and limbus regions 
when comparing the heterogeneous model to the homogeneous model. 
Similarly, the plots of the radial, vertical, and normal displacement\rev{s} on the reference ocular domain $\Omega$ in Figures~\ref{fig:standard_YM}{\bf D}--{\bf F} have similar contour plots when compared to the homogeneous model results (Figures~\ref{fig:standard}{\bf D}--{\bf F}), with a displacement up to $7$ times higher in magnitude for the heterogeneous model scenario. \rev{The normal displacement for the interior points of the domain $\Omega$ in Figure~\ref{fig:standard_YM}{\bf F} is computed in the same manner described for Figure~\ref{fig:standard}{\bf F}.}
Note that accounting for a very small Young's modulus in the center of the eye reduces how much the negative vertical and normal displacements in the center of the eye extends inside the ocular domain of $0.2$~cm, since the contour line that separates the cornea (negative displacement) to the limbus (\rev{positive} displacement) extends to $-0.8$~cm in the homogeneous eye model and only to $-0.6$~cm in the heterogeneous eye model results.
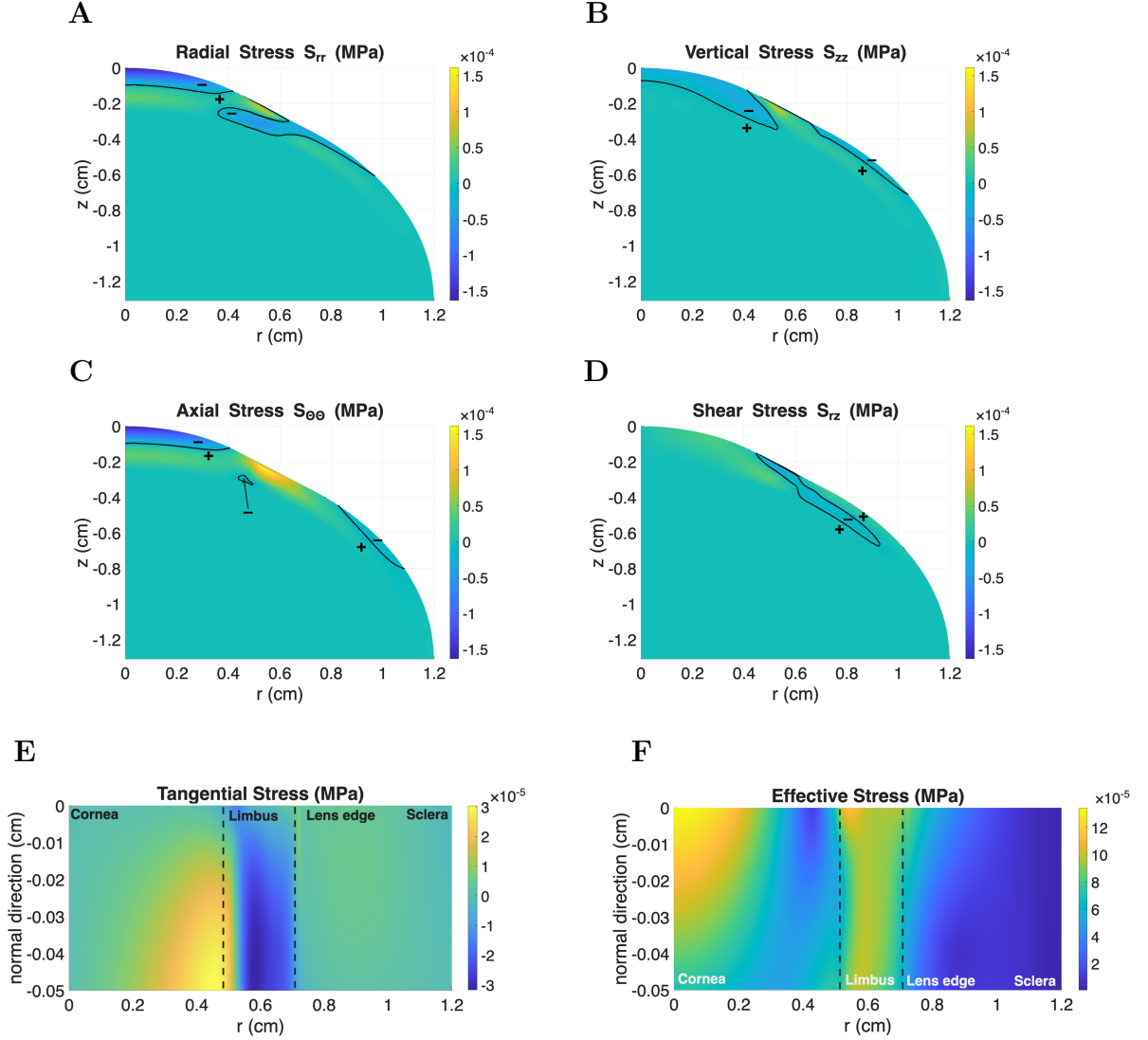
\begin{figure}[h!]
\begin{center}
\begin{tikzpicture}
\node at (-.2,0) {\includegraphics[width=0.4\linewidth]{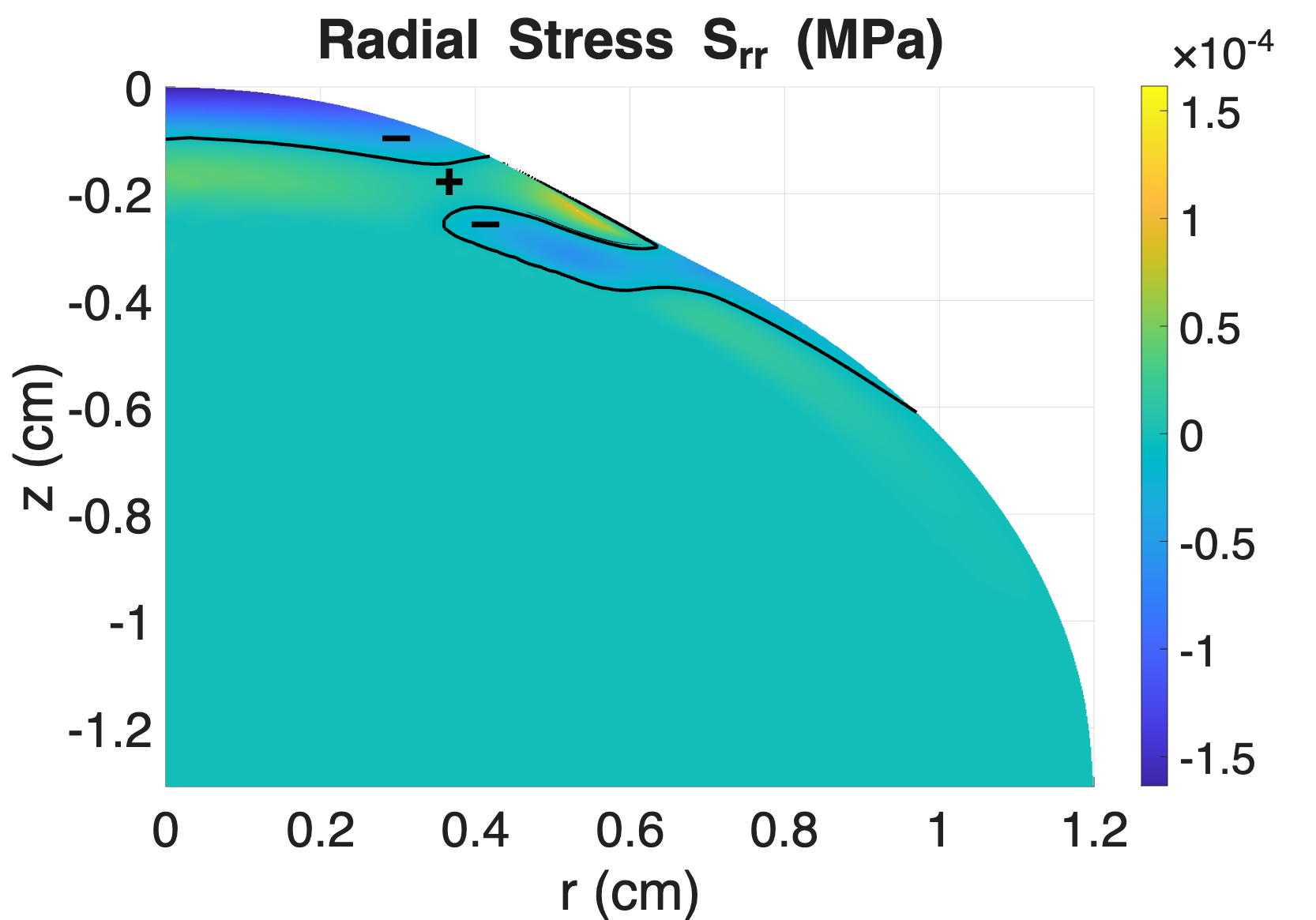}};
\coordinate [label={left:\textcolor{black}{\bf \large A}}] (A) at (-2.8,2.5);
\node at (7,0) {\includegraphics[width=0.4\textwidth]{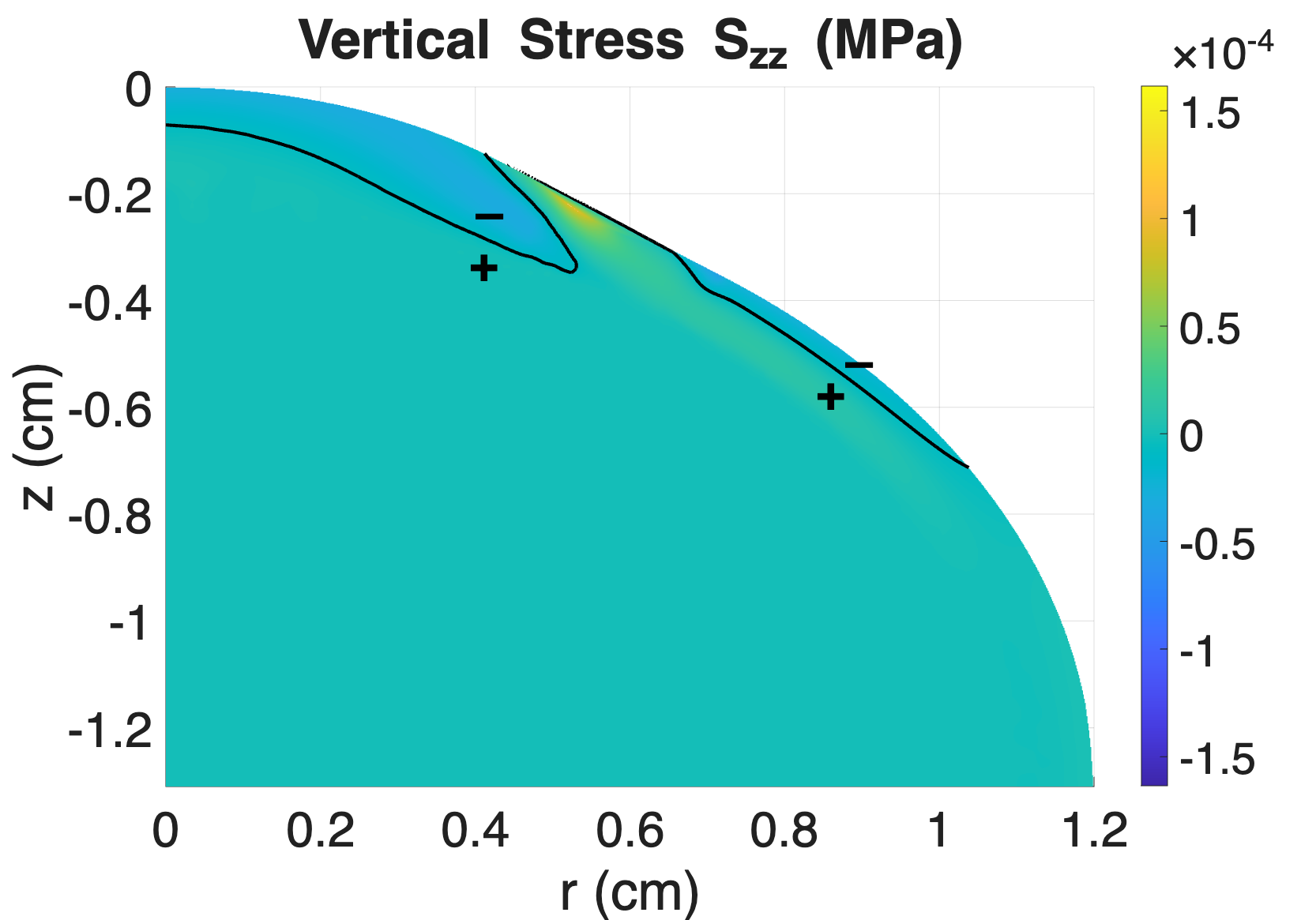}};
\coordinate [label={left:\textcolor{black}{\bf \large B}}] (B) at (4.4,2.5);
\node at (-.2,-5) {\includegraphics[width=0.4\textwidth]{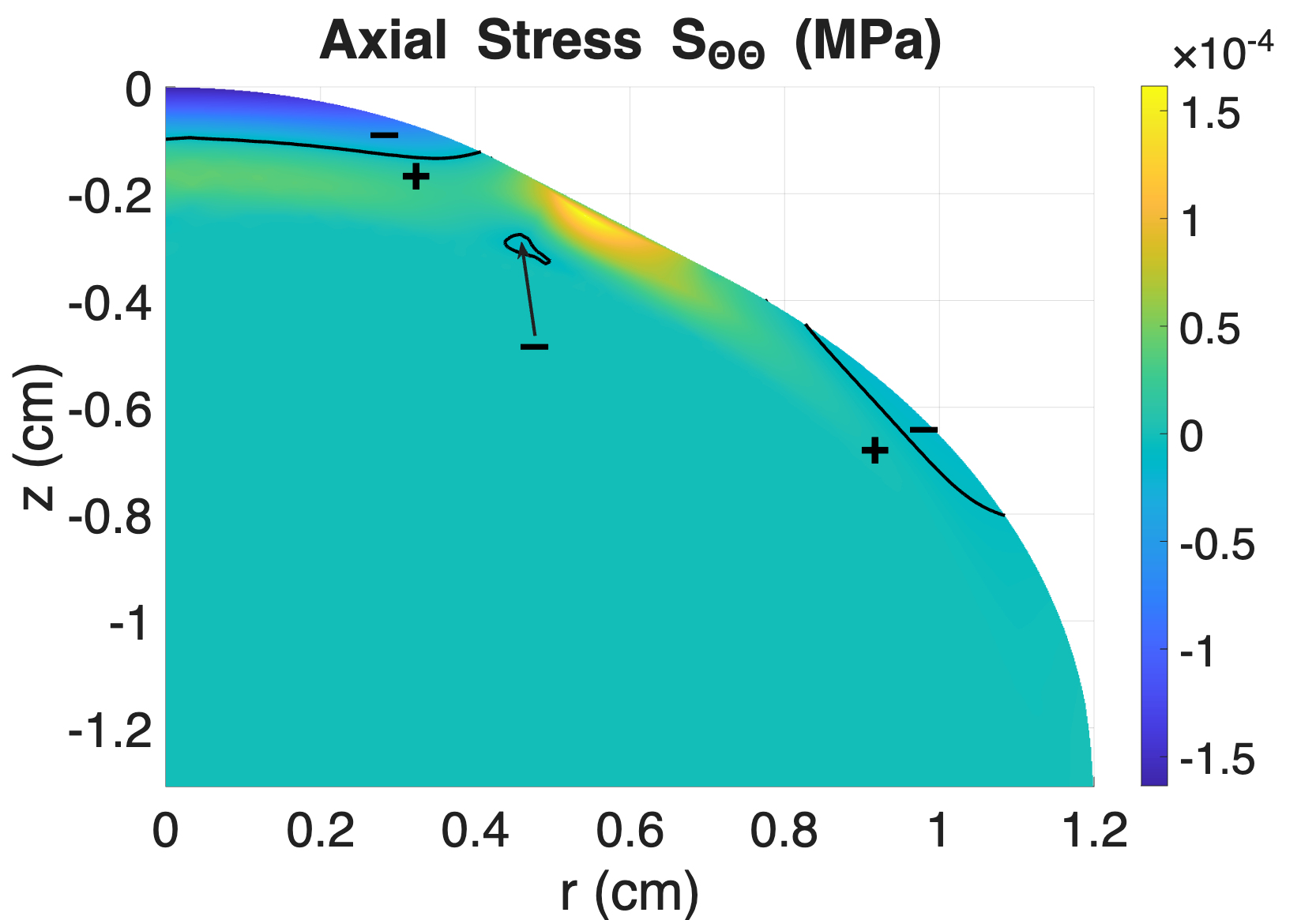}};
\coordinate [label={left:\textcolor{black}{\bf \large C}}] (C) at (-2.8,-2.5);
\node at (7,-5) {\includegraphics[width=0.4\linewidth]{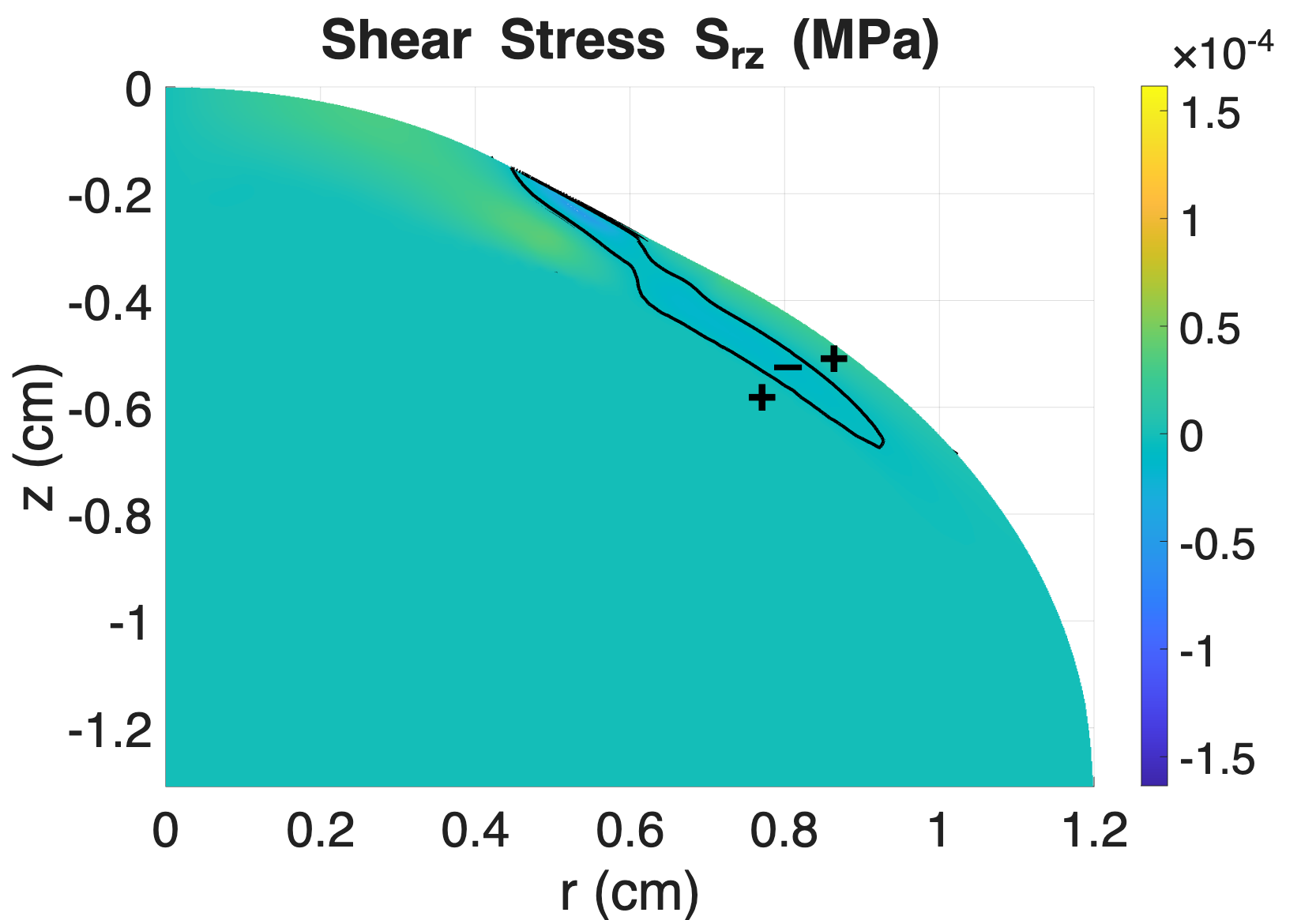}};
\coordinate [label={left:\textcolor{black}{\bf \large D}}] (D) at (4.4,-2.5);
\node at (-.5,-10) {\includegraphics[width=0.49\linewidth]{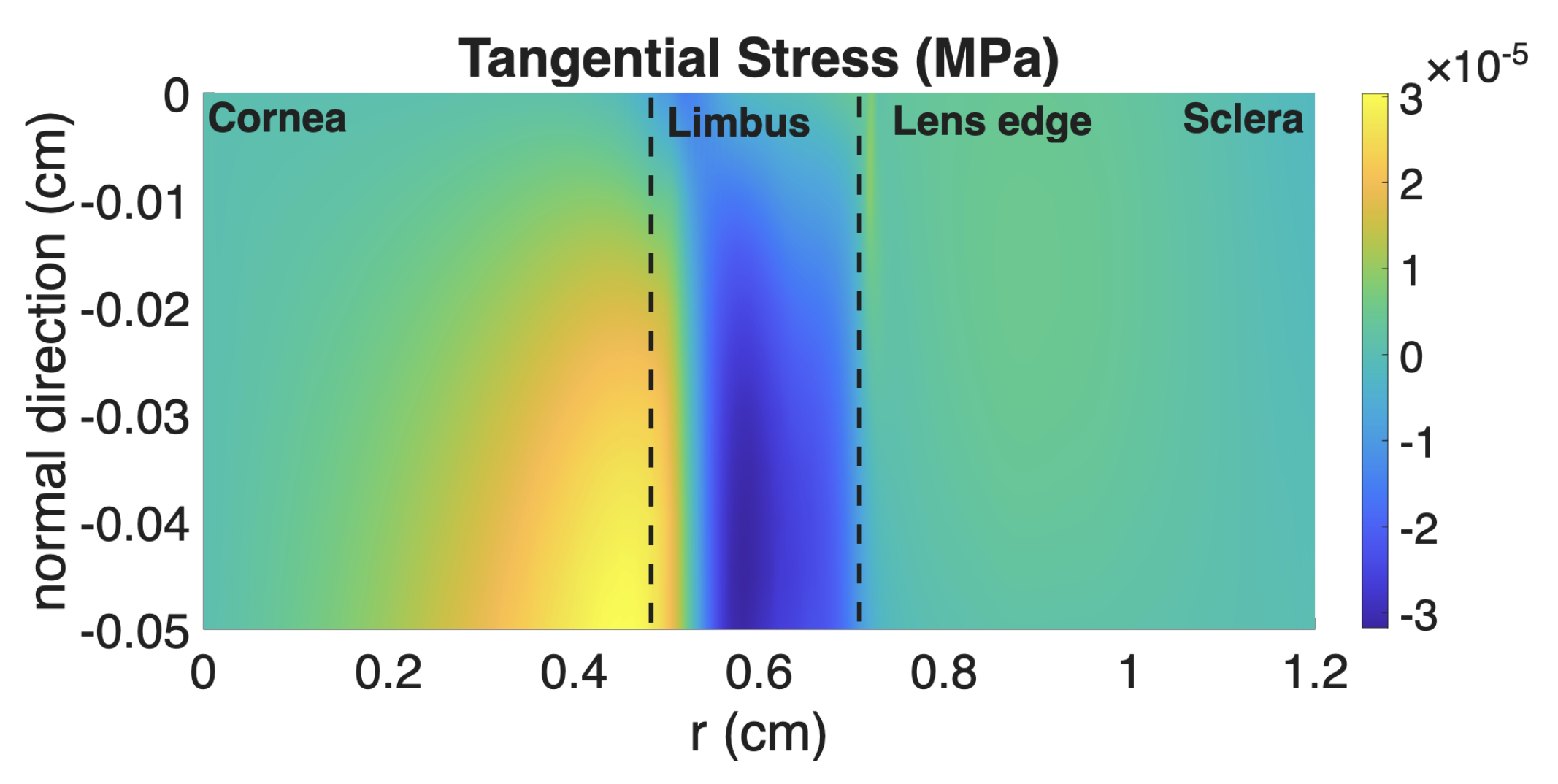}};
\coordinate [label={left:\textcolor{black}{\bf \large E}}] (E) at (-3.6,-7.8);
\node at (8,-10) {\includegraphics[width=0.49\linewidth]{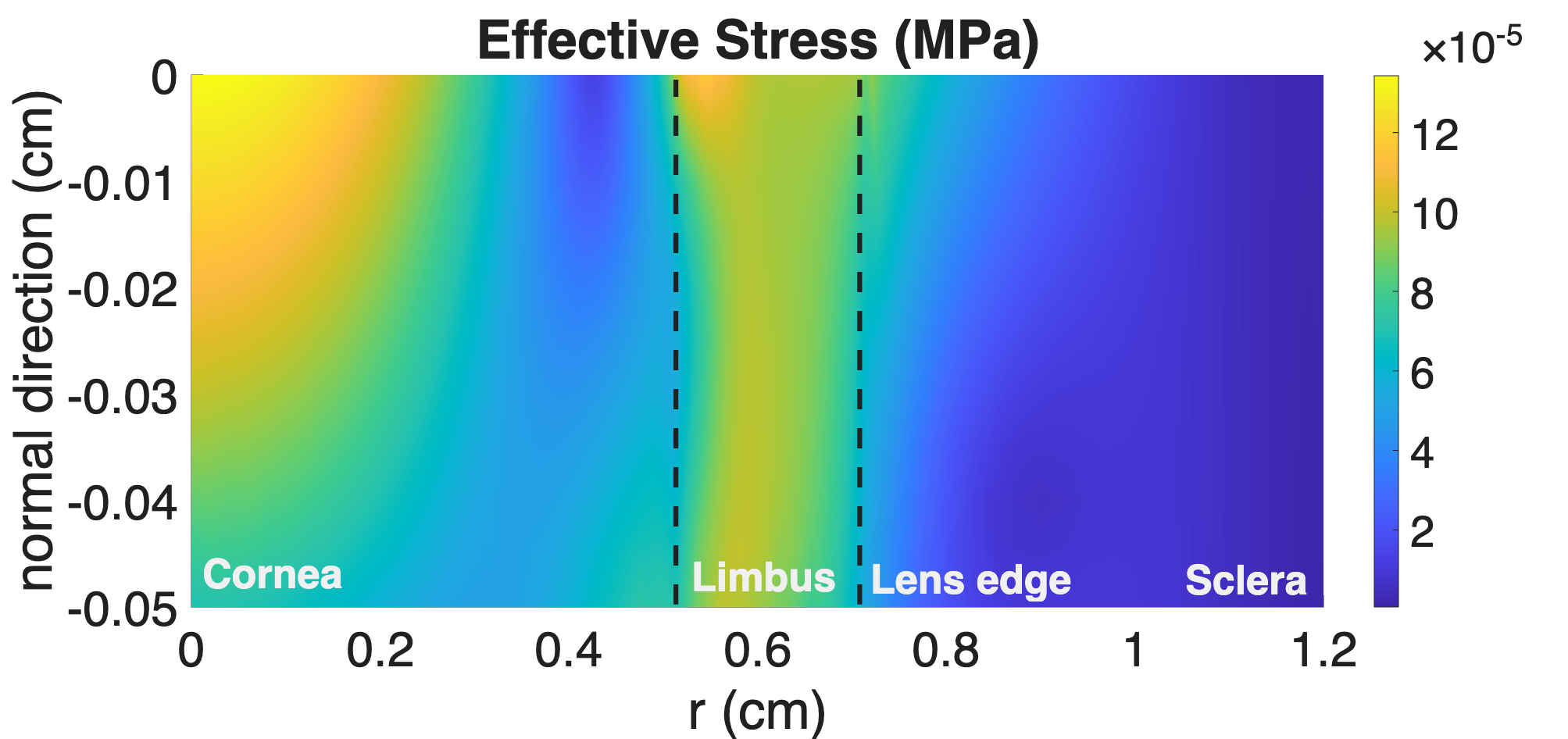}};
\coordinate [label={left:\textcolor{black}{\bf \large F}}] (F) at (5,-7.8);
\end{tikzpicture}
\end{center}
\caption{\textit{Heterogeneous average-shaped eye and average-shaped contact lens of constant thickness.} Ocular stresses:
({\bf A}) radial stress\rev{,} $S_{rr}$; ({\bf B}) vertical stress\rev{,} $S_{zz}$; ({\bf C}) axial stress\rev{,} $S_{\theta \theta}$; and ({\bf D}) shear stress\rev{,} $S_{rz}$. ({\bf E}) The tangential stress $S_{tan}$  and ({\bf F}) the effective stress $S_{eff}$  are shown in the strip of tissue normal to the reference ocular surface $\Gamma_{out}$ (at normal direction value zero) of thickness $0.05$~cm. Note that the ranges of the color bars are the same in {\bf A}--{\bf D}, but are different in {\bf E} and {\bf F}. The solid black lines in the contour plots {\bf A}--{\bf D} are the zero level curves \rev{and the positive or negative signs indicate the sign of the stress near the level curves}.}
\label{fig:stresses_YM}
\end{figure}

Interestingly, in the heterogeneous eye model results, the eye deforms more even though, the suction pressure is equivalent to the homogeneous eye model. Recall that the equation of suction pressure, Eq~\eqref{eq:p}, only depends on the derivatives of the deformed ocular shape $H$. \rev{Even if the heterogeneous eye is deforming more, t}he maximum difference in the first and second derivatives of $H$ in the homogeneous vs the heterogeneous model is $2\times10^{-4}$ and $1\times10^{-2}$, respectively.

Figure~\ref{fig:stresses_YM} displays the predicted ocular stresses for the heterogeneous eye model. In Figures~\ref{fig:stresses_YM}{\bf A}--{\bf D}, the higher stresses are experienced in the cornea and in the limbus, similarly to the homogeneous eye model results (Figures~\ref{fig:stresses}{\bf A}--{\bf D}). In the heterogeneous eye model, such high stresses are limited to the outer region of the eye (see Figure~\ref{fig:simplified_varied_ym}) and do not extend further inside the ocular domain. In the heterogeneous eye model, the eye experiences stresses up to $4$ times higher compared to the homogeneous model. Figures~\ref{fig:stresses_YM}{\bf E}--{\bf F} show the tangential and effective stresses in \rev{the outer} strip of tissue \rev{shown in Figure~\ref{fig:simplified_varied_ym}, which is} normal to the reference ocular surface $\Gamma_{out}$ (at \rev{zero} normal direction \rev{on the vertical axis}) of a thickness of $0.05$~cm. 
Qualitatively, the ocular surface tangential stress in the heterogeneous model, shown in Figure~\ref{fig:stresses_YM}{\bf E}, varies similarly along the coordinate $r$ compared to the homogeneous model, shown in Figure~\ref{fig:stresses}{\bf E}. The effective stress in the heterogeneous model, shown in Figure~\ref{fig:stresses_YM}{\bf E}, reaches its highest values not only in the limbus region but also at the center of the eye. This is not the case for the effective stress predicted in the homogeneous model, which was predicted to be low at the center of the eye, see Figure~\ref{fig:stresses}{\bf F}.

\section{Discussion}
\label{sec:discussion}
For the different ocular and lens shapes and lens thickness\rev{es} considered, our model predicts a lens suction pressure in the range $[-8,6]\times 10^{-5}$~MPa that is positive near the center of the eye and near the edge of the lens, and negative in the limbal region. The predicted suction pressure magnitude and variations along the ocular surface are in agreement with the previous modeling efforts of Maki and Ross~\cite{maki2014new,ross2016existence}, who proposed the contact lens \rev{mechanics} model used here, where the eye is modeled as a rigid solid.

Among the previous modeling efforts that account for a deformable eye, in~\cite{wu2024fea}, the authors reported a distribution of pressure on the ocular surface, see Figure~7(\textbf{f}) in~\cite{wu2024fea}; however, no details of how such pressure was computed were reported in the manuscript. We remind the reader than in~\cite{wu2024fea}, the authors modeled the closed eye interactions between the eye and ortho-k lenses; recall that ortho-k lenses are rigid lenses with a Young's modulus of $1500$~MPa~\cite{wu2024fea} which have a different shape when compared to soft contact lenses. In this work, the highest value of $E_{lens}$ considered was $10$~MPa for $\mathrm{E}=50$. Qualitatively, if we compare Figure~7(\textbf{f}) in~\cite{wu2024fea} to the suction pressure profiles computed in this work, we see that both pressures have areas of higher pressures at the center of the eye and at the edge of the lens. The magnitude of the pressure reported in Figure~7(\textbf{f}) in~\cite{wu2024fea} are in the range $[3,10]\times 10^{-2}$~MPa, one to two order of magnitude higher to the suction pressure predicted by our model for $\mathrm{E}=50$ considered (approximately 100 times the \rev{maximum absolute} suction pressure for $\mathrm{E}=1/2$ as shown in Figure~\ref{fig:material_props}\textbf{B}). \rev{The other theoretical studies we compared our work to did not report or estimate the pressure on the ocular surface when the lens in placed on the eye, see Table~\ref{tab:comp_models_results}.}

In this work, the maximum pressure exerted by the lens on the ocular surface for all the ocular and lens shapes and parameters considered is lower or in the range of the experimental measurements of eye lid pressure (i.e., the pressure that the eye lid exerts on the cornea during blinking). The average eye lid pressure is $16$~mmHg~=~$2.1\times 10^{-3}$~MPa~\cite{shaw2010eyelid}, and such pressure can vary with age between $5$~mmHg~=~$6.6\times 10^{-4}$~MPa to $30$~mmHg~=~$4.0\times 10^{-3}$~MPa~\cite{yamaguchi2018relationship}. The magnitude of the suction pressure predicted by this work if $\mathrm{E}=1/2$ is at most $4\times 10^{-5}$~MPa, thus one order of magnitude smaller than the minimum eyelid pressure measured in~\cite{yamaguchi2018relationship}. As the contact lens becomes stiffer (i.e., as $\mathrm{E}$ increases), the suction pressure magnitude gets closer to measured range of the eyelid pressure. If $\mathrm{E}=50$, then the maximum magnitude of the predicted suction pressure on the eye is equal to the maximum value of the eyelid pressure measured in~\cite{yamaguchi2018relationship}.

Ramasubramanian et al.~\cite{ramasubramanian2024influence} developed a closed eye model with soft contact lenses, and the third principal stress ($\sigma_3$) in the cornea was proposed as an alternative to estimate the contact pressure between the eye and the lens. The authors considered different contact lens Young's modulus from $0.2$ to $1.44$ MPa. They estimated high values of $\sigma_3$ at the center of the cornea, which is similar to the suction pressure predicted in the work, and low values of $\sigma_3$ near the edge of the lens, which is in disagreement with the current work. Recent modeling efforts which studied the interaction between an ortho-k (stiff) lens and the closed eye~\cite{wu2021biomechanical,zhao2023biomechanical,wu2024fea} all predicted an increase in the effective stress in the eye near the edge of the lens, which is similar to this work. These models predict different ranges for the maximum effective stress in the cornea: between $1\times10^{-2}$ to $2\times10^{-2}$~MPa, depending on different corneal thickness, curvature, and refractive change in~\cite{wu2021biomechanical}; between $2.1\times10^{-3}$ and $2.4\times10^{-3}$~MPa, depending on lens refractive power and myopia degree in~\cite{zhao2023biomechanical}; and $1\times 10^{-1}$~MPa in~\cite{wu2024fea}. For the stiffer contact lens considered in this work, namely \rev{$E_{lens}=10$~MPa which corresponds to} $\mathrm{E}=50$, our model predicts a maximum effective stress in the eye of $2\times 10^{-3}$~MPa, similar to what was predicted by ~\cite{zhao2023biomechanical} and lower than what was predicted in~\cite{wu2021biomechanical,wu2024fea}. \rev{We summarize the comparison between this work and \cite{wu2021biomechanical,zhao2023biomechanical,wu2024fea,ramasubramanian2024influence} in Table~\ref{tab:comp_models_results}. For details on the differences between our model and the models in \cite{wu2021biomechanical,zhao2023biomechanical,wu2024fea,ramasubramanian2024influence}, please prefer to Table~\ref{tab:comp_models}.} Note in~\cite{wu2021biomechanical,zhao2023biomechanical,wu2024fea}, the contact lens Young's modulus considered was higher than $100$~MPa.

Ocular displacements due to contact lens wear have been reported in the literature. For example, Ramasubramanian et al.~\cite{ramasubramanian2024influence} (modeling a soft contact lenses interacting with a closed eye) estimated a vertical displacement of the eye between $15$ to $400$~$\mu$m, which is larger than the maximum magnitude of ocular displacement predicted in our heterogeneous eye model scenario of approximately $4$~$\mu$m (see Figure~\ref{fig:standard_YM}\textbf{C}). This difference might be due to the different assumptions between this work and~\cite{ramasubramanian2024influence}; Ramasubramanian et al.~assumed that the eye lid is closed and the lens does not have to conform fully to the eye shape\rev{, see Table~\ref{tab:comp_models}}. The ocular displacements at the center of the cornea predicted in this work are in agreement with the experimental measurements reported by Alonso-Caneiro et al.~\cite{alonso2012using}. In~\cite{alonso2012using}, the authors measure\rev{d} the changes of the ocular surface in the morning vs after wearing a soft contact lens for 6 hours, and they found a change of approximately $4$~$\mu$m in the corneal region and  
$12$~$\mu$m in the scleral/limbal region. Our model predicts a smaller magnitude of the displacement in the limbal region, of the order of $2$~$\mu$m in the heterogeneous eye model. 

\begin{table}[h!]
    \centering
    \footnotesize
     \caption{\rev{Comparison between the prediction of the current work (first column) and previous modeling efforts. The symbol - indicates if a quantity was not reported, not used or not estimated.}}
    \begin{tabular}{l lllll}
    \hline
         & This work & \cite{wu2021biomechanical}& ~\cite{zhao2023biomechanical} & ~\cite{wu2024fea} & ~\cite{ramasubramanian2024influence} \\         
         \midrule 
         Young's modulus~[MPa] & & & &\\
         \quad Cornea  & $0.2$& $[0.1,1.3]$ & $[0.1,1.3]$ & $[0.2,0.3]$& -\\ 
    \quad Contact lens& $[0.1,10]$  &  $100$ & $1160$ & $1500$ & $[0.2,1.44]$ \\
    \hline
      Pressure on the eye & & & & & \\
      \quad Maximum~[MPa] & $[0.04,4]\times 10^{-3}$ & - & - & $[3,10]\times 10^{-2}$ & -  \\
      \quad High pressure near the center of the eye & Yes  & - & - & Yes & -\\
      \quad High pressure near the edge of the lens & Yes  & - & - & Yes & -\\           
      \hline
      Ocular stresses & & & & & \\
      \quad Maximum effective stress [MPa] & $[0.02,2]\times10^{-3}$ & $[1,2]\times10^{-2}$ & $[2.1,2.4]\times10^{-3}$ & $1\times 10^{-1}$ & -\\
      \quad High stresses near the center of the eye & Yes  & Yes & No & No & Yes\\
      \quad High stresses near the edge of the lens & Yes  & Yes & Yes & Yes & No\\    
      \hline
    \end{tabular}
   
    \label{tab:comp_models_results}
\end{table}

\indent \rev{The difference} in limbal/scleral \rev{deformations} between our model and the experimental results reported in~\cite{alonso2012using} might be due to the model assumption that the contact lens fully conforms to the eye also in the limbal region. The contact lens might not necessarily fully conform to the whole ocular surface, but might vault in the limbal region~\cite{shen2011characterization}, as captured by some of other modeling efforts~\cite{ramasubramanian2024influence}.  In further work, capturing a more detailed interaction between the eye and the lens in the limbal region is important to study the relationship between lens and eye interactions on contact lens comfort, since contact lens induced ``limbal indentations" have been observed in the clinical setting for tight-fitting soft contact lenses~\cite{alonso2012using}.\\
\indent The model results show that ocular deformations and stresses are highly dependent on the assumptions made on the ocular material parameters. Deriving the eye model, we assumed that the inside of the eye, which is composed of \rev{a} gel-like fluid called the vitreous humor, deforms either like the cornea (homogeneous model) or like an elastic material with a very small Young's modulus (heterogeneous model). In both ocular models considered, we neglect to consider that the aqueous humor and vitreous humor apply a pressure on the internal surface of the cornea and of the sclera called IOP~\cite{murgatroyd2008intraocular}. The baseline values of IOP are between $10$~mmHg~$=1.3\times 10^{-3}$~MPa and $20$~mmHg~$=2.6\times 10^{-3}$~MPa, of the same order of magnitude of the measured eyelid pressure in~\cite{shaw2010eyelid}, and an order of magnitude higher or equal to the lens suction pressure predicted by our model, depending on the stiffness of the contact lens. Further studies are required to investigate how the IOP affects the interaction between the eye and the contact lens.\\
\indent 
In this work, we modeled the limbus as a transition zone between the corneal and the scleral regions that guarantees that the first and second derivatives of the ocular surface are continuous; this is an important ocular surface requirement to guarantee the stability of the staggered algorithm that couples the lens and eye \rev{mechanics}. We note that the coupling is via the suction pressure and such pressure is computed solely using the derivatives of the deformed ocular surface $H$ in Eq~\eqref{eq:p}, which results in the sensitivity to the derivatives of the ocular surface. Additionally, we assumed axial symmetry of the ocular surface; thus, we neglect the difference in the corneosclearl junction angle between the nasal and temporal sides of the eye. The corneosclearl junction angle has been found to be the sharpest on the nasal side and almost $180^\circ$ at the temporal side (i.e., cornea and sclera are almost tangent~\cite{Halletal2011}). These differences have been found to be correlated with the contact lens fit, in particular, with contact tightness~\cite{Halletal2011}. Further studies are required to investigate how reducing the axial symmetric assumptions of the lens and eye models affect the lens suction pressure, the ocular deformation, and ocular stresses. Relaxing this assumption will increase the numerical complexity of the model, but, at the same time, allow one to investigat\rev{e} the effect of ocular shape irregularities in healthy individuals (accounting for the differences in the limbus) or keratoconus patients (that have non-symmetrical coned shapes corneas) on contact lens wear.\\
\indent 
\rev{In this work, we considered soft to stiff contact lenses, thereby varying the lens Young's modulus from $0.1$ to $10$~MPa. Rigid contact lenses, such as ortho-k lenses, have been modeled using a lens Young's modulus between $100$ and $1500$~MPa~\cite{wu2021biomechanical,zhao2023biomechanical,wu2024fea}. Future studies are required to investigate our model behavior as $E_{lens}$ increases and possibly approaches infinity (rigid contact lens). It would be interesting to explore this limiting scenario from the analytical and numerical points of view, and further the comparison of our model predictions with previous literature on ortho-k modeling~\cite{wu2021biomechanical,zhao2023biomechanical,wu2024fea}.}\\
\indent\rev{The eye model developed only accounts for the isotropic linear elastic response of the ocular tissue to mechanical stresses (i.e., the immediate non-directional dependent deformation response). However, experimental evidence shows that the ocular tissue, when subjected to a load, exhibits non-linear, anisotropic, and poro-viscoelastic response to mechanical stresses~\cite{Wooetal1972,guo2015femtosecond,ashofteh2020characterization}. Cornea's collagen fiber orientation and organization play a crucial role in the corneal ability to resist deformations~\cite{boote2005lamellar}, and the fluid-like component of the corneal tissue plays an important role when studying interactions between the eye and ortho-k lenses~\cite{swarbrick2006orthokeratology,choo2008morphologic}. Previous modeling efforts have considered the non-linear response of the ocular tissue to stresses~\cite{ramasubramanian2024contact}. In future work, it would be interesting to extend the current eye model not only to account for the effect of the IOP and for a non-axially symmetric load, but also to account for a more realistic non-linear response of the ocular tissue to stresses. In particular, accounting for the ocular tissue non-linear response is important to capture the interaction between the eye and ortho-k lenses, where we expect the eye to experience larger stresses than soft contact lenses.}\\
\indent
In between the contact lens and the eye, there is a thin layer of tear film (post-lens layer), with a thickness of the order of $3-4$~$\mu$m~\cite{wang2003precorneal}. Multiple previous modeling efforts, including this work, have either neglected the effect of the tear film \rev{flow} dynamics on the contact lens and eye interactions~\cite{maki2014new,wu2021biomechanical}, or accounted for a constant tear fluid surface tension in between the lens and the eye~\cite{wu2024fea,ramasubramanian2024influence,moore2019simulation} \rev{and a nonzero sliding friction coefficient~\cite{wu2024fea}}. The ratio of the post-lens tear film thickness to the thickness of the contact lens is usually very small; in this work, the \rev{largest} value of this ratio is $3/35 \simeq 0.08$. Various mathematical models have been developed to account for a detailed description of the tear film dynamics in the post- or pre-lens tear film\rev{s}~\cite{ramasubramanian2024contact}, and to the best of our knowledge, all of these model assume that the eye is rigid. In future studies, it would be interesting to incorporate the tear film dynamics into the current model and study how the lens suction pressure and ocular deformations would be affected by the presence of the post-lens tear film \rev{flow}. \rev{Accounting for the post-lens tear film flow in the model, the contact lens would conform to the shape of the eye plus a non-constant tear film thickness. Maki and Ross showed that the thickness of the tear film under a soft contact lens can vary in the radial direction due to the lens suction pressure and presented peaks near the limbal region~\cite{maki2014exchange}. Incorporating the post-lens tear film dynamics could reduce the difference between the current model predictions and the experimental measurements of ocular displacement in the limbus~\cite{alonso2012using}.}

\section{Conclusions}
This paper presented the first mathematical model that \rev{captures} the interaction between the lens and the open eye. \rev{In the model}, the contact lens configuration, the contact lens suction pressure, and the deformed ocular shape are all emerg\rev{ing} properties of the model \rev{and were not imposed a priori}. \rev{The contact lens and the eye were} non-linear\rev{ly} coupl\rev{ed} by assuming that the suction pressure under the lens is directly applied to the ocular surface. We modeled the contact lens \rev{mechanics} using a previous published model~\cite{maki2014new,ross2016existence} \rev{and} considered homogeneous and heterogeneous linear elastic eye model\rev{s.}

\rev{We presented model results for} different ocular shapes, different lens shapes, and lens thickness profiles, and extracted lens deformations, lens suction pressure, ocular deformations, and ocular stresses for all the scenarios considered. The model predicted \rev{that the regions where the} higher ocular deformations and stresses \rev{happen are} the center of the eye and the limbal region. The ocular displacements and stresses non-linearly increase as we increase the stiff\rev{ness} of the contact lens \rev{and as we consider heterogeneous eye material parameters. The maximum ocular effective stress predicted by our model ranged from $0.02\times10^{-3}$ to $2\times10^{-3}$~MPa, depending on the lens and eye material parameters considered.} \rev{The relative difference between the lens curvature and the eye curvature had an important effect on the model results, since} a steeper contact lens on the eye resulted in a reduction of the ocular displacement at the center of the eye and a larger displacement at the edge of the contact lens.\\
\indent
The model predictions were compared to experimental data and previously developed mathematical models. The maximum suction pressures predictions \rev{varied in between $[0.04,4]\times10^{-3}$MPa, depending on the stiff\rev{ness} of the contact lens and the eye model considered, and were either} lower or equal to the measured values of the eye lid pressure \rev{range of $[0.66,4]\times 10^{-3}$~MPa}~\cite{shaw2010eyelid}. The lens induced deformations predicted at the center of the cornea \rev{of approximately $4$~$\mu$m in the heterogeneous model} were in agreement with experimental measurements~\cite{alonso2012using}. The effective stresses predicted in the outer ocular surface were lower or equal to the one predicted by previous closed eye model of soft or stiff contact lenses~\cite{wu2021biomechanical,zhao2023biomechanical,wu2024fea}\rev{, as shown in Table~\ref{tab:comp_models_results}}.

\rev{Despite its limitations, this modeling effort gives insights on the mechanisms driving contact lens and eye interactions and the contact lens fitting process. In the future, lens manufacturers and clinicians, aided by mathematical models, may be able to theoretically predict the performance of a lens design without having to manufacture the lens and then test the lens in an expensive clinical trial, or without having to try the lens on the patient.}  

\section*{Use of AI tools declaration}
The authors declare they have not used Artificial Intelligence (AI) tools in the creation of this article.

\section*{Acknowledgments}
This material is based upon work supported by the National Science Foundation under Award No.~2316951 (PI Carichino) and by the Rochester Institute of Technology, College of Science, Dean’s Research Initiation Seed Funding Grant (PI Carichino).
The authors would like to thank Vedat Kurtay for the help in the initial stages of the project, and George Thurston and Olalekan Babaniyi for the useful discussions. 

\section*{Conflict of interest}
The authors declare there is no conflict of interest.

\appendix
\section*{Appendix}
\setcounter{table}{0}
\setcounter{figure}{0}
\renewcommand{\thefigure}{A\arabic{figure}}
\renewcommand{\thetable}{A\arabic{table}}

\section{Contact lens posterior surface and thickness}
\label{sec:supp_lens}
Funkenbusch and Benson~\cite{funkenbusch1996conformity} studied the conformity of soft contact lenses to \rev{the} eye, and we use similar lens shapes to the lenses defined in their work. The posterior curve of the lens is assumed to be an ellipse as follows:
\begin{equation}
   g(r) = -b_{lens} + b_{lens}\sqrt{1 - \frac{r^2}{K_{lens}b_{lens}}} \ \mathrm{cm}, \quad 0 \leq r < \mathfrak{R}_{lens},
\end{equation}
where 
\begin{equation}
    b_{lens} = -\frac{z_{lens}^2 K_{lens}}{2z_{lens} K_{lens}+\mathfrak{R}_{lens}^2},
\end{equation}
$z_{lens}$ is the sagittal height of the contact lens, $K_{lens}$ is the radius of curvature at $r=0$ of the lens, and $\mathfrak{R}_{lens}$ is the undeformed radius of the lens. Table~\ref{tab:lens_geometry} reports the values of the lens parameters for a flat, average, and steep lens shown in Figure~\ref{fig:eye_geometry}. 
\begin{table}[h]
    \centering
     \caption{Parameters for the contact lens posterior surface.}
    \begin{tabular}{llll }
    \hline
         & Flat lens & Average lens & Steep lens \\
         \hline
         $\mathfrak{R}_{lens}$ [cm] & 0.7000 & 0.7000 & 0.7000 \\
      $K_{lens}$ [cm]  & 0.9000 & 0.8700 & 0.8400 \\
      $z_{lens}$ [cm] & -0.3405 & -0.3630 & -0.3857 \\
      $b_{lens}$ [cm] & 0.8490 & 0.8094 & 0.7910\\
      \hline
    \end{tabular}
   
    \label{tab:lens_geometry}
\end{table}

\begin{figure}[h]
    \centering
    \includegraphics[width=0.45\linewidth]{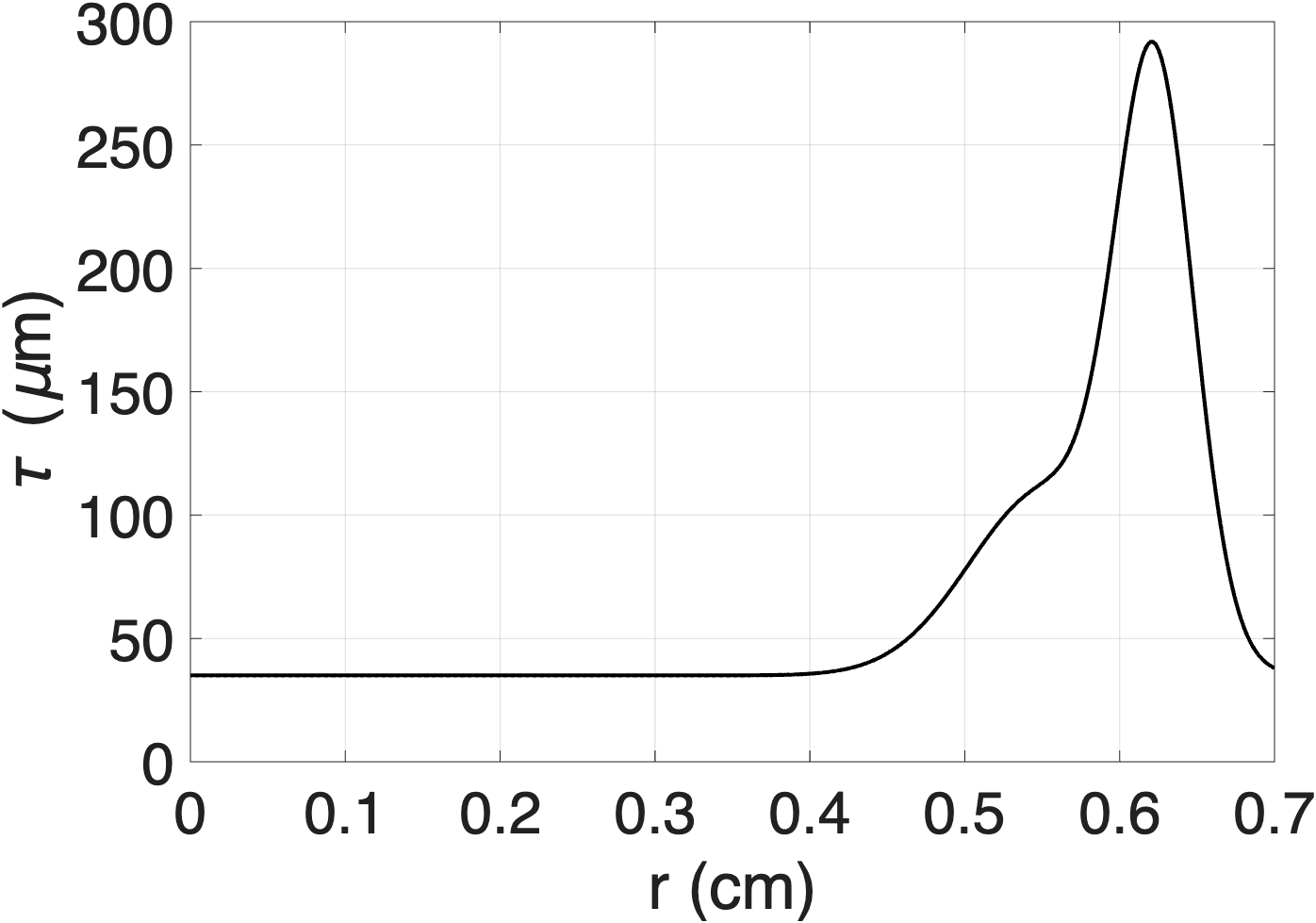}
    \caption{Spatially-dependent contact lens thickness $\tau$ as a function of the radial coordinate~$r$\rev{~\cite{funkenbusch1996conformity}}.}
    \label{fig:tau}
\end{figure}

We consider either a lens with a constant thickness or a lens with a spatially-dependent thickness profile to match the descriptions of the lens thickness profiles studied in Funkenbus\rev{c}h and Benson \cite{funkenbusch1996conformity}. The center and edge lens thicknesses are $35$~$\mu$m. The lens thickness increases to $269$~$\mu$m as the radius increases from $0$ to $0.623$~cm, and then decreases to $35$~$\mu$m at the lens edge. Symbolically,
\begin{equation}\label{eqn:thickness}
    \tau(r) = \tau_{inc}e^{- \dfrac{(r-r_{inc})^2}{2c_{inc}^2}} + \left(\tau_{max} - \tau_{inc} e^{- \dfrac{(r_{max}-r_{inc})^2}{2c_{inc}^2}} -\tau_{cent}\right)e^{- \dfrac{(r-r_{max})^2}{2c_{max}^2}}  +  \tau_{cent},
\end{equation}
where $\tau_{cent}=35$~$\mu$m, $\tau_{inc}=7.5$~$\mu$m, $\tau_{max}=269$~$\mu$m, $r_{max}=0.623$~cm, $r_{inc}=0.553$~cm, $c_{max}=0.025$~cm, and $c_{inc}=0.0.05$~cm, as shown in Figure~\ref{fig:tau}.

\section{Ocular surface}
\label{sec:supp_eye}
\begin{figure}[h]
\centering
\begin{tikzpicture}[scale=0.7]
\def\K{0.775}
\def\ee{0.320}
\def\bc{0.5270}
\def\ac{0.6391}
\def\zs{-1.314}
\def\Re{1.2}
\def\Rs{1.2}
\def\lone{24.847}
\def\ltwo{502.197}
\def\lthree{6099.372}
\def\lfour{-42.384}
\def\lfive{545.984}
\def\lsix{-5042.020}
    \draw[->,thin,color=gray] (0,-6) -- (0,.5) node[left] {$z$};
	\draw[->,thin,color=gray] (-1,0) -- (6,0) node[right] {$r$};
\scalebox{4}{
    \draw[domain=0:0.5, samples = 10, smooth, variable=\x, black,line width=0.1mm] plot ({\x},{-\bc+(\bc^2-(\bc*\x/\ac)^2)^0.5});
    \draw[domain=0.7:\Re-.01, samples=50,smooth, variable=\x, black,line width=0.1mm] plot ({\x},{\zs+(\Rs^2-\x^2)^0.5});
    \draw[domain=0.5:0.7, samples=5, smooth, variable=\x, black,line width=0.1mm] plot ({\x},{(\x-0.7)^3*(\lone+\ltwo*(\x-0.5)+\lthree*(\x-0.5)^2)+(\x-0.5)^3*(\lfour+\lfive*(\x-0.7)+\lsix*(\x-0.7)^2)});
    \draw[color=black, line width=0.1mm] ({0},{\zs+(\Rs^2-\Re^2)^0.5+0.15}) -- ({\Re},{\zs+(\Rs^2-\Re^2)^0.5+0.15});
    \draw[color=black, line width=0.1mm] ({0},{0}) -- ({0},{\zs+(\Rs^2-\Re^2)^0.5+0.15});
}

\draw[dotted,color=gray, thick] (4*\Rs,4*\zs) arc (0:360:4*\Rs);
\filldraw [gray] (0,4*\zs) circle (2pt);

\draw[dotted,color=gray, thick] (4*\ac,-4*\bc) arc(0:360:4*0.6391 and 4*0.527);
\filldraw [gray] (0,-4*\bc) circle (2pt);

\draw[dashed, color=gray] (4*0.5,0.3) -- (4*0.5,-4.7);
\draw[dashed, color=gray] (4*0.7,0.3) -- (4*0.7,-4.7);
\coordinate [label=center:\textcolor{black}{{\small $r = 0.5$}}] (llb) at (1.6,.6);
\coordinate [label=center:\textcolor{black}{{\small $r = 0.7$}}] (lub) at (3.3,.6);
\draw[dashed, color=gray] (4*1.2,0.3) -- (4*1.2,-4.7);
\coordinate [label=center:\textcolor{black}{{\small $\mathfrak{R}_{eye}$}}] (lub) at (5,.6);

\coordinate [label={[rotate=-50]center:\textcolor{black}{\small sclera}}] (a) at (4,-2.15);
\coordinate [label={[rotate=-17]center:\textcolor{black}{\small cornea}}] (a) at (1.3,-.1);
\draw[<->, color=gray] (-0.15,-.05) -- (-0.15,-5.2);
\coordinate [label=center:\textcolor{black}{{$z_{s}$}}] (a) at (-.5,-2.625);
\draw[<->, color=gray] (.1,-4*\bc) -- (4*\ac,-4*\bc);
\coordinate [label=center:\textcolor{black}{{$a_c$}}] (a) at (1.2,-4*\bc+.2);
\draw[<->, color=gray] (0.15,-4*\bc+.05) -- (0.15,-.05);
\coordinate [label=center:\textcolor{black}{{$b_c$}}] (a) at (.4,-1);
\draw[<->, color=gray] (0.1,-5.25) -- (4.75,-5.25);
\coordinate [label=center:\textcolor{black}{{$\mathfrak{R}_s$}}] (a) at (2.3,-5.55);

\coordinate [label=right:\textcolor{black}{{$h(r)$}}] (a) at (4.3,-3.2);

\end{tikzpicture}
\caption{Schematic of the construction of the ocular domain.}
\label{Fig: Eye surface skematic}
\end{figure}
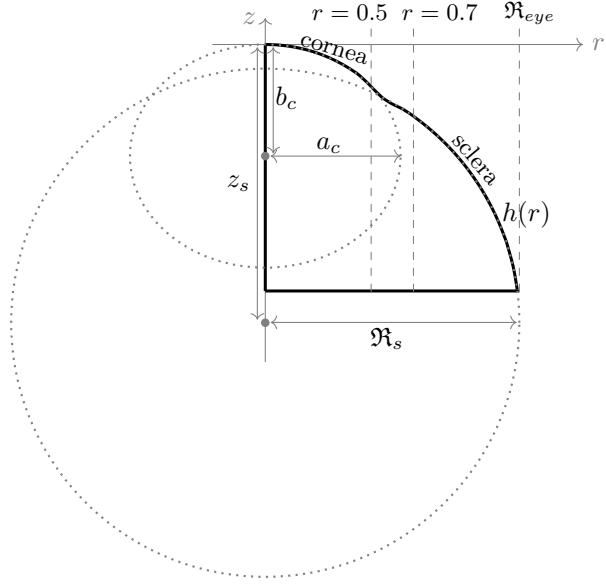

To describe the ocular surface $\Gamma_{out}$ via the function $h(r)$, we divide the ocular surface into the three \rev{segments}: (i) the corneal surface \rev{segment}, (ii) the limbal \rev{surface segment}, and (iii) the scleral \rev{surface segment}, as shown in Figure~\ref{fig:contact}. Here, we first describe how we construct the corneal \rev{segment}, then the scleral  \rev{segment}, and finally the limbal \rev{segment}. A schematic of the construction of the ocular surface is shown in Figure~\ref{Fig: Eye surface skematic}.

The cornea is the clear, dome-shaped window of your eye that focuses light into your eyes.  We characterize the anterior corneal surface as an ellipse given by the following:
\begin{equation}
   z = -b_c + \sqrt{b_c^2 - \left(\frac{b_cr}{a_c}\right)^2} \ \mathrm{cm}, \quad 0 \leq r < 0.5 \ \mathrm{cm},
\end{equation}
centered at $(0,-b_c)$, where $a_c$ and $b_c$ are the axis of the ellipse. The axis of the cornea ellipse can be expressed in terms of the radius of curvature of the ellipse at $r=0$, denoted by $K_c$, and the eccentricity of the ellipse, denoted by $e$, as follows: $a_c^2=K_cb_c$ and $b_c=(1-e^2)K_c$.
The values of $K_c$ and $e$ are chosen to match biometric data present in Hall et al~\cite{Halletal2011}.  
Specifically, $K_c$ is set to the simulated keratometry readings. Recall that keratometry is the measurement of the corneal curvature. Then, the square of the eccentricity $e^2$ is set to the reported corneal shape factor values. Different values of $K_c$ and $e$ are considered to model a flat, average, or steep cornea \rev{surface} shape, as reported in Table~\ref{tab:OcularParameter}. The limit of the corneal region of $0.5$~cm is consistent with the limbal radius of $0.55$~cm reported by Missel \cite{Missel2012} and the mean horizontal visible iris radius of $0.593$~cm reported by Hall et al \cite{Halletal2011}.

The sclera is the white, outer layer of your eye.  We characterize the scleral anterior surface as a circle, centered at $(0,z_s)$, with radius $\mathfrak{R}_s$,
\begin{equation}
    z = z_s + \sqrt{\mathfrak{R}_s^2-r^2} \ \mathrm{cm}, \quad 0.7 \ \mathrm{cm} < r \leq \mathfrak{R}_{eye} ,
\end{equation}
where $\mathfrak{R}_s$ is the radius of the sclera. Different values of $\mathfrak{R}_s$ are considered to model an average or a flat sclera, as reported in Table~\ref{tab:OcularParameter}.  
In the case of an average scleral shape, we choose  $\mathfrak{R}_s=1.2$~cm in agreement with \cite{Missel2012,funkenbusch1996conformity}.
Using the biometric data presented in Hall et al \cite{Halletal2011}, the value of $z_s$ is chosen to match reported experimental values of the ocular sagittal height of a chord at 15 mm. 
We set the ocular domain up to $\mathfrak{R}_{eye}=1.2$~cm. Note that $\mathfrak{R}_{eye}=\mathfrak{R}_s$ when considering an average scleral shape, while $\mathfrak{R}_{eye}<\mathfrak{R}_s$ when considering a flat scleral shape. In any case, $\mathfrak{R}_{eye}$ is always bigger than the undeformed radius of the lens $\mathfrak{R}_{lens}$, reported in Table~\ref{tab:lens_geometry}; thus, we always consider an ocular domain big enough to describe the lens effect on the eye. 
The start of the scleral \rev{segment} (i.e., $r=0.7$~cm) is larger than the mean radius of the corneoscleral junction, $0.656$~cm, reported in \cite{Halletal2011}, to account for the limbal \rev{segment} between the cornea and the sclera.

The limbal \rev{segment} is the transition region between the transparent cornea \rev{surface} and the opaque sclera \rev{surface}. We model this region as a polynomial of the following form: 
\begin{eqnarray}
    z &=& (r-0.7)^3(l_1 + l_2(r-0.5) + l_3(r-0.5)^2) \nonumber\\ 
    &+&(r-0.5)^3(l_4 + l_5(r-0.7)+l_6(r-0.7)^2) \ \mathrm{cm}, \quad 0.5 \ \mathrm{cm} \leq r \leq 0.7 \ \mathrm{cm},
\end{eqnarray}
where the parameters $l_i$, $i=1,...,6$, are chosen to ensure that the ocular surface is a function with continuous first and second derivatives. 

In summary, the ocular surface \rev{$\Gamma_{out}$} is a piecewise function given by the following:
\begin{equation}\label{eqn:ocular_surface}
    h(r) = \begin{cases}
    -b_{c} + \sqrt{b_{c}^2 - \left(\dfrac{b_{c}r}{a_{c}}\right)^2}, \quad & 0 \leq r < 0.5 \ \mathrm{cm} \\
     (r-0.7)^3(l_1 + l_2(r-0.5) + l_3(r-0.5)^2) \\ 
    +(r-0.5)^3(l_4 + l_5(r-0.7)+l_6(r-0.7)^2), \quad &0.5 \ \mathrm{cm} \leq r \leq 0.7 \ \mathrm{cm} \\
    z_s + \sqrt{\mathfrak{R}_s^2-r^2}, \quad &0.7 \ \mathrm{cm} < r \leq \mathfrak{R}_{eye}
    \end{cases}.
\end{equation}
Table~\ref{tab:OcularParameter} reports the parameter values for the four ocular surfaces considered and shown in Figure~\ref{fig:eye_geometry}.
\begin{table}[h]
	\footnotesize
    \centering
      \caption{Ocular surface parameter values with corresponding references.}
    \begin{tabular}{lllll}
    \hline
    & Flat cornea \& & Average cornea \&  & Steep cornea \& & Average cornea \&\\
     & Average sclera & Average sclera  & Average sclera  & Flat sclera\\   
    \hline
    $K_c$ [cm] & 0.785 \cite{Halletal2011} & 0.775 \cite{Halletal2011} & 0.765 \cite{Halletal2011} & 0.775 \cite{Halletal2011} \\
    $e^2$ & 0.430 \cite{Halletal2011} & 0.320 & 0.210 \cite{Halletal2011} & 0.320\\
    $a_c$ [cm] & 0.592 & 0.639 & 0.680 & 0.639 \\
    $b_c$ [cm] & 0.447 & 0.527 & 0.604 & 0.527 \\
    $\mathfrak{R}_s$ [cm] & 1.200 \cite{Missel2012,funkenbusch1996conformity} & 1.200 \cite{Missel2012,funkenbusch1996conformity} & 1.200 \cite{Missel2012,funkenbusch1996conformity}& 4.310 \cite{Halletal2011} \\
    $z_s$ [cm] & -1.314 & -1.314 & -1.314  & -4.621 \\
    $\mathfrak{R}_{eye}$ [cm] & 1.200 \cite{Missel2012} & 1.200 \cite{Missel2012} & 1.200 \cite{Missel2012} & 1.200 \cite{Missel2012}\\
    $l_1$ & 25.902 & 24.847 & 24.349 & 24.847  \\
    $l_2$ & 536.827 & 502.197 & 485.789 & 502.197 \\
    $l_3$ & 6900.432 & 6099.372 & 5760.904 & 6099.372 \\
    $l_4$ & -42.384  & -42.384 & -42.384  & -46.058 \\
    $l_5$ & 545.984 & 545.984 & 545.984 & 670.302 \\
    $l_6$ & -5042.020 & -5042.020 & -5042.020 & -6614.444  \\
    \hline
    \end{tabular}
  
    \label{tab:OcularParameter}
\end{table}

\section{Eye material parameters}

\subsection{Constant young's modulus}
\label{sec:supp_eye_para}
Here is a summary of the resources that we used to estimate the range of the constant Young's modulus of the eye, reported in Table~\ref{tab:material}.
Bryant and McDonnell performed membrane inflation tests on twelve human cornea which measured displacements with a fiber optic displacement probe and estimated material properties by fitting a linear isotropic model to the measured displacements \cite{BryantandMcDonnell1996}. The reported mean Young's modulus was 0.83 MPa with a standard derivation of 0.22 MPa (coefficient of variation is 0.26)~\cite{BryantandMcDonnell1996}. Orssengo and Pye derived a theoretical relationship between  the true intraocular pressure and the Young's modulus for the cornea for an average human cornea~\cite{orssengoandpye1999}. The true intraocular pressure is the pressure in the eye (measured with a Goldmann’s applanation tonometry) divided by a correction factor that takes the dimension, curvature, and material properties of the cornea into account. For a reasonable range of true intraocular pressures, they found that  the Young's modulus can range between $0.2$~MPa (low intraocular pressure) and $1.0$~MPa (high intraocular pressure)~\cite{BryantandMcDonnell1996}.  Woo et al.~found that the sclera can be five stiffer than the cornea at a low IOP~\cite{Wooetal1972}.  

\subsection{Spatially\rev{-d}ependent young's modulus}
\label{sec:supp_YM}
\begin{figure}[htb]
    \centering
    \includegraphics[width=0.45\linewidth]{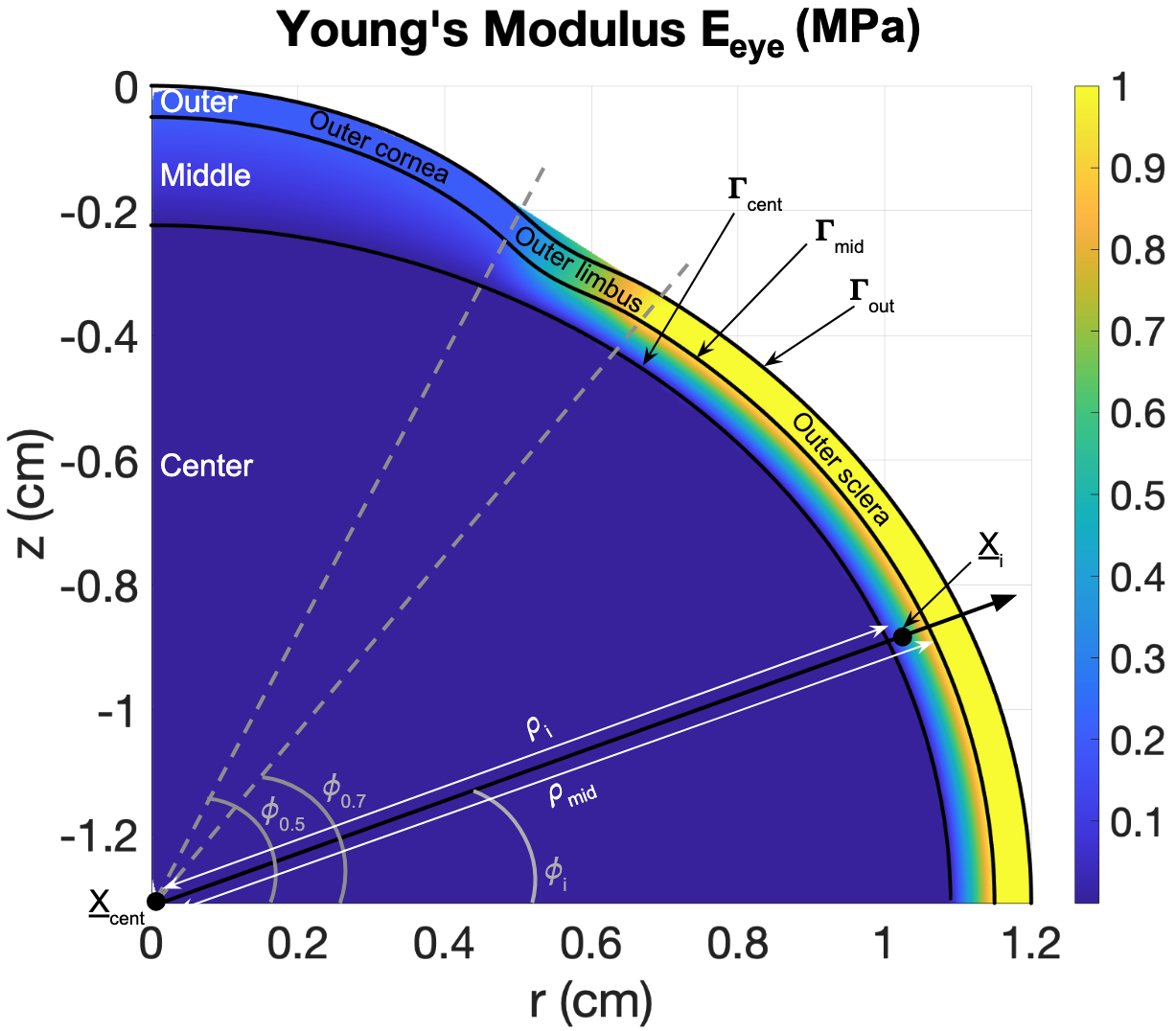}
    \caption{Schematic of the derivation of the spatially\rev{-}dependent Young's modulus values.}
    \label{fig:varied_ym}
\end{figure}

Now, we consider the case in which the Young's modulus of the ocular tissue varies with space. To begin, we partition the ocular domain into three regions, as shown in Figure~\ref{fig:varied_ym}: (i) the outer, (ii) middle, and (iii) center regions. We consider 
an average-shaped eye (Column 2 of Table~\ref{tab:OcularParameter}). 
The outer region mimics the external layer of tissue of the eye. We partition the outer region into three subregions as before: (i) the corneal, (ii) the limbal, and (iii) the scleral regions. The center region represents the inside of the filled with the vitreous humor. The middle region contains anatomical structures such as the aqueous chamber, iris, eye lens, retina, and choroid.

Recall that $\Gamma_{out}$ denotes the outer boundary of the eye. Let $\Gamma_{mid}$ denote the boundary between the outer and middle regions, and let $\Gamma_{cent}$ denote the boundary between the middle and center regions\rev{, see Figure~\ref{fig:varied_ym}}. For the average shaped eye, we assume that the outer ocular tissue is $0.05$~cm thick~\cite{Pandolfi2006}; therefore we define $\Gamma_{mid}$ such that it is a normal distance of $d_{mid} = 0.05$~cm away from $\Gamma_{out}$ toward the inside of the eye.
We assume the choroid and retina tissues to be $0.06$~cm thick at $r = \mathfrak{R}_{eye}$~\cite{Myers2015, Entezari2018}.   
Thus, we model the center region as a circle centered at $(0,z_s)$ with a radius $d_{cent} = \mathfrak{R}_{eye} - 0.05 - 0.06 = 1.09$ cm. 

Let $\vector{X}_{cent} = (0,z_s)$. Each point $\vector{X}_i = (r_i,z_i)$ in the domain $\Omega$ can be alternatively described by two values: (i) its Euclidean distance from $\vector{X}_{cent}$, $\rho_i  = ||\vector{X}_i-\vector{X}_{cent}||_2$; and (ii) the angle between the line $z = z_s$ and the ray connecting $\vector{X}_{cent}$ and $\vector{X}_i$, $\phi_i = \arctan\left(\frac{z_i-z_s}{r_i}\right)$. When describing points in this way, we will denote them as $\vector{\xi}_i = (\rho_i,\phi_i).$ Note that for all the points of $\Gamma_{cent}$, we have that $\rho= d_{cent}$, and for all the points on $\Gamma_{mid}$, we have that $\rho= \rho_{mid}(\phi)$, as shown in Figure~\ref{fig:varied_ym}.

To determine to which region each point $\vector{X}_i$ belongs to and to determine its Young's modulus, we first compute $\rho_i$ and $\phi_i$ to construct $\vector{\xi}_i$. If $0\leq \rho_i < d_{cent}$, then the point $\vector{X}_i$ lies in the center region and has a Young's modulus of $E_{cent}$, which is a constant value throughout the entire center region. 
If $\rho_i \geq \rho_{mid}(\phi_i)$, then the point $\vector{X}_i$ lies in the outer region. The Young's modulus value of a point in the outer region $E_{out}$ depends on whether it lies in the corneal, limbal, or scleral region, thus $E_{out}=E_{out}(\phi_i)$. Points in the cornea have a constant Young's modulus value of $E_{cornea}$, and points in the sclera have a constant value of $E_{sclera} = 5E_{cornea}$ \cite{Wooetal1972}. The Young's modulus value of a point in the limbal region is determined by a linear interpolation in $\phi$ between $E_{cornea}$ and $E_{sclera}$.
If $d_{cent}\leq \rho_i < \rho_{mid}(\phi_i)$, then the point $\vector{X}_i$ lies in the middle region. The Young's modulus of a point in the middle region is determined by a linear interpolation in $\rho$ between $E_{cent}$ and $E_{out}(\phi_i)$. In summary, the Young's modulus value $E_{eye}\left(\vector{\xi}_i\right)$ for a given point $\vector{\xi}_i$ in $\Omega$ is determined as follows: \begin{equation}
    E_{eye}\left(\vector{\xi}_i\right) = E_{eye}(\rho_i,\phi_i) = \begin{cases}
        E_{cent}, & 0\leq \rho_i < d_{cent} \\
        E_{mid}(\rho_i,\phi_i), & d_{cent} \leq \rho_i < \rho_{mid}(\phi_i) \\
        E_{out}(\phi_i), & \rho_{mid}(\phi_i) \leq \rho_i
    \end{cases},
\end{equation}
where \begin{equation}
    E_{mid}(\rho_i,\phi_i) = \left(\frac{E_{cent}-E_{out}(\phi_i)}{d_{cent}-\rho_{mid}(\phi_i)}\right)(\rho_i-d_{cent})+E_{cent},
    \end{equation}
    and \begin{equation} E_{out}(\phi_i) = \begin{cases}
        E_{cornea}, & \phi_i > \phi_{0.5} \\
        \left(\dfrac{E_{cornea}-E_{sclera}}{\phi_{0.5}-\phi_{0.7}}\right)(\phi_i-\phi_{0.5})+E_{cornea}, & \phi_{0.5} \leq \phi_i \leq \phi_{0.7} \\
        E_{sclera}, & \phi_i < \phi_{0.7}
    \end{cases},
\end{equation}
where $\phi_{0.5}$ and $\phi_{0.7}$ are the angles used to describe the points on $\Gamma_{out}$ where $r = 0.5$~cm and $r = 0.7$~cm (the determine the range of the limbus), respectively, \rev{see Figure~\ref{fig:varied_ym}}.

For the results presented here, the Young's modulus in the cornea and sclera are set to $E_{cornea}~=~0.2$~MPa and $E_{sclera} = 5E_{cornea} = 1.0$~MPa, respectively. We define the material properties of the center region to best imitate the vitreous humor, which is often modeled as a viscoelastic hydrogel \cite{Tram2021}. Therefore, we set the Young's modulus in this region to be $E_{center} = 1.17\times10^{-6}$ MPa \cite{Tram2021}. 

\rev{
\section{Mesh grid refinement study}
\label{sec:more_grid_ref}
Since the largest deformations occur along the ocular surface $\Gamma_{out}$, we compare the magnitude of the displacement vector (2 norm) evaluated on $\Gamma_{out}$ for the different meshes used in the grid refinement study. The results of this comparison is shown in Figure~\ref{fig:grid_refinement2}.
\begin{figure}[h!]
\begin{center}
\begin{tikzpicture}
    \node at (0,0) {\includegraphics[width=.5\linewidth]{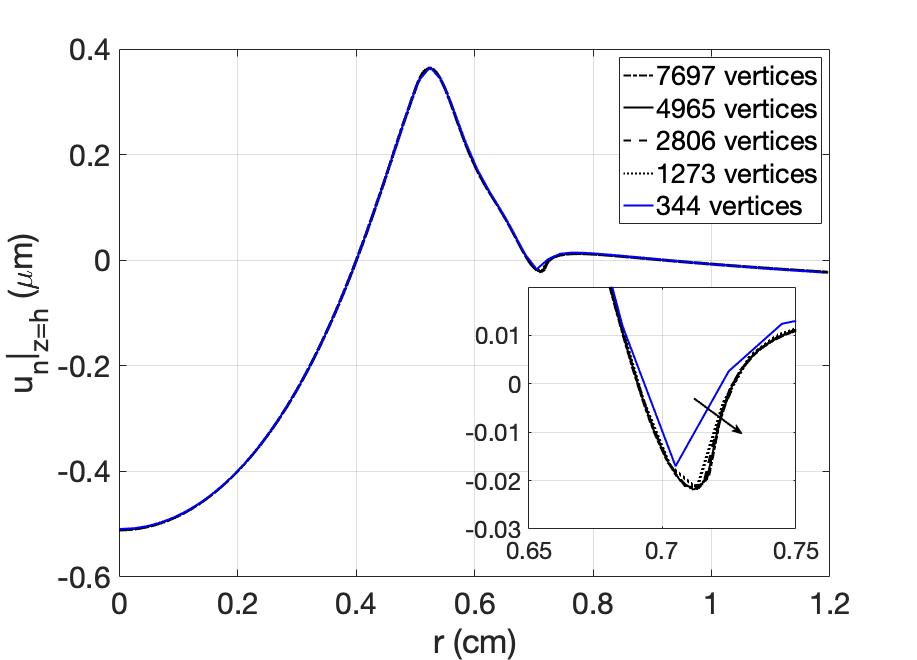}};
\end{tikzpicture}
    \caption{\rev{The magnitude of the ocular displacements of a homogeneous average-shaped eye when an average-shaped contact lens of $100$~$\mu$m constant thickness is inserted, for the different meshes considered in the grid refinement study in Section~\ref{sec:numerical_method}.}}
    \label{fig:grid_refinement2}
\end{center}    
\end{figure}
}

\end{document}